\documentclass[pdftex,12pt]{article}
\usepackage[dvips]{graphicx,color}
\usepackage[pdftex,colorlinks]{hyperref}
\usepackage{amssymb}
\usepackage{latexsym}
\usepackage[all,cmtip]{xy}

\newcommand{\myparagraph}[1]{\paragraph{\bf #1}}
\newcommand{\bit}{\begin{itemize}}
\newcommand{\eit}{\end{itemize}}
\newcommand{\ben}{\begin{enumerate}}
\newcommand{\een}{\end{enumerate}}

\newcommand{\bde}{\begin{definition}}
\newcommand{\ede}{\end{definition}}
\newcommand{\bre}{\begin{remark}}
\newcommand{\ere}{\end{remark}}
\newcommand{\ble}{\begin{lemma}}
\newcommand{\ele}{\end{lemma}}
\newcommand{\bpr}{\begin{proposition}}
\newcommand{\epr}{\end{proposition}}
\newcommand{\bth}{\begin{theorem}}
\renewcommand{\eth}{\end{theorem}}
\newcommand{\bco}{\begin{corollary}}
\newcommand{\eco}{\end{corollary}}

\newcommand{\bmth}{\begin{mainthm}}
\newcommand{\emth}{\end{mainthm}}
\newcommand{\bmpr}{\begin{mainprop}}
\newcommand{\empr}{\end{mainprop}}

\newcommand{\brth}{\begin{refthm}}
\newcommand{\erth}{\end{refthm}}
\newcommand{\brpr}{\begin{refprop}}
\newcommand{\erpr}{\end{refprop}}

\newcommand{\bcon}{\begin{conj}}
\newcommand{\econ}{\end{conj}}

\newcommand{\bex}{\begin{example}}
\newcommand{\eex}{\end{example}}

\newcommand{\barr}{\begin{array}}
\newcommand{\earr}{\end{array}}
\newcommand{\btab}{\begin{tabular}}
\newcommand{\etab}{\end{tabular}}
\newcommand{\beq}{\begin{equation}}
\newcommand{\eeq}{\end{equation}}
\newcommand{\bea}{\begin{eqnarray*}}
\newcommand{\eea}{\end{eqnarray*}}
\newcommand{\bce}{\begin{center}}
\newcommand{\ece}{\end{center}}
\newcommand{\bpi}{\begin{picture}}
\newcommand{\epi}{\end{picture}}
\newcommand{\bfi}{\begin{figure} \begin{center}}
\newcommand{\efi}{\end{center} \end{figure}}

\newcommand{\bsl}{\begin{slide}{}}
\newcommand{\esl}{\end{slide}}

\newcommand{\pf}{\noindent{\bf Proof}\hspace{7pt}}

\newcommand{\hso}[1]{\hspace{-1pt}}

\newcommand{\into}{\hookrightarrow}
\newcommand{\onto}{\twoheadrightarrow}

\newcommand{\sbe}{\subseteq}

\def\<{\langle}
\def\>{\rangle}

\newcommand{\la}{\lambda}

\newcommand{\De}{\Delta}

\newcommand{\La}{\Lambda}
\newcommand{\Om}{\Omega}

\newcommand{\bC}{{\bf C}}

\newcommand{\bK}{{\bf K}}

\newcommand{\bV}{{\bf V}}

\newcommand{\cA}{{\cal A}}

\newcommand{\cB}{{\cal B}}

\newcommand{\cC}{{\cal C}}

\newcommand{\cE}{{\cal E}}

\newcommand{\cG}{{\cal G}}

\newcommand{\cL}{{\cal L}}

\newcommand{\cP}{{\cal P}}

\newcommand{\cR}{{\cal R}}

\newcommand{\cS}{{\cal S}}

\newcommand{\cV}{{\cal V}}

\newcommand{\eI}{I}
\newcommand{\eJ}{J}
\newcommand{\eK}{K}
\newcommand{\eL}{L}

\newcommand{\Aut}{\mathop{\rm Aut}\nolimits}

\renewcommand{\bar}{\overline}

\newcommand{\diag}{\mathop{\rm diag}\nolimits}

\newcommand{\End}{\mathop{\rm End}\nolimits}

\newcommand{\Hom}{\mathop{\rm Hom}\nolimits}
\newcommand{\id}{\mathop{\rm id}\nolimits}

\newcommand{\inc}{\mathop{\rm inc}\nolimits}

\newcommand{\opp}{\mathbin{\rm opp}}
\newcommand{\near}{\mathbin{\rm near}}

\newcommand{\Rad}{\mathop{\rm Rad}\nolimits}

\newcommand{\sign}{\mathop{\rm sign}\nolimits}

\newcommand{\Sym}{\mathop{\rm Sym}\nolimits}

\newcommand{\trace}{\mathop{\rm trace}\nolimits}

\newcommand{\supp}{{\mathop{\rm supp}\nolimits}}

\newcommand{\mn}{\medskip\noindent}

\def\flexbox#1{\mathchoice{\mbox{#1}}{\mbox{#1}}{\mbox{\scriptsize #1}}%
{\mbox{\tiny #1}}}
\newcommand{\GL}{\mathop{\rm GL}}

\newcommand{\SL}{\mathop{\flexbox{\rm SL}}}
\newcommand{\PSL}{\mathop{\flexbox{\rm PSL}}}
\newcommand{\SO}{\mathop{\flexbox{\rm SO}}}
\newcommand{\POm}{\mathop{\rm P\Omega}}

\newcommand{\PSO}{\mathop{\flexbox{\rm PSO}}}

\newcommand{\Spin}{\mathop{\flexbox{\rm Spin}}}
\newcommand{\Sp}{\mathop{\rm Sp}}
\newcommand{\PSp}{\mathop{\flexbox{\rm PSp}}}
\newcommand{\SU}{\mathop{\flexbox{\rm SU}}}

\newcommand{\PG}{\mathop{\flexbox{\rm PG}}}

\newcommand{\wA}{{{}^2\hspace{-.2em}A}}
\newcommand{\wD}{{{}^2\hspace{-.2em}D}}
\newcommand{\wF}{{{}^2\hspace{-.2em}F}}

\newcommand{\Char}{\mathop{\flexbox{\rm Char}}}

\newcommand{\CC}{{\mathbb C}}
\newcommand{\DD}{{\mathbb D}}

\newcommand{\FF}{{\mathbb F}}

\newcommand{\NN}{{\mathbb N}}
\newcommand{\PP}{{\mathbb P}}
\newcommand{\VV}{{\mathbb V}}
\newcommand{\ZZ}{{\mathbb Z}}

\renewcommand{\inc}{\star}
\newcommand{\typ}{\mathop{\rm typ}}%

\newcommand{\dfn}{\em}
\newcommand{\Res}{\mathop{\rm Res}\nolimits}
\newcommand{\proj}{\mathop{\rm proj}\nolimits}%

\newcommand{\gr}{{\mathop{\rm gr}}}%

\newcommand{\sfC}{{\mathsf{C}}}
\newcommand{\sfD}{{\mathsf{D}}}
\newcommand{\sfE}{{\mathsf{E}}}
\newcommand{\sM}{{\mathsf{M}}}

\newcommand{\after}{{\mathbin{\circ}}}

\newcommand{\abs}{\widetilde}
\newcommand{\rabs}{\overline}
\newcommand{\rmin}{\underline}

\newcommand{\Far}{\mathop{\rm Far}\nolimits}

\newcommand{\Stab}{\mathop{\rm Stab}\nolimits}

\newcommand{\AI}[1]{\item[\rm{(#1)}]}

\newcommand{\mmod}[1]{\mathbin{({\rm mod}\ #1)}}

\newcommand{\mfg}{\mathfrak{g}}

\newcommand{\mfh}{\mathfrak{h}}

\newcommand{\mfsl}{\mathfrak{sl}}
\newcommand{\mfgl}{\mathfrak{gl}}

\newcommand{\mfI}{\mathfrak{I}}
\newcommand{\mfL}{\mathfrak{L}}
\newcommand{\mfU}{\mathfrak{U}}

\def\ad{\mathop{\flexbox{\rm ad}}}

\newtheorem{theorem}{Theorem}[section]
\newtheorem{proposition}[theorem]{Proposition}
\newtheorem{corollary}[theorem]{Corollary}
\newtheorem{lemma}[theorem]{Lemma}
\newtheorem{definition}[theorem]{Definition}

\newtheorem{example}[theorem]{Example}
\newtheorem{remark}[theorem]{Remark}

\newtheorem{mainthm}{Theorem}

\newtheorem{mainprop}[mainthm]{Proposition}

\newtheorem{refthm}{Theorem}

\newtheorem{refprop}[refthm]{Proposition}

\newcommand{\qed}{\hfill $\square$}

\newcommand{\Problems}{\noindent{\bf Problems}\quad}

\newcommand{\institute}[1]{\ \newline{#1}}
\newcommand{\email}[1]{\texttt{#1}}
\newcommand{\keywords}[1]{\ \newline{Keywords:} #1}
\newcommand{\subclass}[1]{\ \newline{Mathematics Subject Classification (2010)} #1}

\newcounter{romanlistctr}
{\end{list}}%

\newcommand{\goodsofar}{}

\title{Highest weight modules and polarized embeddings of shadow spaces}
\author{Rieuwert J. Blok\thanks{This research was completed in part while visiting the University of Siena on a grant from the Gruppo Nazionale per le Strutture Algebriche, Geometriche e le loro Applicazioni  in the summers of 2007 and 2008.}\\
\institute{}Department of Mathematics and Statistics,\\
Bowling Green State University,\\
Bowling Green, OH 43403,\\
U.S.A.\\
\email{blokr@member.ams.org}\\
}

\begin{document}
\maketitle
\begin{abstract}
The present paper was inspired by the work on polarized embeddings by Cardinali, De Bruyn, and Pasini~\cite{CaDePa2007}, although some of our results in it date back to 1999.
They study polarized embeddings of certain dual polar spaces, and identify the minimal polarized embeddings for several such geometries.
We extend some of their results to arbitrary shadow spaces of spherical buildings, and make a connection to work of Burgoyne, Wong, Verma, and Humphreys on highest weight representations for Chevalley groups.

Let $\De$ be a spherical Moufang building with diagram $\sM$ over some index set $\eI$, whose strongly transitive automorphism group is a Chevalley group $G(\FF)$ over the field $\FF$.
For any non-empty set $\eK\sbe \eI$ let $\Gamma$ be the $\eK$-shadow space of $\De$.
Extending the notion in~\cite{CaDePa2007,ThVa2006} to this situation, we say that an embedding of $\Gamma$ is polarized if it induces all singular hyperplanes. Here a singular hyperplane is the collection of points of $\Gamma$ not opposite to a point of the dual geometry $\Gamma^*$, which is the shadow geometry of type $\opp_\eI(\eK)$
 opposite to $K$.
We prove a number of results on polarized embeddings, among others the existence of (relatively) minimal
 polarized embeddings.

We assume that $G(\FF)$ is untwisted.
In that case, the point-line geometry $\Gamma$ has an embedding $e_\eK$ into the Weyl module $V(\lambda_\eK)_\FF^0$ of highest weight $\lambda_\eK=\sum_{k\in \eK}\lambda_k$.
We show that this embedding is polarized in the sense described above.
We then prove that the minimal polarized embedding relative to $e_\eK$ exists and equals the unique irreducible $G(\FF)$-module $L(\lambda_\eK)$ of highest weight $\lambda_\eK$.
More precisely we show that the polar radical of $e_\eK$ (the intersection of all singular hyperplanes) coincides with the radical of the contravariant bilinear form considered by Wong to obtain the irreducible (restricted) representations of $G(\FF)$ in positive characteristic.

This viewpoint allows us to ``recognize'' the irreducible $G(\FF)$-modules of highest weight
 $\lambda_\eK$ geometrically as minimal polarized embeddings of the appropriate shadow space.
\subclass{
    Primary 51E24;
    Secondary 20G05, 20G15}
\end{abstract}

\newpage
\pagestyle{plain}

\section{Introduction and preliminaries}\label{sec:introduction}

\subsection{Basic definitions}\label{subsection:basic definitions}
\myparagraph{Vector spaces, group actions and modules}
Throughout this paper $\FF$ shall denote a field.
By a field we shall mean a commutative division ring.
Unless otherwise specified, vector spaces will be left vector spaces.
Actions of groups, or algebras on a vector space shall therefore also be left actions, unless otherwise specified.
\myparagraph{Point-line geometries, hyperplanes and embeddings}
We assume the reader is familiar with the concept of a {\em point-line geometry}
$\Gamma = (\cP, \cL)$ (also called a partial linear rank two
incidence geometry), see e.g.~\cite{Bu1995,Pa1994}.
In a partial linear space we can and shall often identify each line with the set of points incident to it.
By a {\it subspace} of $\Gamma$ we mean a subset $S\sbe \cP$ such that if $l \in \cL$ and $l \cap S$ contains at
least  two points, then $l \subset S$.
Clearly the intersection of subspaces is a subspace and consequently it is natural to define the
subspace generated by a subset $X$ of $\cP, \langle X \rangle_\Gamma,$ to be the intersection of
all subspaces of $\Gamma$ that contain $X$.
A {\dfn hyperplane} of $\Gamma$ is a proper subspace meeting every line.
A set of points is called {\em connected} if its collinearity graph is connected.
We recall the following.
\ble\label{lem:maximal hyperplanes}
If $H$ is a hyperplane such that $\Gamma-H$ is connected, then $H$ is a maximal subspace of $\Gamma$.
\ele
The {\dfn projective point-line geometry} of a vector space $V$ is the point-line
 geometry $\PP(V)=(\cP(V),\cL(V))$ whose points and lines are the $1$-spaces and $2$-spaces of $V$ with incidence given by symmetrized inclusion.

A {\dfn full projective embedding} (or simply {\em embedding}) of $\Gamma$ is a pair $(e,V)$, where
$V$ is a vector space and $e\colon\cP\into\cP(V)$ is an injective map such that
\begin{itemize}
\AI{E1} $\langle e(\cP)\rangle_V=V$, and
\AI{E2} $e$ maps every line of $\Gamma$ onto a line of $\PP(V)$.
\end{itemize}
The {\dfn dimension} of $(e,V)$ is $\dim(V)$. In the literature, this is sometimes called the {\dfn vector dimension} of the embedding to distinguish it from its projective dimension.

The collection $\cE(\Gamma)$ of all full projective embeddings of $\Gamma$ over a division ring (or field) $\FF$, is a category where a morphism between embeddings $(e_1,V_1)$ and
 $(e_2,V_2)$ is an $\FF$-semilinear map $\tau\colon V_1\to V_2$ such that
  $e_2=\tau\after e_1$.
We sometimes indicate this by writing $e_1\ge e_2$.
We have the usual notions of mono-, epi-, and isomorphisms.
An embedding of $\Gamma$ is called {\dfn absolutely universal} or {\dfn absolute} if it is a source in $\cE(\Gamma)$; if it exists, we denote it by $(\abs{e},\abs{V})$.
A source relative to $(e,V)$ always exists by a result due to Ronan~\cite{Ro1987}; it is called the embedding {\em universal relative to $(e,V)$} and will be denoted $(\rabs{e},\rabs{V})$.

An embedding $(e,V)$ of $\Gamma$ is called {\dfn minimal} if it is a sink in $\cE(\Gamma)$ and {\dfn minimal relative to} $(e,V)$ if it is a sink relative to $(e,V)$; the latter will be denoted $(\rmin{e},\rmin{V})$.

Let $(e,V)$ be a full projective embedding for $\Gamma$.
For a point set $X\sbe \cP$, let $\langle X\rangle_e=\langle e(X)\rangle_V$.
The following is well-known and elementary.

\ble\label{lem:induced hyperplanes}
\
\begin{itemize}
\AI{a} If $U$ is a hyperplane of $V$, then $H=e^{-1}(U\cap e(\cP))$
 is a hyperplane of $\Gamma$.
\AI{b} If $H$ is a maximal hyperplane of $\Gamma$, then $\langle H\rangle_e$ either
 equals $V$ or it is a hyperplane of $V$;
  in the latter case $e(H)=e(\cP)\cap \langle H\rangle_e$.
\end{itemize}
\ele
The hyperplane $H$ in (a) is said to be {\dfn induced} by $U$ in  $(e,V)$.

\myparagraph{Buildings and Chevalley groups}
Every geometry we shall study in this paper is derived from a spherical building $\De$.
The building $\De$ has spherical diagram $\sM$ over the index set $I=\{1,2,\ldots,n\}$.
We shall label $\sM$ as in \cite{Bo1975}.
We shall think of $\De$ as a chamber system, also denoted $\De$, with a distance function $\delta\colon\De\times \De\to W$, where $(W,\{r_i\}_{i\in I})$ is a Coxeter system of type $\sM$; we shall use the terminology from~\cite{Ro1989a,We2003}.
Thus, for $i\in I$, we say that two chambers $c$ and $d$ are $i$-adjacent, and write $c\sim_i d$, if $\delta(c,d)=r_i$. For a subset $\eJ\sbe\eI$, the $\eI-\eJ$ residue on $c$ is the collection of chambers
 $\{d\in \Delta\mid \delta(c,d)\in W_{\eI-\eJ}\}$.

Sometimes it will be convenient to talk about $\De$ in terms of the associated incidence geometry $\cG=\cG(\De)$.
This is a diagram geometry with a set of elements $\cE=\cE(\De)$, an incidence relation $\inc$ and a type function $\typ\colon\cE\to I$; 
 for $\cG(\De)$ we shall use the terminology from~\cite[Chap.1]{Bu1995}.
Thus, for $\eJ\sbe \eI$, a flag of type $\eJ$ (or $\eJ$-flag) in $\cG$ is a collection $F=\{f_j\}_{j\in \eJ}$ of pairwise incident elements and $\Res_\De(F)=\{e\in \cE-F\mid e\inc f_j\ \forall j\in \eJ\}$ is the corresponding residue of type $\eI-\eJ$ (or $\eI-\eJ$-residue).
 
For buildings of finite rank these two viewpoints are equivalent (see e.g.\cite{Ti1981} or \cite{Pa1994}); namely
 flags of type $\eJ\sbe \eI$ in the diagram geometry $\cG$ correspond to residues of type $\eI-\eJ$ in the chamber system $\De$ and incidence in the diagram geometry corresponds to inclusion in the chamber system.

Whenever a spherical building $\De$ is Moufang, it has a strongly transitive automorphism group $G$ from which $\De$ can be recovered via a BN-pair in the manner described e.g.\ in~\cite[Ch.~5]{Ro1989a} or ~\cite[Ch.~11]{We2003}.  Namely,
pick an apartment $\Sigma$ and chamber $c\in\Sigma$, then $B=\Stab_G(c)$ and $N=\Stab_G(\Sigma)$
 form a BN-pair for $G$. Conversely, given a BN-pair $(B,N)$, one can construct $\De$ by setting $\De=G/B$
 and defining the $W$-valued distance function using the Bruhat decomposition of $G$.
For the spherical buildings under consideration in the present paper we take $G$ to be a Chevalley group of rank $n\ge 2$. In Sections~\ref{section:polar hyperplanes}~and~\ref{section:covers and quotients}, $G$ can be of any twisted or untwisted type; in Sections~\ref{section:Weyl embedding},~\ref{section:minimal polarized embedding},~\ref{section:minuscule weights}, and~\ref{section:grassmannians}, we require $G$ to be untwisted. We also assume $G$ to be the universal Chevalley group with diagram $\sM$ over $\FF$, denoted $\bar{\sM}(\FF)$, although this does not affect $\De$. In Section~\ref{section:polar grassmannians} we also include the twisted types $\wA_{2n}$, $\wA_{2n+1}$, and $\wD_{2n+2}$.
Finally, Sections~\ref{section:projective grassmannians}~and~\ref{section:projective flag grassmannians} only concern the $A_n$ case.

The correspondence between their diagrams $\sM$, the (commutative) field of definition $\FF$,  and concrete universal Chevalley groups $G$ is given by the following tables:

\begin{table}[h]
\renewcommand{\arraystretch}{1.2}
$$\begin{array}{cc@{}cc}
\begin{array}{@{}l|l@{}}
\sM & G\\
\hline
A_n & \SL_{n+1}(\FF)\\
B_n & \Spin_{2n+1},(\FF)\\
C_n & \Sp_{2n}(\FF) \\
D_n & \Spin^+_{2n}(\FF) \\
\end{array}
&
\begin{array}{@{}l|l@{}}
\sM & G\\
\hline
E_6 & \bar{E}_6(\FF)       \\
 E_7 & \bar{E}_7(\FF)   \\
 F_4 & \bar{F}_4(\FF)   \\
 G_2 & \bar{G}_2(\FF)   \\
\end{array}
&\hspace{1in}
&
\begin{array}{@{}l|l@{}}
\sM & G\\
\hline
\wA_{2n} & \SU_{2n+1}(\FF)\\
 \wA_{2n-1} &\SU_{2n}(\FF) \\
\wD_{n+1} & \Spin^-_{2n+2}(\FF)\\
&\\
\end{array}\\
&&&\\
\multicolumn{2}{c}{\mbox{Untwisted universal Chevalley groups       }} && \mbox{Some twisted Chevalley groups}\\
\end{array}$$
\renewcommand{\arraystretch}{1}
\caption{Diagrams and Groups}\label{table:diagrams and groups}
\end{table}
\myparagraph{Shadow spaces}
For any subset $\eL\sbe \eI$, we denote by $\cS_\eL(\De)=(\cP_\eL(\De),\cL_\eL(\De))$ the {\dfn $\eL$-shadow space of $\De$} (The term ``$L$-Grassmannians" is also used, e.g.\ in~\cite{Pa1994}).
This is the point-line geometry whose point set $\cP_\eL(\De)$ consists of the {\dfn flags of type $\eL$} of $\cG(\De)$ and whose line set $\cL_\eL(\De)$ consists of the collections of points incident to a flag of {\em cotype} $\{l\}$ for some $l\in \eL$ (We call $l$ the type of that line). 
This is a partial linear space (see e.g.~\cite{Co1976a,Co1995}). If $\eL=\{l\}$ for some $l\in \eI$, this is sometimes called the $l$-Grassmannian of $\De$. In case $\De$ is a building of type $A_n(\FF)$, this is the usual Grassmannian of the vector space $\FF^{n+1}$. These single-node shadow spaces are also called {\em Lie incidence geometries}, see e.g.~\cite{Co1976}. Note that if $\eL=\emptyset$, then $\cP_\eL(\De)$ consists of a single point.
At the other extreme, if $\eL=\eI$, then $\cP_\eL(\De)$ is the collection of chambers of $\De$, and $\cL_\eL(\De)$ is the collection of panels of $\De$.
Fixing a non-empty subset $\eK\sbe \eI$, we shall denote $\Gamma=\cS_\eK(\De)=(\cP,\cL)$ (so for this fixed geometry we drop the subscripts $\eK$ and $\De$.

Keeping $\eL$ as above, we define the $\eL$-shadow of an arbitrary set of chambers $X\sbe \De$ to be
$$\cP_{\De,\eL}(X)=\{x\in \cP_\eL(\De)\mid x\cap X\ne \emptyset\}.$$
Here we view $x$ and $X$ as sets of chambers.
We also set
 $$\begin{array}{ll}
 \cS_{\De,\eL}(X)&=(\cP_{\De,\eL}(X),\cL_{\De,\eL}(X))\mbox{, where }\\
 \cL_{\De,\eL}(X)&=\{l\in \cL_\eL(\De)\mid l\cap \cP_{\De,\eL}(X)|\ge 2\}\\
\end{array}$$
Note that $\cS_{\De,\eL}(\De)=\cS_\eL(\De)$.

\mn
We are mostly interested in the case where $X$ is a residue.
Fix $\eJ\sbe \eI$. We let $\sM_\eJ$ denote the subdiagram of $\sM$ induced on the nodes indexed by $\eJ$. Then, 
 $(W_\eJ,\{r_j\}_{j\in\eJ})$ is a Coxeter system of type $\sM_\eJ$ (See e.g.~\cite[Corollary 2.14]{Ro1989a}).
\ble\label{lem:residues are buildings}{\rm(Theorem 3.5 of~\cite{Ro1989a})}
Let $R$ be a $\eJ$-residue of $\De$. Then $R$ is a building of type $\sM_\eJ$; the distance function is given by the restriction $\delta\colon \cC_\De(R)\times \cC_\De(R)\to W_\eJ$.
\ele

\ble\label{lem:restricted shadow spaces} (See also~\cite{Ka2007,Ka2009})
Let $R$ be a $\eJ$-residue of $\De$. 
\begin{itemize}
\AI{a} For any $\eL\sbe \eI$, we have a natural isomorphism $\cS_{\De,\eL}(R) \cong  \cS_{\De,\eL\cap\eJ}(R)$.
\AI{b} If $\eL_1,\eL_2\sbe \eI$ are such that
   $\eL_1\cap\eJ=\eL_2\cap\eJ$, then we have a natural isomorphism
 $\cS_{\De,\eL_1}(R) \cong  \cS_{\De,\eL_2}(R)$.
\AI{c} If $\eL\sbe \eI$, then 
$\cS_{\De,\eL}(R)\cong \cS_{\eL\cap \eJ}(R)$.
\end{itemize}
\ele
\pf
(a) Viewing the elements of $\cP_\eL(\De)$ as sets of chambers, the isomorphism is induced by the map $x\mapsto x\cap R$, for $x\in \cP_{\De,\eL}(R)$.
For more details see the references~\cite{Ka2007,Ka2009}.
Part (b)  is an immediate corollary to part (a).
Part (c) follows from part (b) and the simple observation that 
 $\cS_{\De,\eL\cap \eJ}(R)\cong \cS_{\eL\cap \eJ}(R)$.
\qed

\mn
Note that in part (c) of Lemma~\ref{lem:restricted shadow spaces}, 
 $\cS_{\eL\cap \eJ}(R)$ is a shadow space of the building $R$ as in 
Lemma~\ref{lem:residues are buildings}.
Thus, part (c) allows us to view the $\eL\cap\eJ$-shadow space of $R$ as a subspace of the $\eL$-shadow space of $\De$.
 
\mn
In addition to the references mentioned above, there are a few texts in preparation that deal extensively with shadow spaces:~\cite{BuCo2007,Sh2005}.

\subsection{Main definitions}
In order to present the main results of this paper, we need some new definitions.
\myparagraph{The opposition relation and the dual geometry}
Let $R$ be a residue in $\De$ of type $\eJ\sbe \eI$. 

Conjugation by the longest word $w_\eJ$ of the Coxeter system $(W_\eJ,
\{r_j\}_{j\in \eJ})$ induces the opposition relation $\opp_\eJ$ on the set $\eJ$:
 $i\, \opp_\eJ\, j$ if and only if $r_{i}=r_j^{w_\eJ}$. We set
$$\begin{array}{ll}
\opp_\eJ(\eL)&=\{j\in \eJ\mid j \opp_\eJ l\mbox{ for some }l\in\eL\}.\\
\end{array}$$
In case $\eJ=\eI$, we'll drop it from the notation.

The opposition relation between the chambers of $R$ is given by 
 $x\opp_R y$ if and only if $\delta(x,y)=w_\eJ$.
We extend the opposition relation to residues $S$ and $T$ of $R$ by setting $S\opp_R T$ if every chamber of $S$ is opposite some chamber of $T$ and conversely.

\ble\label{lem:opposite residues}{\rm(see~\cite[Proposition 9.9 and Lemma 9.10]{We2003})}
Let $S$ be an $\eL$-residue of $R$ and let $t\in R$ be a chamber opposite some chamber of $S$.
Let $T$ be the $\opp_\eJ(\eL)$-residue on $t$.
Then $T$ is the unique residue of $R$ opposite to $S$ and containing $t$.
\ele

\mn
Let $\eL\sbe \eJ$ and let $S$ be an $\eL$-residue of $R$. Define
$$\begin{array}{ll}
\opp_R(S)&=\{c\in R\mid c \opp_R d\mbox{ for some chamber }d\in S\}\\
\near_R(S)&= R-\opp_R(S)\\
\end{array}$$
In case $R=\De$, we'll drop it from the notation.
The following is immediate from Lemma~\ref{lem:opposite residues}.
\bco\label{cor:partition of opp}
Let $S$ be an $\eL$-residue of $R$.
Then, $\opp_R(S)$ and $\near_R(S)$ can be partitioned into residues of type $\opp_\eJ(\eL)$.
\eco

\mn
Recall that $\Gamma=(\cP,\cL)$ is the $\eK$-shadow space of $\De$.
The geometry {\em dual to $\Gamma$ in $\De$} will be denoted $\Gamma^*=(\cP^*,\cL^*)$; it is equal to $\cS_{\opp(\eK)}(\De)$.
The {\dfn geometry dual to $\cS_{\eL}(R)$ in $R$} will be denoted $\cS_{\eL}(R)^*$; it is equal to $\cS_{\opp_\eJ(\eL)}(R)$.

We'll call the points and lines of $\Gamma^*$ {\em dual points} and {\em dual lines}.
We shall make use of the fact that $\Gamma^{**}=\Gamma$. All statements made about $\Gamma$ can also be dualized and so we may and shall freely apply results stated for $\Gamma$ also to $\Gamma^*$.

\mn
\myparagraph{Opposite and Far}
We briefly mention a concept related to the opposition relation.
For residues $S$ and $T$ of $R$ we set
 $$S\Far_R T\mbox{ if and only if } s \opp_R t\mbox{ for some chambers }s\in S, t\in T.$$
Using the correspondence $R\leftrightarrow \cG(R)$, this then also defines a relation $\Far_{\cG(R)}$ for flags of $\cG(R)$.

Let $\Far_{\cG(R)}(S)$ be the incidence system of all objects of $\cG(R)$ `far' from $S$ with incidence inherited from $\cG(R)$. It is proved in~\cite{BlBr1998} that this is a transversal geometry with a Buekenhout-Tits diagram whose flag system can be identified with the
 collection of all residues of $R$ intersecting $\opp_R(S)$ (as chamber sets) non-trivially. In particular, the set
  $\opp_R(S)$ is exactly the chamber system of $\Far_{\cG(R)}(S)$.
This means that two flags that are far from $S$ are incident exactly if they share a chamber that belongs to $\opp_R(S)$.
In particular, for any type set $\eL\sbe \eJ$, the subspace 
$\cS_{R,\eL}(\opp_R(S))$ of $\cS_\eL(R)$ coincides with the $\eL$-shadow space defined by $\Far_{\cG(R)}(S)$ in the obvious way.

\myparagraph{Singular hyperplanes and polarized embeddings}
We continue the notation from above.
There are two important geometric structures associated to $\opp_R(S)$ and $\near_R(S)$.
$$
\begin{array}{ll}
\cP_{R,\eL}(\opp_{R}(S))&=\{x\in \cP_\eL(R)\mid x\cap \opp_R(S)\ne \emptyset\}\\
\cP_{R,\eL}(\near_{R}(S))&=\{x\in \cP_\eL(R)\mid x\cap \near_R(S)\ne \emptyset\}\\
\end{array}$$
If $\eL$ is clear  from the context, we shall omit it to unburden the notation. 
These sets are not always disjoint.
However, we have the following simple observation.

\ble\label{lem:H disjoint from opp}
Let $R$ be a $\eJ$-residue of $\De$.
Let  $\eL,\eL^*\sbe \eJ$ and let $S$ be a residue of  type $\eJ-\eL^*$ in $R$. Then $\cP_{R,\eL}(\near_R(S))\cap \cP_{R,\eL}(\opp_R(S))=\emptyset\mbox{  if and only if  } \opp_\eJ(\eL^*)\sbe \eL.$
In particular, the $\eK$-shadows of $\opp_\De({p^*})$ and $\near_\De({p^*})$
 are disjoint subsets of $\cP$ if $p^*\in \cP^*$.
\ele
\pf
We have $\cP_{R,\eL}(\near_R(S))\cap \cP_{R,\eL}(\opp_R(S))=\emptyset$ if and only if every $(\eJ-\eL)$-residue of $R$ that meets
 $\opp_R(S)$ is entirely contained in $\opp_R(S)$.
By Corollary~\ref{cor:partition of opp}, this happens if and only if $\eJ-\eL \sbe \opp_\eJ(\eJ-\eL^*)= \eJ-\opp_\eJ(\eL^*)$, that is, if $\opp_\eJ(\eL^*)\sbe \eL$.
\qed

\bde\label{dfn:H}
Given a $\eJ$-residue $R$, a subset $\eL\sbe \eJ$ and a dual point
 $x^*$ in $\cP_\eL(R)^*$, we shall define the following subspace of $\cS_\eL(R)$:
 $$H_{R,\eL}(x^*)=\cP_{R,\eL}(\near_R(x^*)).$$
Note that by Lemma~\ref{lem:H disjoint from opp}, we have
$H_{R,\eL}(x^*)=
\cP_\eL(R)-\cP_{R,\eL}(\opp_R(x^*)).$
In case $R=\De$ or $\eL=\eK$, we shall drop that subscript from the notation.
As we shall see in Proposition~\ref{close=maxhyp}, $H(x^*)$ is often a (maximal) hyperplane of $\Gamma$. Hyperplanes of the form $H(x^*)$ are called {\dfn singular} or {\dfn attenuated}.
\ede

\bde\label{dfn:polarized embedding}
We call a full projective embedding $(e,V)$ of $\Gamma$ {\dfn polarized} if every singular hyperplane is induced by $(e,V)$.
In that case the {\dfn polar radical } of $e$ is the subspace
 $$\cR_e=\bigcap_{p^*\in\Gamma^*} \langle H(p^*)\rangle_e.$$
\ede

\bre
Our notion of ``polarized'' specializes to the notion of ``polarized'' defined in~\cite{CaDePa2007,ThVa2006} if we let $\eK=\{k\}$ refer to an end-node of a diagram $\sM$ of type $B_n$, $C_n$, or $\wA_n$.
Since in those cases the map $\opp_I$ is equal to the identity it happens that $\Gamma^*=\Gamma$ in loc. cit.. In the present paper we do not restrict ourselves to that situation. In particular, when $\sM=A_n$, $D_n$ ($n$ odd), or $E_6$ it may happen that $\Gamma^*$ and $\Gamma$ are different, if isomorphic, geometries.
\ere
\subsection{Main results and organization of the paper}
In Sections~\ref{section:polar hyperplanes}~and~\ref{section:covers and quotients} we consider a shadow space $\Gamma$ of a spherical Moufang building $\De$ associated to a twisted or untwisted Chevalley group $G$. The opposition relation gives us the right perspective on polarized embeddings of shadow spaces as we shall see; many fundamental properties of spherical buildings translate transparently into properties of polarized embeddings.
In particular, the following is proved in Section~\ref{section:polar hyperplanes}.
The case where $|\eK|=1$ was proved in~\cite{Bl1999,BlBr1998}.
\bmpr\label{close=maxhyp}
For any dual point $p^*$ of $\Gamma^*$,
\begin{itemize}
\AI{a} $H(p^*)$ is a hyperplane of $\Gamma$, and
\AI{b} $\cS_{\De,\eK}(\opp_\De(p^*))$ is a connected subgeometry of $\Gamma$, except if $\sM_{n,\eK}(\FF)$ is one of the following:
\begin{itemize}
\AI{i} $G_{2,\{1\}}(\FF_2)$ ($1$ denoting the short root), $G_{2,\{1,2\}}(\FF_2)$, $G_{2,\{1,2\}}(\FF_3)$,$\wF_{4,\{1\}}(\FF_2)$, $\wF_{4,\{1,2\}}(\FF_2)$,
\AI{ii} $C_{n,\eK}(\FF_2)$, or $F_{4,\eK}(\FF_2)$, where $n\ge 2$ and $\eK$ contains both nodes of the double bond in the diagram.
\end{itemize}
As a consequence, except in cases {\rm (i)} and {\rm (ii)}, $H(p^*)$ is a maximal subspace of $\Gamma$.
\end{itemize}
\empr
Thus, apart from a few exceptions, the geometry $\Gamma$ contains a maximal singular hyperplane for each point of the dual geometry $\Gamma^*$.
As a step up to the major results of the paper, Theorem~\ref{thm:residually polarized} says that singular hyperplanes of $\Gamma$ are ``residually singular''. As a consequence, polarized embeddings of $\Gamma$ are ``residually polarized''.

\bmth\label{thm:residually polarized}
Let $\De$ be a spherical building with type set $I$ and let $\Gamma$ be its $K$-shadow space for a non-empty set $K\sbe I$.
Let $R$ be a $\eJ$-residue of $\De$.
\begin{itemize}
\AI{a} For any dual point $p^{*_R}$ of $\cP_{\eK\cap\eJ}(R)^*$, there is some dual point $p^*$ of $\Gamma^*$ such that
$H_{R,\eK\cap\eJ}(p^{*_R})$ corresponds to $\cP_{\De,\eK}(R)\cap H(p^*)$ under the isomorphism $\cS_{\eK\cap\eJ}(R)\cong\cS_{\De,\eK}(R)$.
\AI{b} If $(e,V)$ is a polarized embedding of $\Gamma$, then
  the isomorphism $\cS_{\eK\cap \eJ}(R)\cong \cS_{\De,\eK}(R)$ composed with  $e$ yields a full polarized embedding of $\cS_{\eK\cap\eJ}(R)$ into $\PP(V_R)$; here $V_R=\langle e(x)\mid x\in \cP_{\De,\eK}(R)\rangle$.
\end{itemize}
\emth
Theorem~\ref{thm:residually polarized} is proved in Subsection~\ref{subsec:residually polarized}. In short, part (b) follows directly from part (a).  Part (a) follows from a fundamental relation between the opposition and projection maps of spherical buildings. This relation is conveyed in Theorem~\ref{thm:opispo}, which was proven in~\cite{BlBr1998} to show that in a ``Far away'' geometry obtained from $\Gamma$ by removing the hyperplane $H(p^*)$, the residues are themselves ``Far away'' geometries and that, as a consequence, this far away geometry has a Buekenhout-Tits diagram.

\mn
We also prove the following general property of polarized embeddings.
\bmpr\label{mprop:abspol implies minpol}
If a point-line geometry $\Theta$ possesses a polarized embedding $e$ and an absolute embedding $\abs{e}$, then $\Theta$ possesses a unique minimal polarized embedding $\check{e}$. Namely we have
 $\check{e}=\abs{e}/\cR_{\abs{e}}$.
\empr
We note that many shadow spaces of spherical buildings do have an absolute embedding~\cite{KaSh2001,BlPa2003}.

\mn
In Sections~\ref{section:Weyl embedding}~and~\ref{section:minimal polarized embedding} we consider the following setup.
Let $\De$ be obtained from an untwisted Chevalley group $G$ over a field $\FF$ (see Table~\ref{table:diagrams and groups}).
Fix a non-empty subset $\eK\sbe \eI$ and let $\Gamma$ be the $\eK$-shadow space of $\De$.
Let $V(\lambda_\eK)_\FF$ be the weak Weyl module for $G$ with maximal vector $v^+$ of highest weight $\lambda_\eK$, as defined in Subsection~\ref{subsec:Weyl module}, and let $V^0=\FF G v^+$.
Let $p$ be the point of $\Gamma$ corresponding to the parabolic group stabilizing $v^+$.
Then, the map
 $$\begin{array}{rl}
 e_\eK\colon \Gamma&\to \PP(V^0)\\
  gp&\mapsto \langle gv^+\rangle\\
  \end{array}$$
defines a full projective embedding for $\Gamma$.
\bmth\label{thm:Weyl embedding = polarized}\ 
\begin{itemize}
\AI{a} The embedding $e_\eK$ is polarized.
\AI{b} The codomain $V^0/\cR_{e_\eK}$ of the minimal polarized embedding relative to $e_\eK$ is the unique irreducible $G$-module $L(\lambda_\eK)_\FF$ of highest weight $\lambda_\eK$.
\end{itemize}
\emth

The general idea of the proof is the following. In Subsection~\ref{subsec:contravariant form}, following~\cite{Hu2006,Wo1972b,Ve1975} we define a $\tau$-contravariant bilinear form $\beta$ on $V(\lambda_\eK)_\FF$ with the property that $\beta(v^+,v^+)=1$.
This means that there is an involution $\tau$ of $G$ that interchanges two opposite Borel groups $B^+$ and $B^-$, such that, for $g\in G$ and $u,v\in V(\lambda_K)_\FF$, we have
 $$\beta(gv,u)=\beta(v,g^\tau(u)).$$
This form is symmetric and non-degenerate on $V^0$ and has the property that weight spaces corresponding to distinct weights are orthogonal with respect to $\beta$.
Then, in Subsection~\ref{subsec:Weyl is polarized}, we show that the subspace $\ker(\beta(v^+,-))$ of $V^0$ induces a singular hyperplane $H(p^*)$ of
 $\Gamma$, for some dual point $p^*$ opposite to $p$, and part (a) follows by contravariance.
From the preceding discussion it follows that the polar radical $\cR_{e_\eK}$ coincides with the
 radical of $\beta$. It is known, and not hard to prove that this radical is the unique largest submodule of
  $V^0$.
Thus, the quotient $V^0/\cR_{e_\eK}$ is an irreducible $G$-module of highest weight $\lambda_\eK$.

\bmth\label{thm:L = minimal polarized}
Let $\De$ be a spherical building obtained from an untwisted Chevalley group and let $\Gamma$ be its $\eK$-shadow space,
  for some non-empty type set $\eK\sbe \eI$.
If $\Gamma$ possesses an absolute embedding, then  the unique irreducible $G$-module $L(\lambda_\eK)_\FF$ of highest weight $\lambda_\eK$ affords the unique minimal polarized embedding for $\Gamma$.
\emth

\myparagraph{Motivation for this paper}
The main motivation for the style in which this paper is written is to exhibit a connection between the geometry of shadow spaces of buildings and the representation theory of groups of Lie type. 
This builds on the connection between a number of results that seems to be known only to a handful of colleagues, but is, to the best of my knowledge, not written down anywhere.
I've made an attempt to bring these results  (on embeddings, buildings, Chevalley groups and their highest-weight representations) together in the hope that this makes this connection as well as the new results presented here accessible to a wider audience.

As for the results presented here themselves, let  us mention at least two motivations. First of all, Theorem~\ref{thm:Weyl embedding = polarized} allows us to ``recognize'' the fundamental weight modules geometrically as the minimal polarized embeddings of $\Gamma$. This opens the door to studying (certain) modular representations of Chevalley groups by geometric means. A pilot project formulated by Blok, Cardinali and Pasini to study the decomposition series of the
 Weyl modules of the symplectic groups from this perspective and reinterpret the results from~\cite{PrSu1983a,Ad1986} is in progress.

Secondly, we'd like to show that the notion of a polarized embedding of a shadow space is rather fundamental, as it relates directly to fundamental properties of the corresponding spherical buildings. This is for instance evidenced in
 Theorem~\ref{thm:residually polarized}.

\Problems
\begin{itemize}
\AI{1} Modify the above results to include twisted Chevalley groups.
\AI{2} Give a geometric proof that if a geometry $\Gamma$ with sufficiently transitive automorphism group $G$ has a minimal polarized embedding, then that embedding is irreducible for the automorphism group $G$, assuming that the automorphisms in $G$ lift to linear isomorphisms of that embedding.
\AI{3} Use the above connection to find decomposition series of $G$-modules that afford embeddings for $\Gamma$.
\end{itemize}

\mn
In  Sections~\ref{section:minuscule weights},~\ref{section:projective grassmannians},~\ref{section:polar grassmannians},~and \ref{section:projective flag grassmannians} we study polarized embeddings of various shadow spaces in detail.

\section{Polarized embeddings of shadow spaces}\label{section:polar hyperplanes}
In this section we prove Proposition~\ref{close=maxhyp} which tells us when the set $H(p^*)$ is a
 maximal hyperplane of $\Gamma$. Then we prove Theorem~\ref{thm:residually polarized} which says that
 a full polarized embedding of $\Gamma$ induces full polarized embeddings on each residue of $\Gamma$.
Both are closely related to fundamental properties of buildings.

\subsection{In which shadow spaces is $H(p^*)$ a maximal hyperplane?}\label{subsection:singular}
Recall from Lemma~\ref{lem:H disjoint from opp}  that, for a dual point $p^*\in \Gamma^*$, the subsets $H(p^*)=\cP_{\De,\eK}(\near(p^*))$ and $\cP_{\De,\eK}(\opp(p^*))$ are disjoint subsets of $\Gamma$.
We begin by proving Proposition~\ref{close=maxhyp}.

\mn
\pf (of Proposition~\ref{close=maxhyp})
(a) Choose a line $l$. By definition this is the $\eK$-shadow of a $k$-panel $\pi$ for some $k\in \eK$.
Either $\pi\sbe \near(p^*)$ or $\pi\cap\opp(p^*)\ne\emptyset$.
In the latter case $\proj_\pi(p^*)$ is the unique chamber of $\pi\cap\near(p^*)$.
By Lemma~\ref{lem:H disjoint from opp}, a point $q$ of $\cP$ meeting $\pi$ is in $H(p^*)$ if and only if $q\cap \pi\in \near(p^*)$. Thus either one or all points of $l$ are in $H(p^*)$. This means that $H(p^*)$ is a hyperplane of $\Gamma$.

(b) It follows from~\cite{Br1993,Ab1996,BlBr1998,Bl1999} that under the restrictions (i) and (ii), the collection
 $\opp(p^*)$ of chambers opposite to $p^*$ is connected as a chamber system.
By Corollary~\ref{cor:partition of opp},  $\cP_{\De,\eK}(\opp(p^*))$ is exactly the collection of points all of whose chambers are opposite to some chamber of $p^*$. Thus, connectedness of the chamber system implies  connectedness of the subgeometry $\cS_{\De,\eK}(\opp(p^*))$ of $\Gamma$.
That $H(p^*)$ is a maximal hyperplane of $\Gamma$ now follows from the above and Lemma~\ref{lem:maximal hyperplanes}.
\qed

\subsection{Polarized embeddings are residually polarized}\label{subsec:residually polarized}
We prove Theorem~\ref{thm:residually polarized}, which says that singular hyperplanes of $\Gamma$ are ``residually singular''. That is, given a residue $R$ of $\De$, the $\eK$-shadow of a singular hyperplane intersected with the $\eK$-shadow $\cP_\eK(R)$ of $R$, is a singular hyperplane, or all, of $\cP_\eK(R)$.
As a consequence, polarized embeddings of $\Gamma$ are ``residually polarized''.
This result stems from Theorem~\ref{thm:opispo}, proved in~\cite{BlBr1998}, which describes
 the interaction of the opposition and projection maps, two of the most fundamental maps in the theory of spherical and twin-buildings.
For convenience we quote this result here.
\bth\label{thm:opispo}{\rm(Lemmas 3.5 and 3.6 of~\cite{Bl1999} and Lemmas 6.2 and 6.4 of~\cite{BlBr1998})}
\begin{itemize}
\AI{a} Let $W$ be a Coxeter building of spherical type. Then
for any two residues $R$ and $S$ of $W$ we have
$\opp_R(\proj_R(S)) = \proj_R(\opp_W(S))$.
\AI{b}
Let $R$ and $S$ be residues containing opposite chambers in the spherical building $\Delta$.
Then the set of residues meeting $\opp_{\Delta}(S)$ and
contained in $R$ equals the set of residues in $\opp_R(\proj_R(S))$.
In particular, every object in $\cG(R)$ meets $\opp_{\Delta}(S)$ (as a set of chambers) 
if and only if $\proj_R(S) = R$, that is, if and only if
$\opp_I(\typ S) \subseteq \typ R$.
\end{itemize}
\eth

\medskip\noindent
We shall now prove Theorem~\ref{thm:residually polarized}. 

\mn
\pf (of Theorem~\ref{thm:residually polarized})
(a) Let $p^{*R}$ be a dual point of $\cP_{\eK\cap\eJ}(R)^*$. For some apartment $\Sigma$ on $p^{*R}$,
 let $p^*$ be the unique dual point in $\Gamma^*$ meeting 
  $\opp_\Sigma(\opp_{R\cap\Sigma}(p^{*R}))$.
Then, by Theorem~\ref{thm:opispo} (a) we have $p^{*R}=\proj_R(p^*)$ and so by Theorem~\ref{thm:opispo} (b)  $p^*$ satisfies the claim.

(b) Clearly $(e,V_R)$ is a full embedding for $\cS_\eK(R)$.
Now consider a dual point $p^{*R}$ of $\cP_{\eK\cap\eJ}(R)^*$ and let $\Sigma$ and $p^*$ be as in (a).
Moreover, let $p$ be the point of $\Gamma$ meeting $\opp_{R\cap \Sigma}(p^{*R})\sbe \opp_{\Sigma}(p^*)$.
Then, $p\not\in H(p^*)$ and since $H(p^*)$ is induced by $V$, also $e(p)\not\in\langle H(p^*)\rangle_e$ .
It follows that also $H_R=V_R\cap \langle H(p^*)\rangle_e$  is a proper hyperplane of $V_R$.

To see that $H_R$ induces $H_R(p^{*R})$, note that, by (a) $\langle H_R(p^{*R})\rangle_e \le H_R$.
As both the former and the latter are hyperplanes of $V_R$, we must have equality. That is, $H_R(p^{*R})$ is induced by $(e,V_R)$.
\qed

\section{Covers and quotients of polarized embeddings}\label{section:covers and quotients}
In this section we obtain some general properties of polarized embeddings.
In particular we consider absolute and minimal polarized embeddings.
In~\cite{CaDePa2007} several of these results were obtained for dual polar spaces. As it turns out, many of these can be generalized to arbitrary shadow spaces.

\subsection{The polar radical}
We first consider an arbitrary point-line geometry $\Gamma$ and a full projective embedding $(e,V)$.

\bde\label{dfn:factoring subspace}
Call $R\le V$ a {\dfn factoring subspace} for $(e,V)$ if
\begin{itemize}
\AI{QE1} $R\cap e(p)=\{0\}$ for every point $p$ of $\Gamma$, and
\AI{QE2} for any two distinct points $p,q\in \Gamma$ we have $\langle R, e(p)\rangle_V\ne \langle R,e(q)\rangle_V$.
\end{itemize}
For any factoring subspace $R\le V$, we define the mapping $e/R$, called the {\em quotient of $e$ over $R$}, as follows:
$$\begin{array}{rl}
e/R\ \colon \Gamma & \to \PP(V/R)\\
 p& \mapsto \langle R,e(p)\rangle.\\
 \end{array}$$
\ede
The following are immediate (cf.~\cite{CaDePa2007}):
\ble\label{lem:QE}
Let $(e,V)$ be a projective embedding of a point-line geometry $\Gamma$ and let
 $R$ be a factoring subspace.
If $e$ is a full embedding, then so is $e/R$ and if $e/R$ is polarized, then so is $e$.
\ele

\bco\label{cor:polarized absolute embeddings}
If $\Gamma$ has a polarized embedding, then its absolutely universal embedding, if it exists, is also polarized.
\eco

We now return to the standard situation of the paper, where $\Gamma$ is the $\eK$ shadow space of the spherical building $\De$. Recall from Definition~\ref{dfn:polarized embedding} that the {\dfn polar radical} of a polarized embedding $(e,V)$ is the subspace
$$\cR_e=\bigcap_{p^*\in\Gamma^*} \langle H(p^*)\rangle_e.$$
Any subspace of $\cR_e$ is called a {\dfn radical subspace}. The significance of radical subspaces is given by the following result.

\bpr\label{prop:quotient polarized embedding}
Let $(e,V)$ be a polarized embedding of a shadow space $\Gamma$ and let $R\le V$.
\begin{itemize}
\AI{a} If $R\le \cR_e$, then $R$ is a factoring subspace and $e/R$ is polarized.
\AI{b} If $R$ is a factoring subspace such that $e/R$ is polarized, then $R\le \cR_e$.
\AI{c} $e/\cR_e$ is the minimal polarized quotient of $e$.
\end{itemize}
\epr
\pf
(a)
(QE1) Let $p$ be any point of $\Gamma$. Then it is opposite some dual point $p^*\in \Gamma^*$ so that
 $p\not\in H(p^*)$. Since $H(p^*)$ is induced by $e$ this means that $e(p)\cap R\sbe e(p)\cap \langle H(p^*)\rangle_e=\{0\}$ and so (QE1) is satisfied.

(QE2) Clearly $\langle R,e(p)\rangle\le \langle \cR_e,e(p)\rangle\le \langle H(x^*)\rangle_e$ for any dual point $x^*$ satisfying $p\in H(x^*)$.
Now let $p$ and $q$ be distinct points of $\Gamma$.
Pick an apartment $\Sigma$ on $p$ and $q$ and let $p^*=\opp_\Sigma(p)$ and $q^*=\opp_\Sigma(q)$.
Then $p^*$ is not opposite to $q$  and $q^*$ is not opposite to $p$.
Hence, $e(q)\not\in\bigcap_{x^*\in H^*(p)} \langle H(x^*)\rangle_e\ge \langle \cR_e,e(p)\rangle\ge \langle R,e(p)\rangle$  and
$e(p)\not\in \bigcap_{y^*\in H^*(q)}\langle H(y^*)\rangle_e\ge \langle \cR_e,e(q)\rangle\ge \langle R,e(q)\rangle$.
This proves (QE2).
By Lemma~\ref{lem:QE} $e/R$ is again a full projective embedding.

Moreover, since $R\le H(p^*)$ for every dual point $p^*$ and $\langle H(p^*)\rangle_e$ is a hyperplane of $V$,
 $\langle H(p^*)\rangle_{e/R}$ is again a hyperplane of $V/R$.
That is, $e/R$ is again polarized.

(b)
Suppose $R\not\le \cR_e$.
Then there is some dual point $p^*$ such that $R\not\le \langle H(p^*)\rangle_e$.
As a consequence, $\langle (e/R) (H(p^*))\rangle=V/R$, contradicting that $e/R$ is polarized.

(c) This follows immediately from (b).
\qed

\subsection{Absolutely universal embeddings}
Combining Corollary~\ref{cor:polarized absolute embeddings} with Theorem~\ref{thm:Weyl embedding = polarized} we find the following.

\bco\label{cor:abs is polarized}
Let $\De$ be a spherical building obtained from an untwisted Chevalley group and let $\Gamma$ be its $\eK$-shadow space, for some non-empty type set $\eK\sbe \eI$.
Then, the absolutely universal embedding of $\Gamma$, if it exists, is polarized.
\eco

In~\cite{KaSh2001} and~\cite{BlPa2003} it is shown that for many buildings $\De$ and type sets $\eK$, the shadow space $\Gamma$ does possess an absolutely universal embedding.
So Corollary~\ref{cor:abs is polarized} ensures that this embedding is polarized whenever that building is obtained from an untwisted Chevalley group.

For several special single-node shadow spaces $\Gamma$ even more is known. Namely, that a particular embedding induces {\em all} hyperplanes of $\Gamma$.  Then, since Veldkamp lines exist in those cases (See~\cite{Sh1995}), that particular embedding is the absolutely universal embedding. Clearly, in those cases the absolute embedding is polarized.
Table~\ref{table:fundemb affords all hyperplanes} lists some of these shadow spaces. The entries 1-6 come from~\cite{ShTh1992}.
In this table we find all embeddable non-degenerate polar spaces, (including those coming from untwisted Chevalley groups with diagrams $\sM_{n,k}$ equal to  one of $B_{n,1}$, $C_{n,1}$, $D_{n,1}$).
Note that all minuscule weight geometries of untwisted Chevalley groups are present: they are the geometries of type  $A_{n,k}$ (any $k$), $B_{n,n}$ ($n\ge 2$), $C_{n,1}$ ($n\ge 3$), $D_{n,1}$, $D_{n,n-1}$, $D_{n,n}$ ($n\ge 4$), $E_{6,1}$ and $E_{6,6}$, and $E_{7,7}$ 
(For a definition of a minuscule weight see Section~\ref{section:minuscule weights}).

\noindent
\begin{table}[h]
\begin{tabular}{lll}
\hline
 & Type of Geometry & Source \\
\hline
1) &  Desarguesian projective spaces of finite rank &~\cite{VeYo1910}\\
2) &  Embeddable non-degenerate polar spaces of  & \cite{BuLe1974,Di1980,Le1981,Ti1986a} \\
    &  rank at least $3$ and embeddable generalized  & \\
    &   quadrangles which do not  possess ovoids.  &   \\
3) &  All embeddable point-line geometries with & \cite{Ro1987} \\
 &  three points per line & \\
4)  & The Grassmannian of projective lines over a  & \cite{CoSh1991,CoSh1997b}\\
    & field ($A_{n,2}(\FF)$, $n>3$,  $\FF$ a field) and low rank  & \\
     &  geometries $D_{5,5}$ and $E_{6,1}$.  & \\
5)  & All embeddable Grassmannians $A_{n,k}(\FF)$, $n>2$,
 &\cite{Sh1992}\\
 &  $1<k<n-1$, $\FF$ a field. & \\
 6) & All half-spin geometries $D_{n,n}(\FF)$, $n>4$. & \cite{Sh1994}\\
 7) & All orthogonal spin geometries $B_{n,n}(\FF)$, where
  & \cite{ShTh1992}\\
  &  $\FF$ is such that $B_{2,2}(\FF)$ has no ovoids.  & \\
  8)  & The exceptional geometry $E_{7,7}(\FF)$, $\FF$ a field  & \cite{Sh1997}\\
\hline
\end{tabular}
\caption{Shadow spaces whose universal embedding affords all hyperplanes}\label{table:fundemb affords all hyperplanes}
\end{table}

\newpage
\subsection{Minimal polarized embeddings}

We first prove Proposition~\ref{mprop:abspol implies minpol}.

\pf (of Proposition~\ref{mprop:abspol implies minpol})
Let $e$ be any polarized embedding of the point-line geometry $\Theta$.
Let $\abs{\cR}=\cR_{\abs{e}}$ be the radical of the absolute embedding of $\Theta$ and set
  $\check{e}=\abs{e}/\abs{\cR}$ and $\check{V}=\abs{V}/\abs{\cR}$.
By Corollary~\ref{cor:polarized absolute embeddings} $\abs{e}$ is polarized as well. Since $\abs{e}$ is absolute, there is some subspace $R$ such that $\abs{e}/R=e$.  Since $e$ is polarized, it follows from Proposition~\ref{prop:quotient polarized embedding} that $R\le \abs{\cR}$.
The canonical maps corresponding to the inclusions $\{0\}\le R\le \abs{\cR}$ yield morphisms $\abs{e}\to e\to\check{e}$. 
Since $\abs{e}$ does not depend on $e$, $\check{e}$ is the 
 unique minimal polarized embedding of $\Theta$.
\qed

\mn
We return to the situation where $\Gamma$ is the $\eK$-shadow space of a spherical building $\De$ for some non-empty subset $\eK$ of $\eI$.
In Theorem~\ref{thm:residually polarized} it was shown that any full polarized embedding $e$ of $\Gamma$ induces a full polarized embedding on each of the residues of $\Gamma$.
The next result shows that the same is true if we replace ``full polarized'' by ``minimal full polarized''.
\bth\label{thm:residually minimal}
A minimal polarized embedding of $\Gamma$ induces a minimal polarized embedding on each residue of $\Gamma$.
\eth
\pf
Let $(e,V)$ be a minimal polarized embedding for $\Gamma$. That is, $\cR_e=\{0\}$.
Let $R$ be a $\eJ$-residue of $\De$.
Identify $\cS_{\eK\cap \eJ}(R)\cong\cS_{\De,\eK}(R)$.
Let $(e',V')$ be the embedding induced by $e$ on $\cS_{\De,\eK}(R)$,
 that is, $e'$ is the restriction of $e$ to $\cP_{\De,\eK}(R)$, with codomain $V'=\langle \cP_{\De,\eK}(R)\rangle_e$.
By Theorem~\ref{thm:residually polarized}, $(e',V')$ is polarized. It suffices to show that the polar radical $\cR_{e'}$ is trivial.

\mn
Instead we prove a slightly stronger claim  from which the result immediately follows. Namely
$$\cR_e\cap V'=\cR_{e'}.$$
First of all, if $x^*$ is not opposite to any point of $\cP_{\De,\eK}(R)$, then $$\langle H(x^*)\rangle_e\ge e(\langle \cP_{\De,\eK}(R)\rangle)=V'$$ and so 
$$\cR_e\cap V'=\bigcap_{x^*\in\cP^*}\langle H(x^*)\rangle_e\cap V' 
=\bigcap_{x^*}\langle H(x^*)\rangle_e\cap V',$$ 
where $x^*$ runs over all dual points in $\cP^*$ that are opposite some point of $\cP_{\eK}(R)$. 
 
Now suppose that $x^*$ is opposite some point $x$ of $\cP_{\De,\eK}(R)$.
By Theorem~\ref{thm:opispo} part (b), we have
 $H(x^*)\cap \cP_{\De,\eK}(R) = H_{R,\eK\cap \eJ}(x^{*R})$, 
  where $x^{*R}=\proj_R(x^*)$ and by part (a) of Theorem~\ref{thm:opispo}, $x^{*R}\in \cP_{\eK\cap\eJ}(R)$.
Moreover, by Theorem~\ref{thm:residually polarized}, (a), as $x^*$ runs over all dual points in $\cP^*$, $x^{*R}$ runs over all dual points in $\cP_{\eK\cap\eJ}(R)^*$.
Combining this with the previous equality, we find 
$$\cR_e\cap V'=\bigcap_{x^{*R}\in \cP_{\eK\cap \eJ}(R)^*} \langle H_{R,\eK\cap\eJ}(x^{*R})\rangle_{e'}=\cR_{e'}.$$
\qed

\section{The Weyl embedding}\label{section:Weyl embedding}
Let $\Gamma$ be the $\eK$-shadow space of the spherical building $\De$, 
 for some subset $\eK$ of $\eI$.
We prove Theorem~\ref{thm:Weyl embedding}, stated below, which says that if $\De$ is obtained from an untwisted Chevalley group $G$, then $\Gamma$ can be embedded into a subspace of the Weyl module of suitably chosen highest weight. We shall call this embedding the Weyl embedding. A special case, where $\Gamma$ is a single node shadow space was considered in~\cite{Bl1999,BlBr1998a}.

\subsection{$\Gamma$ obtained from a $BN$-pair of $G$}\label{subsec:Gamma from BN pair}
We shall assume that $\De$ is obtained from the universal Chevalley group $G=\bar{\sM}(\FF)$ over the field $\FF$ with Dynkin diagram $\sM$ over the set $I=\{1,2,\ldots,n\}$ via its thick $BN$-pair in the canonical way (see~\cite{Ro1989a,Ti1986a,We2003}).
That is, we let $\De$ be the chamber system whose set of chambers is $G/B$ and in which the distance function is given by
 $$\begin{array}{rll}
 \delta\colon \De\times\De & \to W=N/H,&\mbox{ where }H=B\cap N,\\
 (gB,hB)&\mapsto w,&\mbox{ where } Bg^{-1}hB=BwB.\\
 \end{array}$$
Let $c$ be the chamber corresponding to $B$ and let $\Sigma$ be the apartment
 corresponding to the collection of $B$ cosets $NB$.
For every $\eJ\sbe \eI$, let $P_\eJ$ be the standard parabolic subgroup of $G$ of type $\eJ$. This is the stabilizer of the $\eJ$-residue and the $\eI-\eJ$-flag on $c$.

We now construct the $\eK$-shadow space $\Gamma=(\cP,\cL)$ of $\De$ as follows:
$$\begin{array}{ll}
 \cP&=G/P_{I-\eK}\\
 \cL&=\bigcup_{k\in \eK}\cL_k\mbox{, where}\\
  \cL_k&=\{\{ghP_{I-\eK}\mid h\in P_{k}\}\mid g\in G\}\mbox{ for each }k\in \eK.\\
\end{array}$$
So $\cL$ is the union of $G$-orbits of the "fundamental" lines $\{h P_{I-\eK}\mid h\in P_{\{k\}}\}$. Incidence is symmetrized containment.

\bre
For the above construction there is no need to restrict to the case where $G$ is universal. As is proved for instance in~\cite[\S 3]{St1967} if $G'$ is any other Chevalley group of type $\sM$ over $\FF$, then 
 $G'$ is a central quotient of the universal Chevalley group $G$. Since the center of $G$ is contained in $B\cap N$, replacing $G$ by $G'$ in the above construction yields canonically isomorphic $\De$ and $\Gamma$.
\ere

\goodsofar
\subsection{The Weyl module $V(\lambda)_\FF^0$}\label{subsec:Weyl module}
Given a weight $\lambda$ of the complex semi-simple Lie algebra of type $\sM$ and a field $\FF$, we construct the Weyl module $\bV=V(\lambda)_\FF^0$ and an associated Chevalley group $G_{\lambda}(\FF)$ along the lines of~\cite[Ch.~4]{Hu2006}. 
Another, ultimately equivalent, approach was taken in~\cite{Wo1972b}.

\myparagraph{The Kostant $\ZZ$-form of the universal enveloping algebra.}
Let $\mfg_\CC$ be the semi-simple Lie algebra with Dynkin diagram $\sM$ indexed by $\eI=\{1,2,\ldots,n\}$ over $\CC$. Let $\mfU=\mfU(\mfg_\CC)$ be its universal enveloping algebra.
Fix a Cartan subalgebra $\mfh_\CC$ of $\mfg_\CC$ and a choice of positive roots $\Phi^+$ with fundamental roots
 $\alpha_1$, \ldots, $\alpha_n$, along with a Chevalley basis $\sfC$ in $\mfg_\CC$ consisting of
  co-roots $H_1,\ldots,H_n\in\mfh_\CC$, $H_i$ being the co-root associated to $\alpha_i$, as well as a positive root vector $X_\alpha$ and a negative root vector $Y_\alpha$ for each $\alpha\in \Phi^+$.
Let $\mfg_\ZZ$ be the $\ZZ$-span of this Chevalley basis $\sfC$; it is stable under the Kostant $\ZZ$-form $\mfU_\ZZ$ of $\mfU$, which is generated by $1$ along with all products $H_i(H_i-1)\cdots (H_i-a+1)/a!$ and all $X_\alpha^a/a!$ and $Y_\alpha^a/a!$. Then we can find a triangular decomposition $\mfU=\mfU^- \mfU^0\mfU^+$ relative to a fixed ordering of the set of all positive roots $\alpha_1,\ldots,\alpha_r$ compatible with $\mfU_\ZZ$. For example, $\mfU_\ZZ^+$ has as a $\ZZ$-basis all products
 $X_{\alpha_1}^{a_1}/(a_1!)\cdots X_{\alpha_r}^{a_r}/(a_r!)$.

\myparagraph{The Verma module $M(\lambda)_\CC$ and admissible lattices.}
Let $\sfE$ be the real subspace of the dual space $\mfh_\CC^*$ of $\mfh_\CC$ spanned by the root system $\Phi$. The set $\La_w =\{\lambda\in \sfE\mid \la(H_j)\in \ZZ\ \forall j\in \eI\}$
 forms a lattice and is called the {\em weight lattice}.
It has a $\ZZ$-basis $\{\lambda_i\mid i\in \eI\}$, where $\lambda_i$ is the fundamental dominant weight associated to the $i$-th node of $\sM$ (labeled as in~\cite{Bo1981}) and has the property that $\lambda_i(H_j)=\delta_{ij}$ for all $i,j\in \eI$.
Thus $\La_w$ is the $\ZZ$-span of $\{\la_i\mid i\in \eI\}$.
Roots are also weights and, accordingly, the {\em root lattice} $\La_r\le \La_w$ is the $\ZZ$-span of the fundamental system $\Pi=\{\alpha_i\mid i\in I\}$. 
There is a natural ordering on weights given by setting $\mu\prec \nu$, whenever $\nu-\mu$ is a sum of positive roots ($\mu,\nu\in \Lambda_w$).

Now given any weight $\lambda$, there is a $\mfg_\CC$-module  $M(\lambda)_\CC$, called the Verma module, having a unique maximal submodule $M'$ and unique simple quotient $V(\lambda)_\CC$ (also denoted $L(\lambda)_\CC$) (See e.g.\cite[Ch. 10]{Car2005}).
It is a well-known result that this simple module is finite-dimensional if and only if $\lambda$ is dominant. That is, $\lambda=\sum_{i\in \eI}n_i\lambda_i$, where $n_i\in \NN$ are not all zero.
Since this is the case we're interested in, we shall henceforth assume that $\lambda$ is dominant.  Of particular importance for us are the following weights  
$$\lambda_K=\sum_{k\in K}\lambda_k\mbox{ for some non-empty subset }\eK\sbe\eI. $$

The module $V(\lambda)_\CC$ is generated, as a $\mfg_\CC$-module, by a vector $v^+$ of weight $\lambda$.
We pause here to insert an important observation on weight spaces
 of a finite-dimensional $\mfg_\CC$-module $V$.
For any $\mu\in\mfh_\CC^*$, let $V_\mu=\{v\in V\mid H\cdot v=\mu(H)v\mbox{ for all }H\in \mfh_\CC\}$ be the weight space of $V$ corresponding to $\mu$.
The weight-lattice of $V$ is defined by 
 $\Lambda(V)=\{\mu\in\mfh_\CC^*\mid V_\mu\ne 0\}$. It satisfies
  $\La_r\le \Lambda(V)\le \La_w$~\cite[\S 3]{St1967}.
We shall denote $\Lambda(\lambda)=\Lambda(V(\lambda)_\CC)$.
Now $\lambda$ is the highest weight of $V(\lambda)_\CC$ in the sense that it is maximal in $\Lambda(\lambda)$ with respect to the ordering $\prec$. Accordingly we call $v^+$ a maximal or highest weight vector.

Following~\cite[Ch. 27]{Hu1972} an admissible lattice in $V(\lambda)_\CC$ is a finitely generated $\ZZ$-submodule of $V(\lambda)_\CC$ that spans $V(\lambda)_\CC$ over $\CC$, has $\ZZ$-rank at most $\dim_\CC(V(\lambda)_\CC)$ and is invariant under $\mfU_\ZZ$. The minimal choice would be $A_{\min{}}=V(\lambda)_\ZZ=\mfU_\ZZ v^+$. 
There is also a unique maximal admissible lattice, denoted $A_{\max{}}$; it can be obtained as the dual of a minimal admissible lattice in the dual module (See also Example~\ref{ex:sl2 beta}). Any other admissible lattice $A$ for $V(\lambda)_\CC$ satisfies $A_{\min{}}\le A\le A_{\max{}}$.
For any $\mu\in\Lambda(\lambda)$, we set $A_\mu=A\cap (V(\lambda)_\CC)_\mu$.

\mn
In the remainder of this section, unless otherwise specified, we shall work with minimal admissible lattice $A_{\min}=V(\lambda)_\ZZ$.

\bpr\label{prop:U_ZZ-modules}~\  
\begin{itemize}
\AI{a} $V(\lambda)_\CC$ is an irreducible $\mfg_\CC$-module spanned by a vector $v^+$ of highest weight $\lambda$.
\AI{b} $V(\lambda)_\CC$ is the direct sum of its weight spaces. That is, 
$V(\lambda)_\CC=\bigoplus_{\mu\in\Lambda(\lambda)} V(\lambda)_{\CC,\mu}$, where $V(\lambda)_{\CC,\lambda}=\CC v^+$.
\end{itemize}
Let $A$ be any admissible lattice in $V(\lambda)_\CC$
 and let $A_\mu=A\cap V(\lambda)_{\CC,\mu}$.
Then we have
\begin{itemize}
\AI{c} We have $A=\bigoplus_{\mu\in\Lambda(\lambda)} A_\mu$ and $A_\lambda=\ZZ v^+$ for suitable choice of $v^+$. 
\AI{d} For each weight $\mu$ of $V(\lambda)_\CC$, $A_\mu$ spans $V(\lambda)_{\CC,\mu}$.
\end{itemize}
\epr
\pf
Proofs can be found in the following references.
 (a)~\cite[Ch. 10]{Car2005},\cite[Ch. 21]{Hu1972};
 (b)~\cite[Ch. 2]{Car2005},\cite[Ch. 20]{Hu1972};
 (c)~\cite[Ch. 27]{Hu1972} In particular Theorem 27.1 and its proof.
 (d) Clearly, for each weight $\mu$ of $V(\lambda)_\CC$, we have
  $\langle A_\mu\rangle\le V(\lambda)_{\CC,\mu}$. Now note that, by definition, $A$ spans $V(\lambda)_\CC$. Combining (b) and (c), we see that $A$ spans $V(\lambda)_\CC$ if and only if, for each $\mu$, $\langle A_\mu\rangle=V(\lambda)_{\CC,\mu}$.
\qed

\mn
As an example the minimal admissible lattices for $\mfg_\CC=\mfsl_2(\CC)$ are described in Lemma~\ref{lem:sl2C-modules}
\ble~\label{lem:sl2C-modules}{\rm \cite[Ch. 7]{Hu1972}}
Let $\mfg_\CC=\mfsl_2(\CC)$ and let $\lambda=m\lambda_1$ for some $m\in\NN$. Also write $\alpha=\alpha_1$.
Then, we have 
\begin{itemize}
\AI{a} $A_{\min}=V(\lambda)_\ZZ$ is a free $\ZZ$-module with basis $\{v_0,v_1,\ldots,v_m\}$, where
 $v_i=(Y_\alpha^{i}/i!)v^+$.
 \AI{b} Setting $v_{-1}=v_{m+1}=0$, we have
$$\begin{array}{ccc}
 H_\alpha v_i=(m-2i)v_i,&
 X_\alpha v_i=(m-i+1)v_{i-1},&
 Y_\alpha v_i=(i+1)v_{i+1}.
 \end{array}$$
\AI{c}
The maximal lattice $A_{\max}$ corresponding to $A_{\min}$ is a free $\ZZ$-module with basis $\{f_0,f_1,\ldots,f_m\}$, where the action is given by
$$\begin{array}{ccc}
 H_\alpha f_i=(m-2i)f_i,&
 X_\alpha f_i=if_{i-1},&
 Y_\alpha f_i=(m+1-i)f_{i+1}.
 \end{array}$$
Here $f_{m+1}=f_{-1}=0$.
In fact we have $v_i={m \choose i}f_i$ for all $i$.
\end{itemize}
\ele

\myparagraph{Passage to an arbitrary field $\FF$.}
We now pass to an arbitrary field $\FF$. We follow Chapter 27 from~\cite{Hu1972}. Let $\mfg(\lambda)_\ZZ$ be the stabilizer in $\mfg_\CC$ of $V(\lambda)_\ZZ$; it contains the $\ZZ$-span of the Chevalley basis $\sfC$ for $\mfg_\CC$. In fact
 $$\mfg(\lambda)_\ZZ=\mfh(\lambda)_\ZZ\oplus\bigoplus_{\alpha\in\Phi^+} \ZZ X_\alpha \oplus\bigoplus_{\alpha\in\Phi^+} \ZZ Y_\alpha,$$
where $\mfh(\lambda)_\ZZ=\{H\in \mfh\mid \forall
 \mu\in\Lambda(\lambda): \mu(H)\in \ZZ \}$.
Now we pass to the field $\FF$ by setting
 $$\begin{array}{ll}
 V(\lambda)_\FF &= \FF\otimes_\ZZ V(\lambda)_\ZZ,\\
\mfg(\lambda)_\FF&=\FF\otimes_\ZZ\mfg(\lambda)_\ZZ.\\
 \end{array}$$
The module $V(\lambda)_\FF$ shall be called the {\em weak Weyl module} of highest weight $\lambda$ over $\FF$.

\bre\label{rem:variation of admissible lattice}
The construction of $\mfg(\lambda)_\ZZ$, $\mfh(\lambda)_\ZZ$ and $V(\lambda)_\FF$ can be done using any admissible lattice $A$ instead of the minimal lattice $V(\lambda)_\ZZ$.
It is shown in loc. cit.\ that up to isomorphism, $\mfg(\lambda)_\ZZ$ and $\mfh(\lambda)_\ZZ$  only depend on $V(\lambda)_\CC$ (or in fact $\Lambda(\lambda)$), not on the choice of $A$.
However, it does affect the  $\mfg(\lambda)_\FF$-action on $V(\lambda)_\FF$, as one can deduce from Lemma~\ref{lem:sl2C-modules}.
\ere

\bre\label{rem:reduced weight space}
For any weight $\mu\in \La(\lambda)$ and $a\in A_\mu$, we have  $(1\otimes H_i)\cdot (1\otimes a)=\mu(H_i)\otimes a$, where
 $\mu(H_i)\in \FF$.
Thus, if $\FF$ has characteristic $0$,  distinct weights of $V(\lambda)_\CC$ induce distinct weights of $V(\lambda)_\FF$.
On the other hand, if $\FF$ has positive characteristic $p$, then it may occur for distinct weights $\mu$ and $\nu$ of $V(\lambda)_\CC$ that 
 $\mu(H_i)\equiv \nu(H_i)\mmod{p}$ for all $i=1,2,\ldots,n$.
In particular $\FF\otimes_\ZZ A_\mu$ and $\FF\otimes_\ZZ A_{\nu}$ belong to the same weight space of $V(\lambda)_\FF$.
\ere

\bde\label{dfn:reduced weight space}
For each weight $\mu\in \La(\lambda)$, we shall set
$V(\lambda)_{\FF,\mu}=\FF\otimes_\ZZ A_\mu$
and call this the reduced $\mu$-weight space.
\ede

\mn
We continue the example from Lemma~\ref{lem:sl2C-modules} in Lemma~\ref{lem:sl2F-modules}.
\ble~\label{lem:sl2F-modules}
Let $\mfg(\lambda)_\FF=\mfsl_2(\FF)$ and let $\lambda=m\lambda_1$ 
 for some $m\in\NN$.
\begin{itemize}
\AI{a} Then, $V(\lambda)_\FF$ is an $\mfsl_2(\FF)$-module of dimension $m+1$. 
\AI{b}
The module $V(\lambda)_\FF$ has a basis $\{v_0,v_1,\ldots,v_m\}$ such that the formulas describing the $\mfsl_2(\FF)$-action are as in Lemma~\ref{lem:sl2C-modules}{\rm (b)}.

\AI{c} 
If $\Char(\FF)=0$ or $\Char(\FF)=p$ and $m$ is restricted, i.e.\ $0\le m < p$, then $V(\lambda)_\FF$ is cyclic as an $\mfsl_2(\FF)$-module.
\end{itemize}

\ele
\pf
Parts (a) and (b) follow immediately from Lemma~\ref{lem:sl2C-modules}. Part (c) follows from the fact that under these conditions all coefficients in the action formulas of Lemma~\ref{lem:sl2C-modules} are non-zero in $\ZZ$ and remain so on reduction modulo $p$.
\qed

\bre
If one constructs the module $V(\lambda)_\FF$ and the $\mfsl_2(\FF)$-action from the maximal lattice $A_{\max}$, then there is a basis $\{f_0,\ldots,f_{m}\}$ such that the formulas describing the $\mfsl_2(\FF)$-action are as in Lemma~\ref{lem:sl2C-modules}{\rm (c)}.
\ere

\myparagraph{The Chevalley groups}\label{subsec:Chevalley Group}
For any $\alpha\in \Phi^+$, and $n\in \NN$, the elements $x_{\alpha,a}=X_\alpha^a/a!$ and $y_{\alpha,a}=Y_\alpha^a/a!$ belong to $\mfU_\ZZ$ and hence leave $V(\lambda)_\ZZ$ invariant under their action on the module $V(\lambda)_\CC$. Hence they also induce endomorphisms of $V(\lambda)_\FF=\FF\otimes_\ZZ V(\lambda)_\ZZ$, via a representation that we shall call $\rho$. 
Note that for large $a$, these elements represent the null operator on $V(\lambda)_\CC$.
Therefore, for any constant $t\in \FF$, the element
$$x_\alpha(t)=\sum_{a=0}^\infty t^a \rho(x_{\alpha,a})$$
is an endomorphism of $V(\lambda)_\FF$. It is a purely formal fact that $x_\alpha(t)^{-1}=x_\alpha(-t)$.
We now let
$$\sM_{\Lambda(\la)}(\FF)=\langle x_\alpha(t),y_\alpha(t)\mid \alpha\in \Phi^+,t\in\FF\rangle\le \SL(V(\lambda)_\FF)\le\End(V(\lambda)_\FF).$$
The group operation is the multiplication of $\End(V(\lambda)_\FF)$, i.e.\ composition of endomorphisms.
This group is a Chevalley group of type $\sM$ over $\FF$.
Up to isomorphism it only depends on $\Lambda(\lambda)$, but not on the choice of $A$. We shall also use the notation $\sM_{\lambda}(\FF)=\sM_{\La(\lambda)}(\FF)$

The group $\sM_{\La_r}(\FF)$ is called the adjoint group and denoted $\sM_{\rm ad}(\FF)$. It can be obtained from the adjoint representation of 
 $\mfg_\CC$, whose highest weight is the "highest root", denoted $\alpha^*$ (see column 3 in Table~\ref{table:Chevalley groups}).
The group $\sM_{\La_w}(\FF)$ is called the universal group and denoted $\bar{\sM}(\FF)$.  
Given $\lambda$ we have $\La_r\le\La(\lambda)\le \La_w$ and central surjective homomorphisms  $\bar{\sM}(\FF)\to \sM_{\Lambda(\la)}(\FF)$ and $\sM_{\La(\lambda)}(\FF)\to \sM_{\rm ad}(\FF)$.
We also have $Z(\bar{\sM}(\FF))\cong \Hom(\La_w/\La_r,\FF^*)$ \cite[\S 3]{St1967}. 
In particular, we can always view $V(\lambda)_\FF$ as a $\bar{\sM}(\FF)$-module, regardless of the choice of $\lambda$.

For types $E_8$, $F_4$ and $G_2$ we have $\La_r=\La_w$ so that $\sM_{\rm ad}(\FF)=\bar{\sM}(\FF)$.
Table~\ref{table:Chevalley groups} lists the possible Chevalley groups for the other spherical types, along with some weights giving rise to these groups.
\begin{table}[h]
\renewcommand{\arraystretch}{1.2}
$$\begin{array}{@{}c@{}||@{}c@{}|l@{}l@{}|@{}l@{}l@{}|@{}l@{}l}
\sM & \La_w/\La_r & \multicolumn{2}{c|}{\sM_{\rm ad}}  & \multicolumn{2}{c|}{\sM_{\La'}  \mbox{ with }}  & \multicolumn{2}{c}{\bar{\sM} } \\
 & & & & \multicolumn{2}{c|}{\La_r<\La'<\La_w}  & & \\
\hline
A_n (n\ge  2) &\ZZ_{n+1} & \PSL_{n+1}, & \lambda_1+\lambda_n; & & &  \SL_{n+1}, & \lambda_1 \\ 
\hline
B_n (n\ge 3) & \ZZ_2 & \PSO_{2n+1}, &\lambda_2; \lambda_i (i\ne n) & && \Spin_{2n+1}, & \lambda_n \\
\hline
C_n (n\ge 2) &\ZZ_2 &  \PSp_{2n}, &2\lambda_1; \lambda_i (i\mbox{ even}) & &  &\Sp_{2n}, & \lambda_i (i\mbox{ odd}) \\
\hline
D_{4}  & \ZZ_2\times\ZZ_2 &\PSO_{8}^+, & \lambda_2;  & \SO_{8}^+, & \lambda_1,\lambda_3,\lambda_4 & \Spin_{8}^+ & \\
\hline
D_{n}  (n\ge 5, \mbox{odd}) & \ZZ_4 &\PSO_{2n}^+, & \lambda_2; \lambda_{i} (i\le l-2, \mbox{even}) & \SO_{2n}^+, & \lambda_{i} (i\le l-2, \mbox{odd}) & \Spin_{2n}^+, & \lambda_{n-1},\lambda_n \\
\hline

D_n  (n\ge 6, \mbox{even}) &\ZZ_2\times\ZZ_2 & \PSO_{2n}^+, & \lambda_2; \lambda_{i} (i\le l-2,\mbox{even}) & \SO_{2n}^+, & \lambda_{i} (i\le l-2,\mbox{odd}) & \Spin_{2n}^+, & \lambda_{n-1},\lambda_n \\
\hline
E_6  &\ZZ_3 & {E_6} ,& \lambda_2; \lambda_4 & && \bar{E}_6,&\lambda_1, \lambda_3,\lambda_5,\lambda_6 \\
\hline
E_7 & \ZZ_2 & E_7, & \lambda_1; \lambda_3,\lambda_4,\lambda_6 &  & & \bar{E}_7, & \lambda_2,\lambda_5,\lambda_7\\
\end{array}$$
\renewcommand{\arraystretch}{1}
\caption{Chevalley groups}\label{table:Chevalley groups}
\end{table}
Most of Table~\ref{table:Chevalley groups} comes from~\cite[\S 3]{St1967}. In columns 3, 4, and 5, we list the group along with some weights $\lambda$ for which $\sM_{\La(\lambda)}$ is the desired group.
Dynkin diagram labelings are as in~\cite{Bo1981}. 
The first weight in column 3 is $\alpha^*$.
The remaining weights can be found from Chapters 8 and 13 in~\cite{Car2005} with an easy calculation.

\bre
Steinberg's table~\cite[\S 3]{St1967} also points out that $\PSO_{2n+1}\cong\SO_{2n+1}$.
From Ree~\cite{Ree1957}, we also have  $B_{n,{\rm ad}}=\POm_{2n+1}$, $\bar{B}_n=\Om_{2n+1}$, $D_{n,{\rm ad}}=\POm_{2n}^+$, and $\bar{D}_n=\Om_{2n}^+$.

There does not seem to be a standard notation distinguising the adjoint and universal Chevalley groups of types $E_6$ and $E_7$. We shall write $E_n(\FF)=E_{n,{\rm ad}}(\FF)$ (and $\bar{E}_n(\FF)=\bar{E}_n(\FF)$).
\ere

\mn
From now on we shall write 
$G_{\rm ad}=\sM_{\rm ad}(\FF)$, $G=\bar{\sM}(\FF)$, and $G_{\lambda}(\FF)=\sM_{\La(\lambda)}(\FF)$.

\bde\label{dfn:Weyl module}
The module $V(\lambda)_\FF$ is in general not cyclic. In some cases it may depend on the choice of $A$. However, when for instance $\infty>\dim(V(\lambda)_\CC)>|G_\lambda(\FF)|$, so that $V(\lambda)_\FF$ has  $\FF$-dimension strictly greater than the size of the group algebra, the module cannot be cyclic.

We shall denote the $G_{\lambda}(\FF)$-module generated by $v^+$ by $\bV=V(\lambda)_\FF^0$ and call this the {\em Weyl module} of highest weight $\lambda$. For $V(\lambda)_\FF$ itself we reserve the name {\em weak Weyl module}. Note that in some cases, such as when $\Char(\FF)=0$, these modules coincide. 
\ede

We illustrate what may happen using the example from Lemma~\ref{lem:sl2F-modules}.
\ble~\label{lem:SL2F-modules}
Let $\sM=A_1$. Write $\alpha=\alpha_1$
and let $\lambda=m\lambda_1$ for some $m\in\NN$.
Let $V=V(\lambda)_\FF$ be obtained from $V(\lambda)_\ZZ$.
Then all of the following hold.
\begin{itemize}
\AI{a}  The $G_\lambda(\FF)$-module $V$ has dimension $m+1$ and $\FF$-basis $\{v_0,v_1,\ldots,v_m\}$, where $v_i=(Y_\alpha^{i}/i!)v^+$.

\AI{b} If $\FF$ has sufficiently many elements, then $V$ is cyclic and generated by $v^+$ as a $G_\lambda(\FF)$-module.
This is the case if $\FF$ is infinite or if $\FF$ is finite, but $|\FF|\ge m+1$.

\AI{c} If $\Char(\FF)=0$, or $\Char(\FF)=p>0$ and $m$ is restricted, then $V$ is irreducible as a $G_\lambda(\FF)$-module.
\end{itemize}
 \ele
\pf
(a) Follows immediately from Lemma~\ref{lem:sl2F-modules}.

(b) Consider the element $y_{\alpha}(t)=\sum_{i=0}^m \left(\frac{Y_{\alpha}^i}{i!}\right) t^i$.
We have $y_{\alpha}(t)v^+=\sum_{i=0}^m t^i v_i$. Therefore, if $\FF$ contains distinct elements $t_0,\ldots,t_m$, then $(t_j^i)_{i,j=0}^m$ is a Vandermonde matrix. Hence the elements
 $y_{\alpha}(t_0)v^+$, \ldots, $y_{\alpha}(t_m)v^+$ are linearly independent over $\FF$.

(c) 
If $V'$ is any submodule then it has a highest weight vector. Thus this highest weight vector is (a scalar multiple of) $v_i$ for some $i$.
Note that, due to the fact that $m$ is restricted, all of the coefficients appearing in part (b) of Lemma~\ref{lem:sl2C-modules} are non-zero over $\FF$. Hence an argument similar to that of part (b) shows that
 $V'$ also contains $v_{i+1},\ldots,v_m$.
Applying the argument of (b) to $x_\alpha(t) v_m$, we see that
 $V'$ also contains $v_0,\ldots,v_{i-1}$ and so $V'=V$.
 \qed

\bre\label{rem:not cyclic}  Suppose that we pick some finite field $\FF_p$, $p$ prime and 
 a weight $\lambda=m\lambda_1$ with $m>|\SL_2(\FF_p)|$.
Then, from Lemma~\ref{lem:SL2F-modules} we see that while $V(\lambda)_\CC$
 is irreducible for $\SL_2(\CC)$, the module $V(\lambda)_{\FF_p}$ is not cyclic since its dimension exceeds that of the group algebra $\FF_p\SL_2(\FF_p)$.
So in this case $\bV=V(\lambda)_{\FF_p}^0$ is a proper submodule of 
 $V(\lambda)_{\FF_p}$.
\ere

\myparagraph{The $BN$-pair of $G_\lambda(\FF)$}
To finish this subsection we identify some relevant subgroups of $G_\lambda(\FF)$, following~\cite{Ca1972}.

The group 
$$G_\alpha=\langle x_\alpha(t),y_\alpha(t)\mid t\in \FF\rangle$$
 is a Chevalley group over $\FF$ with diagram $A_1$.
By the above discussion therefore, $G_\lambda(\FF)$ is a quotient of $\SL_2(\FF)$, the universal group of type $A_1$, and it has the adjoint group 
 $\PSL_2(\FF)$ as a quotient.
More precisely, for each $\alpha$ we have a surjective homomorphism
$$\begin{array}{rl}
 \phi_\alpha\colon \SL_2(\FF)&\to G_\alpha\\
 \left(\begin{array}{@{}cc@{}} 1 & t \\ 0 & 1\end{array}\right)&\mapsto x_\alpha(t)\\
 \left(\begin{array}{@{}cc@{}} 1 & 0 \\ t & 1\end{array}\right)&\mapsto y_\alpha(t)\\
 \end{array}$$
The Steinberg presentation theorem~\cite{St1967} says that in case $G_\lambda(\FF)$ is universal, $G_\lambda(\FF)$ can be viewed as being generated abstractly by elements of the form
 $x_\alpha(t)$ and $y_\alpha(t)$, ($\alpha\in \Phi^+$, $t\in \FF$) subject to the Chevalley commutator relations, and such that the above homomorphisms are in fact isomorphisms.

\mn
For each $\alpha\in\Phi^+$ and $t\in \FF$, let 
$$\begin{array}{rlrl}
h_\alpha(t)&=\phi_\alpha\left(\begin{array}{@{}cc@{}}t & 0 \\ 0 & t^{-1}\\\end{array}\right),&
n_\alpha(t)&=\phi_\alpha\left(\begin{array}{@{}cc@{}}0 & t \\ -t^{-1} & 0\\\end{array}\right).
\end{array}$$
For each $\alpha\in\Phi^+$, let 
$$\begin{array}{llrl}
U_\alpha&=\langle x_\alpha(t)\mid t\in \FF\rangle, &
H_\alpha&=\langle h_\alpha(t)\mid t\in \FF\rangle,\\ 
U_{-\alpha}&=\langle y_\alpha(t)\mid t\in \FF\rangle, &
N_\alpha&=\langle n_\alpha(t)\mid t\in \FF\rangle. 
\end{array}$$
We now have the following distinguished subgroups of $G_\lambda(\FF)$:
$$\begin{array}{lll}
 U^+=\langle U_\alpha\mid \alpha\in \Phi^+\rangle,&
 H=\langle H_\alpha\mid \alpha\in \Phi^+\rangle,&
 B^+ =U^+H,\\
 U^-=\langle U_\alpha\mid \alpha\in \Phi^-\rangle, &
 N=\langle N_\alpha\mid \alpha\in \Phi^+\rangle, &
 B^- =U^-H.\\
 \end{array}$$
We also set $n_i=n_{\alpha_i}(1)$, for all $i\in \eI$.
\bth\label{thm:BN pair} (see~\cite{Ca1972})
Let $B=B^+$. Then, the pair $(B,N)$ is a BN-pair of type $\sM$ for $G_\lambda(\FF)$.
More precisely, setting $W=N/H$ and $r_i=n_i H$, for $i\in \eI$, the pair $(W,\{r_i\}_{i\in \eI})$ is a Coxeter system of type $\sM$. 
\eth

\subsection{The Weyl embedding}
Continuing the notation from this section we now let 
 $$\begin{array}{lll}
 \lambda=\lambda_\eK,\ & \bV=V(\lambda)_\FF^0,\ & G=\bar{\sM}(\FF).
 \end{array}$$
Let $U=U^+$, $U^-$, $B=B^+$, $B^-$, $N$ and $H$ be defined as in Subsection~\ref{subsec:Chevalley Group}, starting with the universal Chevalley group $G$.
We assume that $\De$ is obtained from the $BN$-pair $(B,N)$ of $G$.
We also assume that $\Gamma$ is the $\eK$-shadow space obtained from $\De$ as in Subsection~\ref{subsec:Gamma from BN pair}

\bth\label{thm:Weyl embedding}
The shadow space $\Gamma=(\cP,\cL)$ has a full projective embedding into the Weyl module $\bV$ for $G$ as follows
$$\begin{array}{rl}
e_W \colon \cP & \to \PP(\bV)\\
 g p & \mapsto \langle g v^+\rangle,\\
\end{array}$$
for any choice of a point $p$ of $\Gamma$.
\eth
\pf
The proof is the same as the proof of Lemma 3.1 in~\cite{BlBr1998a}, after noting that the weight $\lambda=\lambda_\eK$ is restricted for any characteristic, and $\bV$ is cyclic by definition.
\qed

\bre
Clearly if we consider an arbitrary weight $\lambda=\sum_{i\in \eI} m_i\lambda_i$, where $m_i\in\NN$, and we define the {\em support} of $\lambda$
 to be the set $\supp(\lambda)=\{i\in \eI\mid m_i>0\}$, then 
  $\bV=V(\lambda)_\FF^0$ affords some kind of embedding of the shadow space of 
   type $\eI-\supp(\lambda)$.
The arguments used above show that the point-set is mapped into the 
    set of $1$-spaces of $\bV$.
The standard line of type $k$ is now mapped to an arc in the $(m_k+1)$-dimensional subspaces $\bV_k=\langle g v^+\mid g\in P_{i}\rangle_\bV$. Such modules are well-known, so studying such embeddings is within reach and some of them might be interesting.
\ere

\goodsofar

\section{The minimal polarized embedding obtained from the Weyl embedding}\label{section:minimal polarized embedding} 
In this section we show that the Weyl embedding is polarized and that its minimal polarized quotient is the unique irreducible module 
  $L(\lambda_\eK)_\FF$ of
 highest weight $\lambda_\eK$.
\subsection{Opposite structures, actions and representations}
The reader familiar with opposite structures can skip this subsection.
In general, if $\bK$ is a category, then the opposite category $\bK^{\opp}$ is the category with the same objects as $\bK$, but with arrows reversed. Thus the identity map $\id\colon \bK\to \bK^{\opp}$  is a contravariant functor.
For example, if $G$ is a group with operation $\ast$, then the opposite group $G^{\opp}$ is the set $G$ equipped with the opposite operation $\ast^{\opp}$ given by $x\ast^{\opp} y=y\ast x$ for all $x,y\in G$.
Note that if $G$ is commutative, then $\id\colon G\to G^{\opp}$ is an isomorphism. 
Similarly, if $R$ is a ring with multiplication $\ast$, then the opposite ring $R^{\opp}$ is the additive group of $R$ equipped with the opposite multiplication $\ast^{\opp}$.  

Now let $R$ be a commutative ring with $1$.
The definition of the opposite of an associative $R$-algebra $A$ now follows from the definitions above. If $\mfL$ is a Lie algebra over $R$ with bracket $[\cdot,\cdot]$, then the opposite Lie algebra $\mfL^{\opp}$ is the $R$-module $\mfL$ equipped with the opposite bracket given by $[x,y]^{\opp}=[y,x]$.
For an associative $R$-algebra $A$ with multiplication $\ast$, let $\mfL(A)$ denote the Lie algebra on the $R$-module $A$ equipped with the bracket
 $[x,y]=x\ast y- y\ast x$. Then, we have
  $\mfL(A^{\opp})=\mfL(A)^{\opp}$.

An isomorphism between a group, ring, algebra, or Lie algebra $X$, and its opposite $X^{\opp}$ is called an anti-automorphism of $X$.

\subsection{The contravariant form on $V(\lambda)_\FF$
 and the irreducible quotient $L(\lambda)_\FF$.}\label{subsec:contravariant form}
In this section we follow the approach sketched in~\cite[Section 3.8]{Hu2006,Wo1972b,Ve1975}.
We first define an  automorphism $\tau$ of the vector space underlying $\mfg_\CC$ by setting, for each $\alpha\in\Phi$,  
(see e.g.~\cite{Hu1977}):
 $$\begin{array}{rlrlrl}
X_\alpha^\tau&=Y_\alpha, & H_\alpha^\tau=H_\alpha, & Y_\alpha^\tau&=X_\alpha.\\
 \end{array}$$
and extending $\CC$-linearly.
It follows from the isomorphism theorem for semi-simple complex Lie algebras~\cite[Ch. IV]{Jac1962} that $\tau$ induces an anti-automorphism of $\mfg_\CC$.
Naturally this extends to an isomorphism, also denoted $\tau$, between $\mfU=\mfU(\mfg_\CC)$ and $\mfU(\mfg_\CC^{\opp})$.  
We can view $\mfU(\mfg_\CC^{\opp})$ as $\mfU^{\opp}$ as follows.

First note that, as vector spaces, $\mfg_\CC=\mfg_\CC^{\opp}$ and that the tensor algebra $T(\mfg_\CC)$ and its opposite $T^{\opp}(\mfg_\CC^{\opp})$ are isomorphic via the map $r$, given on pure vectors by $x_1\otimes\cdots \otimes x_n\mapsto x_n\otimes \cdots\otimes x_1$ ($n\in \NN$, $x_i\in \mfg_\CC$).
Thus we can construct $\mfU=T(\mfg_\CC)/\mfI$ and 
  $\mfU(\mfg_\CC^{\opp})=T^{\opp}(\mfg_\CC^{\opp})/\mfI^{\opp}$
 (see e.g.~\cite{Car2005}), where 
 $\mfI=\langle x\otimes y-y\otimes x -[x,y]\mid x,y\in \mfg_\CC\rangle$  and $\mfI^{\opp}=\langle x\otimes^{\opp} y-y\otimes^{\opp} x - [x,y]^{\opp}\mid x,y\in \mfg_\CC^{\opp}\rangle$ are two-sided ideals
  of the associative algebras.
Now note that $\mfI=\mfI^{\opp}$ as subspaces of  the vector space $T(\mfg_\CC)=T^{\opp}(\mfg_\CC^{\opp})$, so that 
 as vector spaces $\mfU=\mfU(\mfg_\CC^{\opp})$. Hence $\mfU(\mfg_\CC^{\opp})$ can alternatively be constructed by taking the opposite 
    associative algebra of $\mfU$ and taking its Lie algebra, or, equivalently, by taking the opposite Lie algebra structure of $\mfU$.
Thus, $\tau$ extends to a proper anti-automorphism of 
the universal enveloping algebra  $\mfU$ and its underlying associative algebra structure inherited from $T(\mfg_\CC)$.
Since $\tau$ preserves $\sfC$, it restricts to a proper anti-automorphism of the Kostant $\ZZ$-form $\mfU_\ZZ$.

From $\tau$ one creates a symmetric bilinear form $\beta_\ZZ$ on $V(\lambda)_\ZZ$.
Namely, one first defines a twisted $\mfg_\CC$-module ${}^\tau V(\lambda)_\CC$. The module ${}^\tau V(\lambda)_\CC$ is the dual vector space  $V(\lambda)_\CC^*$ with a twisted action defined as follows: $$(g\cdot f)(v)=f(g^\tau(v))\mbox{ for all }g\in\mfg_\CC, f\in V(\lambda)_\CC^*, v\in V(\lambda)_\CC.$$
One verifies that ${}^\tau V(\lambda)_\CC$ has highest weight vector $f^+$ of weight $\lambda$, defined by $$f^+(v)=c  \mbox{ if } v\in cv^+ + \bigoplus_{\mu\ne \lambda}V(\lambda)_{\CC,\mu}.$$
The map $\phi\colon V(\lambda)_\CC\to {}^\tau V(\lambda)_\CC$ given by 
 $uv^+\mapsto uf^+$ ($u\in \mfU$) induces a unique $\mfU$-module isomorphism.
Clearly it restricts to a $\mfU_\ZZ$-module isomorphism $\phi\colon \mfU_\ZZ v^+\to \mfU_\ZZ f^+$.
We now define a bilinear form $\beta$ on $V(\lambda)_\CC$ by setting
$$\beta(v_1,v_2)=\phi(v_1)(v_2)\mbox{ for all }v_1,v_2\in V(\lambda)_\CC.$$
The form $\beta$ has the following properties (see~\S 3.8~and~\S 2.4 of~\cite{Hu2006}).
\ble\label{lem:contravariant form}\ 
\begin{itemize}
\AI{a} The bilinear form $\beta$ is contravariant. That is,
$$\beta(g u,v)=\beta(u,g^\tau v)\mbox{ for all }g\in \mfU,\mbox{and all } u,v\in V(\lambda)_\CC.$$
\AI{b} The form $\beta$ is symmetric and non-degenerate on $V(\lambda)_\CC$.
\AI{c}  Weight spaces corresponding to distinct weights of $V(\lambda)_\CC$
 are orthogonal with respect to $\beta$.
\AI{d} We have $\beta(v^+,v^+)=1$. Hence, it restricts to
 $\beta\colon A_{\min{}}\times A_{\min{}}\to \ZZ$.
Mutatis mutandis, (a), (b), and (c) also hold for this restriction.
\end{itemize}
\ele

We note that in Lemma~\ref{lem:contravariant form} part (c) is a consequence of (a).

\mn
Tensoring over $\ZZ$ with $\FF$, we obtain a symmetric $\FF$-bilinear contravariant 
 form $\beta$ on $V(\lambda)_\FF$.
\bco\label{cor:contravariant form}\
\begin{itemize}
\AI{a}
The form $\beta$ induced on $V(\lambda)_\FF$ is a symmetric $\FF$-bilinear 
 contravariant form such that $\beta(v^+,v^+)=1$. 
 
\AI{b}
Reduced weight spaces of $V(\lambda)_\FF$ corresponding to distinct weights of $V(\lambda)_\CC$ are orthogonal with respect to $\beta$.
Hence, 
 weight spaces corresponding to distinct weights of $V(\lambda)_\FF$ are orthogonal with respect to $\beta$.
\end{itemize}
Statements (a) and (b) also hold when $V(\lambda)_\FF$ is replaced with the submodule $V(\lambda)_\FF^0$ generated by $v^+$.
\eco

\bex\label{ex:sl2 beta}
Let $\mfg_\CC=\mfsl_2(\CC)$ and $\lambda=m\lambda_1$, for some $m\in \NN$.
The minimal admissible lattice $A_{\min}=V(\lambda)_\ZZ$ is the $\ZZ$-span of the basis $\cV=\{v_0,\ldots,v_m\}$ from Lemma~\ref{lem:sl2C-modules} part (b).
We now construct ${}^\tau V(\lambda)_\CC$.
Let $\{f_0,f_1,\ldots,f_m\}$ be the basis of $V(\lambda)_\CC^*$ dual to $\cV$, i.e.\ such that $f_i(v_j)=\delta_{i,j}$ for $0\le i,j\le m$. 
Using the formulas from Lemma~\ref{lem:sl2C-modules} part (b), together with the fact that $\tau$ fixes $H_{\alpha_1}$ and interchanges $X_{\alpha_1}$ and $Y_{\alpha_1}$ we find that the action of $\mfsl_2(\CC)$ on ${}^\tau V(\lambda)_\CC$ is given by the formulas in part (c) of that lemma.
The isomorphism $\phi\colon V(\lambda)_\CC\to {}^\tau V(\lambda)_\CC$ 
 is given by $v_i\mapsto{m\choose i}f_i$ for all $i=0,1,\ldots,m$.
Hence, with respect to the basis $\cV$, $\beta$ is given by the diagonal matrix with entry ${m\choose i}$ in the i-th row and column.
We have $\La(\lambda)=\{(m-2i)\lambda_1\mid i=0,1,\ldots,m\}$ with $A_{\min,(m-2i)\lambda_1}=\ZZ v_i$, for each $i$.
Lemma~\ref{lem:contravariant form} and Corollary~\ref{cor:contravariant form} are easily verified in this case.
Also, the maximal lattice in $V(\lambda)_\CC$ corresponding to $A_{\min}$ is
 $A_{\max}=\{f\in V(\lambda)_\CC\mid \beta(f,a)\in \ZZ\ \forall a\in A_{\min}\}$, which is the $\ZZ$-span of $\{v_i/{m\choose i} \mid i=0,\ldots,m\}$.
\eex

\goodsofar

\mn
We now extend the preceding results to the Chevalley groups. 
\bpr\label{prop:tau on G}
The map $\tau$ induces an anti-automorphism of $G_\lambda(\FF)$  that satisfies 
  $$x_\alpha(t)^\tau = x_{-\alpha}(t)\mbox{ for all }\alpha\in\Phi, t\in \FF.$$ 
\epr
Before we prove Proposition~\ref{prop:tau on G}, we introduce some notation connecting the various module structures involved with the corresponding representations.

Let $\rho\colon\mfg_\CC\to\mfgl(V(\lambda)_\CC)$ denote the representation corresponding to the left-module structure:
$\rho(x)(v)=x\star v$. 
Let $\langle\cdot,\cdot\rangle\colon V(\lambda)_\CC^*\times V(\lambda)_\CC\to \CC$ be the standard pairing given by $\langle f, v\rangle=f(v)$.
One verifies that the $\CC$-linear map
 $\rho^{\dagger}\colon\mfg_\CC^{\opp}\to\mfgl(V(\lambda)^*_\CC)$ denoted  
  $\rho^\dagger(x)(f)=x\star^\dagger f$ and given by 
   $\langle x\star^\dagger f, v\rangle=\langle f,x\star v\rangle$, for all $x\in \mfg_\CC^{\opp}$, $f\in V(\lambda)_\CC^*$, and $v\in V(\lambda)_\CC$ is a homomorphism.
The composition $\rho^\tau=\rho^\dagger\after\tau$ is the 
 representation $\rho^\tau\colon \mfg_\CC\to \mfgl(V(\lambda)_\CC^*)$  that turns $V(\lambda)^*_\CC$ into the $\mfg_\CC$-module ${}^\tau V(\lambda)_\CC$. We denote the action $x\star^\tau f=\rho^\tau(x)(f)$. It satisfies $x\star^\tau f=x^\tau\star^\dagger f$.
 
We then have natural extensions to the universal enveloping algebras $\mfU$, $\mfU^{\opp}$ that we denote $\rho$, $\rho^\tau$, and $\rho^\dagger$ as well.
The isomorphism $\phi$ sends $x\star v^+\mapsto x\star^\tau f^+$, for all $x\in\mfU$.
Setting $V(\lambda)_\ZZ=\mfU_\ZZ \star v^+$ 
 and $V(\lambda)^*_\ZZ=\mfU_\ZZ^{\opp}\star^\dagger f^+=\mfU_\ZZ\star^\tau f^+$, we see that $\phi$ restricts to 
 $\phi\colon V(\lambda)_\ZZ \to V(\lambda)_\ZZ^*$.
Since we also have $\langle f^+,v^+\rangle=1$, the pairing restricts to 
 $\langle \cdot,\cdot\rangle\colon V(\lambda)_\ZZ^*\times V(\lambda)_\ZZ\to \ZZ$.
Therefore we have restrictions
 $\rho\colon \mfU_\ZZ\to \mfgl(V(\lambda)_\ZZ)$, 
 $\rho^\dagger\colon \mfU^{\opp}_\ZZ\to \mfgl(V(\lambda)^*_\ZZ)$, and 
 $\rho^\tau\colon \mfU_\ZZ\to \mfgl(V(\lambda)^*_\ZZ)$, such that $\phi$ induces an isomorphism between $\rho$ and $\rho^\tau$, and, for every
  $f\in V(\lambda)^*_\ZZ$, $x\in \mfU_\ZZ$, and $v\in V(\lambda)_\ZZ$, we have
   $\langle x\star^\dagger f,v\rangle=\langle f,x\star v\rangle$.
Combining all this we find that, for $v_1,v_2\in V(\lambda)_\ZZ$, we have 
 $\beta(x\star v_1,v_2)=\langle \phi(x\star v_1),v_2\rangle
  = \langle x\star^\tau\phi(v_1),v_2\rangle=\langle x^\tau\star^\dagger\phi(v_1),v_2\rangle=\langle \phi(v_1),x^\tau\star v_2\rangle = \beta(v_1,x^\tau \star v_2)$, i.e.\ $\beta$ is $\tau$-contravariant.

As stated in Corollary~\ref{cor:contravariant form}, when we pass to the field $\FF$ these properties are preserved.
First tensor $\tau$ to get an isomorphism 
 $\tau_\FF\colon\mfU_\FF\to \mfU_\FF^{\opp}$.
Next, we tensor the representations $\rho$, $\rho^\dagger$ and $\rho^\tau$. For example, we define 
 $\rho_\FF\colon \mfU_\FF\to\mfgl(V(\lambda)_\FF)$ by $c\otimes x\mapsto c\otimes\rho(x)$ and with corresponding action $(c\otimes x)\star (d\otimes v)=cd\otimes x\star v$.
The others are defined similarly. Note that $\rho_\FF$ also defines a homomorphism between the underlying associative algebras.
One verifies that $\rho^\tau_\FF=\rho^\dagger_\FF \after \tau_\FF$.
Tensoring the $\mfU_\ZZ$-module isomorphism $\phi$ gives an isomorphism of $\FF$-vector spaces 
 $\phi_\FF\colon V(\lambda)_\FF\to V(\lambda)_\FF^*$. Conjugation by $\phi_\FF^{-1}$ gives an isomorphism of associative algebras 
  $\phi^*_\FF\colon \End(V(\lambda)_\FF)\to \End(V(\lambda)_\FF^*)$, which also induces an isomorphism between the corresponding Lie algebras
  $\phi^*\colon  \mfgl(V(\lambda)_\FF)\to \mfgl(V(\lambda)_\FF^*)$.
Thus we have the following commutative diagram of associative algebra morphisms.

\beq\label{eqn:rho diagram}
\xymatrix{
 \mfU_\FF  \ar[r]^{\rho_\FF} \ar[d]_{\tau_\FF}^{\cong} \ar[dr]^{\rho^\tau_\FF} & \End(V(\lambda)_\FF) \ar[d]^{\phi^*_\FF}_{\cong} \\
  \mfU^{\opp}_\FF \ar[r]_{\rho^\dagger_\FF} &   \End(V(\lambda)_\FF^*)\\ 
}
\eeq

\mn

\pf (of Proposition~\ref{prop:tau on G}) 
We shall construct the Chevalley group $G_\lambda(\FF)$
as well as  the group $G_\lambda^*(\FF)$  for $\mfg_\CC^{\opp}$.
We then show that $G_\lambda^*(\FF)$ is naturally isomorphic to the opposite group $G_\lambda(\FF)^{\opp}$ and that 
 $\tau$ induces an isomorphism $G_\lambda(\FF)^*\cong G_\lambda(\FF)^{\opp}$.
 
\mn
The group $G_\lambda(\FF)$ is generated as a subgroup of the multiplicative group of $\End(V(\lambda)_\FF))$ by  images under $\rho_\FF$ of elements of the form
$\bar{x_\alpha}(t)=\sum_{a=0}^\infty t^a\otimes x_{\alpha,a}$.
Thus we have $\rho_\FF(\bar{x_\alpha}(t))=\sum_{a=0}^\infty t^a\otimes\rho(x_{\alpha,a}) $ and the action on $V(\lambda)_\FF=\FF\otimes_\ZZ V(\lambda)_\ZZ$ is therefore given by $x_\alpha(t)\star (1\otimes v)=\sum_{a=0}^\infty t^a\otimes (x_{\alpha,a}\star v)$, for any $\alpha\in \Phi$, $t\in \FF$, and $v\in V(\lambda)_\ZZ$. 

If $\VV$ is the category of finite dimensional $\FF$-vector spaces and linear transformations and $\sfD\colon \VV\to \VV$ denotes the duality functor, sending each $V\in \VV$ to its dual $V^*\in\VV$ and
 each linear transformation $E$ to $E^*\colon f\mapsto f\after E$, then $\sfD$ is contravariant and $\sfD^{\opp}\colon \VV\to \VV^{\opp}$ is covariant. 
Now view  $\End(V(\lambda)_\FF)$ as the subcategory with single object $V(\lambda)_\FF$ in which the morphisms are those of $\VV$ having $V(\lambda)_\FF$ as starting point and end point.
Since $V(\lambda)_\FF$ is finite dimensional, we have $V(\lambda)_\FF^{**}=V(\lambda)_\FF$ and it follows that  $\sfD^{\opp}$ restricts to an isomorphism between 
 $\End(V(\lambda)_\FF)$ and $\End(V(\lambda)_\FF^*)^{\opp}$.

Composing $\sfD^{\opp}$ with $\rho_\FF\colon \mfU_\FF\to \End(V(\lambda)_\FF)$ we get a dual representation $\rho_\FF^{*\opp}=\sfD^{\opp}\after \rho_\FF\colon \mfU_\FF\to \End(V(\lambda)_\FF^*)^{\opp}$.
Taking opposites in the domain and codomain of $\rho_\FF^{*\opp}$, one verifies that we recover the representation 
$\rho^\dagger_\FF\colon \mfU_\FF^{\opp}\to\End(V(\lambda)^*_\FF)$ introduced above.
Thus all triangles in Diagram~(\ref{eqn:F diagram}), which includes Diagram~(\ref{eqn:rho diagram}), commute.

\beq\label{eqn:F diagram}
\xymatrix{
 \mfU_\FF  \ar[r]^{\rho} \ar[d]_{\tau}^{\cong} \ar[dr]^{\rho^\tau} & \End(V(\lambda)_\FF) \ar[d]^{\phi^*}_{\cong} \\
  \mfU^{\opp}_\FF \ar[r]_{\rho^\dagger}\ar[dr]_{\rho} &   \End(V(\lambda)_\FF^*)\\ 
  & \End(V(\lambda)_\FF)^{\opp} \ar[u]_{\sfD^{\opp}}^{\cong} \\
}
\eeq
Let $G_\lambda^*(\FF)$ be the  Chevalley group for $\mfg_\FF^{\opp}$ obtained using $\rho_\FF^\dagger$.
Clearly $\sfD\colon \End(V(\lambda)_\FF)\to \End(V(\lambda)^*_\FF)$ satisfies $\rho_\FF(\bar{x_\alpha}(t))\mapsto \rho^\dagger_\FF(\bar{x_\alpha}(t))$, for all $\alpha\in \Phi$ and $t\in \FF$, and, being bijective and contravariant, it restricts to an anti-isomorphism between $G_\lambda(\FF)$ and $G_\lambda^*(\FF)$
 .

From the commutativity of Diagram~(\ref{eqn:F diagram}), we deduce that $\tau$ induces an isomorphism $\phi^*\colon G_\lambda(\FF)\to 
 G_\lambda(\FF)^*$ satisfying  $\rho(\bar{x_\alpha}(t))\mapsto \rho^\dagger(\bar{x_{-\alpha}}(t))$ as well as an anti-automorphism: $\sfD^{\opp}\after\phi^*\colon G_\lambda(\FF)\to G_\lambda(\FF)^{\opp}$ which satisfies $\rho(\bar{x_\alpha}(t))\mapsto \rho(\bar{x_{-\alpha}}(t))$. This is the sought anti-automorphism.
\qed

\mn
We shall denote the anti-automorphism of $G_\lambda(\FF)$ induced by $\tau$ also by $\tau$.

It is now straightforward to verify the following.
\ble\label{lem:contravariant G-module}
For any $g\in G_\lambda(\FF)$ and $u,v\in V(\lambda)_\FF$ we have
 $$\beta(gu, v) = \beta(u, g^\tau v).$$
\ele

\mn
Recall from Subsection~\ref{subsec:Chevalley Group} that 
 $n_i=n_{\alpha_i}(1)$, for $i\in I$.
\ble\label{lem:W is beta isometry group}
The subgroup $N_0=\langle n_i\mid i\in \eI\rangle\le N$ of the universal Chevalley group $G$ is an isometry group for the form $\beta$.
\ele
\pf
First we note that for each $\alpha\in \Phi$ and $t\in \FF$ we have
 $x_\alpha(t)^\tau=x_{-\alpha}(t)$ and it follows that 
$n_\alpha(t)^\tau=(x_{\alpha}(t)x_{-\alpha}(-t^{-1})x_{\alpha}(t))^\tau=n_\alpha(-t^{-1})$. In particular, 
$n_{\alpha}(1)^\tau=n_{\alpha}(1)^{-1}$.
Recall also that $\tau$ is an anti-automorphism, meaning that
 $(gh)^\tau=h^\tau g^\tau$ for any $g,h\in G_\lambda(\FF)$.
It follows that if $n=n_{i_1}\cdots n_{i_l}$, then 
 $n^\tau = n^{-1}$.
Hence, for $u,v\in V(\lambda)_\FF$, we have
 $\beta(n u, nv)=\beta(u,n^\tau n v)=\beta(u,v)$.
\qed

\mn
It follows from Lemma~\ref{lem:W is beta isometry group} that weight vectors in one $N_0$-orbit have the same length with respect to $\beta$. 
It should be pointed out that in Lemma~\ref{lem:W is beta isometry group} $N_0$ cannot be replaced by $N$ in general.

\mn
Our motivation for introducing the contravariant form $\beta$ is 
 that it gives a connection between the Weyl module $\bV=V(\lambda)_\FF^0$ and the unique irreducible module $L(\lambda)_\FF$ of highest weight $\lambda$ (see See \cite{Hu2006}).

\bpr\label{prop:radical is maximal submodule}{\rm (cf.~\cite[\S 3.8]{Hu2006})}
The Weyl module $\bV=V(\lambda)_\FF^0$ has a unique maximal $G_\lambda(\FF)$-submodule and this submodule equals the radical of $\beta$ in $\bV$.
As a consequence, $\beta$ induces a non-degenerate contravariant form on the simple quotient $L(\lambda)_\FF$.
\epr
\pf The proof is the same as that of Proposition 3.8 of~\cite{Hu2006} after noting that $\bV=V(\lambda)_\FF^0$ is cyclic.
\qed

\subsection{The Weyl embedding is polarized}\label{subsec:Weyl is polarized}
Recall that $G_\lambda(\FF)$ is a central quotient of the universal Chevalley group $G$. Hence any $G_\lambda(\FF)$-module is automatically a $G$-module.

\bpr\label{prop:Weyl embedding is polarized}
The Weyl embedding of $\Gamma$ into $\PP(\bV)$ is polarized.
\epr
\pf
Let $\Sigma$ be the apartment of $\De$ corresponding to $N$ 
 and let $c^+$ be the chamber corresponding to $B$.
Let $p$ be the point of $\Gamma$ on $c^+$.

Recall that $\bV=V(\lambda)_\FF^0$.
For any weight $\mu\in \La(\lambda)$, let $\bV_\mu$ denote its reduced weight space as defined in Definition~\ref{dfn:reduced weight space}.
Recall that $e_W(p)$ is the subspace $\bV_\lambda$ of $\bV$ spanned by $1\otimes v^+$, where $v^+$ is the highest weight vector.
Also note that for any 
point $q$ of $\Sigma$ we have $q=wp$ for some $w\in W-\{1\}$ so that $e_W(q)=w \bV_\lambda= \bV_\mu$. Now if $p\ne q$, then since $\dim(\bV_\lambda)=1$, we must have $\mu=w\lambda\ne \lambda$.
Thus, $\beta(e_W(p),e_W(q))=0$ for all points $q$ on $\Sigma$ different from $p$.

Note that the chamber $c^-$ opposite to $c^+$ in $\Sigma$
 corresponds to $B^-=w_0 B^+ w_0$, where $w_0$ is the longest word in the Coxeter system $(W,\{r_i\}_{i\in \eI})$.
Now let  $p^*$ be the dual point of $\Gamma^*$ on $c^-$ and let 
 $H(p^*)$ be the hyperplane of $\Gamma$ consisting of points not opposite to $p^*$ in $\De$.
Then for every point $q'\in H(p^*)$ there is an apartment $\Sigma'$ on $q'$ and $c^-$. The group $B^-$ stabilizes $p^*$ while being regular on such apartments so there is some $b\in B^-$ with
 $b\Sigma=\Sigma'$ and there is some point $q$ on $\Sigma$ such that $bq=q'$. Since $q'$ was not opposite to $p^*$, $q\ne p$.

Turning back to the embedding we note that by contravariance and since $b^\tau\in B^+=B$, we have 
$\beta(e_W(p),e_W(q'))=\beta(e_W(p),b e_W(q))=\beta(b^\tau e_W(p),e_W(q))=\beta(e_W(p),e_W(q))=0.$
Thus the hyperplane $e_W(p)^\perp=\ker(\beta(v^+,-))$ contains 
$\langle H(p^*)\rangle_{e_W}$. 
Since $H(p^*)$ is a maximal subspace of $\Gamma$
 and using Lemma~\ref{lem:induced hyperplanes}, we find 
$$e_W(p)^\perp=\langle H(p^*)\rangle_{e_W}.$$
Part (b) of Lemma~\ref{lem:induced hyperplanes} moreover tells us that the hyperplane $H(p^*)$ is induced by $\bV$.
Since $G$ and hence also $G_\lambda(\FF)$ are transitive on dual points, and since $\beta$ is contravariant, the same holds for all other dual points of $\Gamma^*$.
Thus $\bV$ is polarized.
\qed

\bco\label{cor:contravariant radical is polar radical}
Let $e_W$ be the Weyl embedding of $\Gamma$ into the 
 Weyl module $\bV$ and let $\cR$ be the polar radical of $e_W$.
 Then the codomain of the minimal polarized embedding with respect to $e_W$ is the unique irreducible $G_\lambda(\FF)$-module $L(\lambda)_\FF$ of highest weight $\lambda$.
\eco
\pf
From Proposition~\ref{prop:Weyl embedding is polarized} it follows that 
$$\cR=\bigcap_{p^*\in \Gamma^*}\langle H(p^*)\rangle_{e_W}=
 \bigcap_{p\in \Gamma} e_W(p)^\perp=\Rad(\beta).$$
Therefore the codomain of the minimal polarized embedding
 with respect to $e_W$ is 
  $\bV/\cR=V(\lambda)_\FF^0/\Rad(\beta)$, which by 
Proposition~\ref{prop:radical is maximal submodule} is the unique irreducible $G_\lambda(\FF)$-module of highest weight $\lambda$.   
\qed

\mn
Theorem~\ref{thm:Weyl embedding = polarized} now follows from 
 Proposition~\ref{prop:Weyl embedding is polarized} and Corollary~\ref{cor:contravariant radical is polar radical}.

\section{Minuscule weight geometries}\label{section:minuscule weights}
Let $\lambda=\lambda_k$ be a fundamental dominant weight that is minuscule. This means by definition that the weight lattice of $V(\lambda)_\CC$ equals the orbit of $\lambda$ under the action of the Weyl group $W$. In particular, all weight spaces have dimension $1$.
Recall that $\lambda_k$ is a minuscule weight for the diagram $\sM$ if $\sM_k$ is one of the following: $A_{n,k}$ (any $k$), $B_{n,n}$ ($n\ge 2$), $C_{n,1}$ ($n\ge 3$), $D_{n,1}$, $D_{n,n-1}$, $D_{n,n}$ ($n\ge 4$), $E_{6,1}$ and $E_{6,6}$, or $E_{7,7}$. 
For names and dimensions of these embeddings see Table~\ref{table:minuscule weight embeddings}.
 
Call $e=e_W$ the embedding of $\Gamma=\Delta_k$ into $\bV=V(\lambda)_\FF^0$.
The weight spaces of $V$ are precisely the images of the point set of the apartment  $\Sigma=WP_{\eI-\{k\}}$ of $\Gamma$ corresponding to $W$, and in almost all cases these points generate $\Gamma$ (see~\cite{BlBr1998} for a precise statement).
By Theorem~\ref{thm:Weyl embedding = polarized} the embedding $\bV$ is polarized.
Moreover, $V(\lambda)_\FF=\bV=L(\lambda)_\FF$ since if $V(\lambda)_\FF$ had any proper submodule, it would be the direct sum of its weight spaces. But the weight spaces of
 $V(\lambda)_\FF$ are all of dimension $1$ and form a single orbit under $W$.
Hence, no proper submodule exists.
In view of Lemma~\ref{lem:W is beta isometry group} it also implies that $\bV$ has a basis of weight-vectors that is orthonormal with respect to $\beta$. Therefore, $\beta$ has trivial radical, which by Theorem~\ref{thm:Weyl embedding = polarized} again implies that $V(\lambda)_\FF=\bV$ is irreducible.

\mn
We finish this section with a brief remark on generating singular hyperplanes. Let $p$ be a point of $\Sigma$ opposite some dual point $p^*$ also on $\Sigma$. We now see that the hyperplane $\langle H(p^*)\rangle_e$ is exactly the
hyperplane of $\bV$ spanned by the set $\{e(q)\mid p\ne q\mbox{ a point  of }\Sigma\}$. It is proved in Blok~\cite{Bl1999} that often the hyperplane $H(p^*)$
 itself is generated, as a subspace of $\Gamma$ by the set of points of $\Sigma$ different from $p$; this is the case for instance if the diagram $\sM$ is one of
  $A_n$, $D_n$, $E_6$ or $E_7$.

\mn
In Table~\ref{table:minuscule weight embeddings} we list the 
 Weyl embeddings of the minuscule weight geometries.
Here $\Gamma$ is the $k$-shadow space of a building  associated to the universal Chevalley group $G=\bar{\sM}(\FF)$ and $V=V(\lambda_k)_\FF=\bV=L(\lambda_k)_\FF$.
\begin{table}[h]
$$
\begin{array}{c|c|l|c|c|c}
\sM & G &  k &  \bV & \dim(\bV) \\
\hline
A_n&  \SL_{n+1}(\FF) & k   & \mbox{ Grassmann }&  { n+1\choose k} \\
\hline
B_{n} & \Spin_{2n+1}(\FF) & n  & \mbox{ spin }    & 2^n \\
\hline
C_{n} & \Sp_{2n}(\FF)  & 1 &   \mbox{natural}    & 2n \\
\hline
D_{n} &  \Spin_{2n}^+(\FF)  & 1          & \mbox{natural} & 2n \\
\hline
           &                 & n,n-1 &      \mbox{half-spin} & 2^{n-1} \\
\hline
E_6          & \bar{E}_6(\FF)           & 1,6  & L(\lambda_k) & 27\\
\hline
E_7          & \bar{E}_7(\FF)           & 7 &  L(\lambda_7) & 56\\
\end{array}$$
\caption{Minuscule weight embeddings}\label{table:minuscule weight embeddings}
\end{table}

\goodsofar
\section{Grassmannians}\label{section:grassmannians}
As a preliminary to Sections~\ref{section:projective grassmannians},~\ref{section:polar grassmannians}~and~\ref{section:projective flag grassmannians}, we collect some information on tensor products and exterior powers of modules for the Lie algebra $\mfg=\mfg(\lambda)_\FF$ and its associated Chevalley group $G_\lambda(\FF)$. In particular, we shall study the form $\beta$ and the automorphism $\tau$.

\mn
It is well known (see e.g.~\cite{Car2005}) and easy to check that whenever $V_1,\ldots,V_l$ are $\mfg$-modules, then so is $V_1\otimes \cdots \otimes V_l$ under the action
\beq\label{eqn:Lie action on tensors}
g\cdot \otimes_{i=1}^l v_i=\sum_{i=1}^l v_1\otimes \cdots\otimes v_{i-1}\otimes gv_i\otimes v_{i+1}\otimes \cdots\otimes v_l,
\eeq
for all $g\in\mfg$ and $v_i\in V_i$ for all $i=1,\ldots,l$.
Similarly, if $V$ is a $\mfg$-module, then so is $\bigwedge^k V$
under the action
\beq\label{eqn:Lie action on exterior powers}
g\cdot \wedge_{i=1}^k v_i=\sum_{i=1}^k v_1\wedge \cdots\wedge v_{i-1}\wedge gv_i\wedge v_{i+1}\wedge \cdots\wedge v_k,
\eeq
for all $g\in\mfg$ and $v_i\in V$ for all $i=1,\ldots,k$.
As for the action of $G_\lambda(\FF)$, it is well known and it follows easily from Equation~(\ref{eqn:Lie action on tensors}) and the definition of $G_\lambda(\FF)$, that whenever $V_1,\ldots,V_l$ are $G_\lambda(\FF)$-modules, then so is $V_1\otimes \cdots \otimes V_l$ under the action
\beq\label{eqn:group action on tensors}
g\cdot \otimes_{i=1}^l v_i=g v_1\otimes \cdots\otimes  gv_i\otimes  \cdots\otimes gv_l,
\eeq
for all $g\in G_\lambda(\FF)$ and $v_i\in V_i$ for all $i=1,\ldots,l$.
Similarly, if $V$ is a $G_\lambda(\FF)$-module, then so is $\bigwedge^k V$
under the action
\beq\label{eqn:group action on exterior powers}
g\cdot \wedge_{i=1}^k v_i=gv_1\wedge \cdots\wedge  gv_i\wedge \cdots\wedge gv_k,
\eeq
for all $g\in G_\lambda(\FF)$ and $v_i\in V$ for all $i=1,\ldots,k$.

Next, we describe how a covariant or contravariant form on a collection of modules for $\mfg$ or $G_\lambda(\FF)$ induces a similar form on their tensor product or exterior powers. \ble\label{lem:form on tensor product}
Let $V_1,\ldots,V_l$ be finite dimensional $\FF$-vector spaces, let 
 $\sigma$ be an automorphism of $\FF$ 
 of order at most $2$ 
 and let 
  $\zeta_i$ be a $\sigma$-sesquilinear form on $V_i$. Then
\begin{itemize}
\AI{a} there is a unique $\sigma$-sesquilinear form $\zeta^{\otimes}$
 on $V_1\otimes \cdots \otimes V_l$ given by \\
 $\zeta^{\otimes}(u_1\otimes\cdots\otimes u_l,v_1\otimes\cdots\otimes v_l)=\Pi_{i=1}^l \zeta_i(u_i,v_i)$;
\AI{b}
if each $\zeta_i$ is non-degenerate, so is $\zeta^{\otimes}$;
\AI{c}
if each $\zeta_i$ is symmetric bilinear, so is $\zeta^{\otimes}$;
\AI{d}
if each $\zeta_i$ is skew-symmetric, then $\zeta^{\otimes}$ is 
 skew symmetric if $l$ is odd and symmetric otherwise;
\AI{e}
if each $V_i$ is a module for $\mfg$ and $\zeta_i$ is $\tau$-contravariant, then so is $\zeta^{\otimes}$;
\AI{f}
if each $V_i$ is a module for $G_\lambda(\FF)$ and $\zeta_i$ is $\tau$-contravariant, then so is $\zeta^{\otimes}$;
\AI{g}
if each $V_i$ is a module for $G_\lambda(\FF)$ and $G_\lambda(\FF)$ preserves $\zeta_i$, then $G_\lambda(\FF)$ preserves $\zeta^{\otimes}$.
\end{itemize}
\ele
\pf
(a) For each $i$ we have a $\sigma$-semilinear map
 $\phi_i\colon V_i\to V_i^*$ so that $\zeta_i(u,v)=\langle \phi_i(u),v\rangle$, where
 $\langle f,v\rangle=f(v)$ is the standard pairing $V_i^*\times V_i\to \FF$.
Note that $\phi^\otimes=\phi_1\otimes\cdots\otimes\phi_l\colon V_1\otimes \cdots \otimes V_l\to V_1^*\otimes \cdots \otimes V_l^*$ is again a $\sigma$-semilinear map.
We can compose this map with the standard pairing
  $V_1^*\otimes \cdots \otimes V_l^*\times V_1\otimes \cdots \otimes V_l\to \FF$
  given by $\langle f_1\otimes\cdots\otimes f_l,v_1\otimes\cdots\otimes v_l\rangle=\Pi_{i=1}^lf_i(v_i)$ to get the form $\zeta^\otimes$.
 It is immediate from this construction that $\zeta^\otimes$ is  $\sigma$-sesquilinear.
 (b) The standard pairing is non-degenerate and in this case $\phi_i$ and $\phi^{\otimes}$ are isomorphisms. 
(c) and (d) are trivial observations.
(e), (f), and (g) are easily seen to follow from Equations~(\ref{eqn:Lie action on tensors})~and~(\ref{eqn:group action on tensors}).
\qed

\ble\label{lem:form on exterior power}\label{lem:G preserves zeta wedge}
Let $V$ be a finite dimensional $\FF$-vector space, let 
 $\sigma$ be an automorphism of $\FF$ of order at most $2$ and let 
  $\zeta$ be a $\sigma$-sesquilinear form on $V$. Let $k\in \NN$ with $1\le k\le \dim(V)$.
Then
\begin{itemize}
\AI{a} there is a unique $\sigma$-sesquilinear form $\zeta^{\wedge }$
 on $\bigwedge^k V$ given by \\
  $\zeta^{\wedge }(u_1\wedge \cdots\wedge u_k,v_1\wedge\cdots\wedge v_k)=\det(\zeta(u_i,v_j))$;
\AI{b} if $\zeta$ is non-degenerate, so is $\zeta^{\wedge}$;
\AI{c}
if $\zeta$ is symmetric bilinear, so is $\zeta^{\wedge }$;
\AI{d}
if $\zeta$ is skew-symmetric, then $\zeta^{\wedge }$ is 
 skew symmetric if $k$ is odd and symmetric otherwise;
\AI{e}
if $V$ is a module for $\mfg$ and $\zeta$ is $\tau$-contravariant, then so is $\zeta^{\wedge }$;
\AI{f}
if $V$ is a module for $G_\lambda(\FF)$ and $\zeta$ is $\tau$-contravariant, then so is $\zeta^{\wedge }$;
\AI{g}
if $V$ is a  module for $G_\lambda(\FF)$ and $G_\lambda(\FF)$ preserves $\zeta$, then $G_\lambda(\FF)$ preserves $\zeta^{\wedge }$.
\end{itemize}
\ele
\pf
(a)
There is a $\sigma$-semilinear map $\phi\colon V\to V^*$ so that
 $\zeta(u,v)=\langle \phi(u),v\rangle$, where $\langle f,v\rangle=f(v)$ is the standard pairing $V^*\times V\to \FF$.
Note that $\phi^{\wedge}=\phi\wedge\cdots\wedge\phi\colon \bigwedge^k V\to \bigwedge^k V^*$ is again a $\sigma$-semilinear map.
We can compose this $\sigma$-semilinear map with the standard pairing
  $\bigwedge^k V^*\times \bigwedge^k V\to \FF$
  given by 
  $\langle f_1\wedge\cdots\wedge f_l,v_1\wedge\cdots\wedge v_l\rangle= \det(f_i(v_j))$ to get the form $\zeta^\wedge$.
Noting that $\det(\zeta(u_i,v_j))=\sum_{\rho\in\Sym(k)} {\sign(\rho)}\Pi_{i=1}^k \zeta(u_i,v_{\rho(i)})$ we see that $\zeta^\wedge$ is  $\sigma$-sesquilinear 
 because $\zeta^\otimes$ as defined in Lemma~\ref{lem:form on tensor product}, is $\sigma$-sesquilinear. 
(b) The standard pairing is non-degenerate and in this case $\phi$ and $\phi^\wedge$ are isomorphisms.
(c) This is because $\det(f_i(v_j))=\det(f_j(v_i))$. (d) Same as in Lemma~\ref{lem:form on tensor product}.
(e), (f) and (g). These follow from Equations~(\ref{eqn:Lie action on exterior powers})~and~(\ref{eqn:group action on exterior powers}) together with
the definition of the determinant as in (a).
\qed

\mn
Thus, we see that $\beta^\wedge$ and $\beta^\otimes$ are symmetric bilinear $\tau$-contravariant forms. Orthogonality of distinct weight spaces follows from contravariance as in Lemma~\ref{lem:contravariant form}.

\section{The projective Grassmannians}\label{section:projective grassmannians}
Let $\De$ be the building of type $A_n$ over the field $\FF$.
The universal Chevalley group is $G=\SL(V)$, where
 $V$ is a vector space of dimension $n+1$ over $\FF$.
Picking an ordered basis $\cA=\{a_1,\ldots,a_{n+1}\}$ for $V$,
 we identify $G$ with $\SL_{n+1}(\FF)$.
A BN-pair for $G$ is given by letting $B$ be the upper triangular matrix group
 and $N$ the monomial matrix group.
Then $H=B\cap N$ is the diagonal matrix group and the Weyl group $W=N/H\cong\Sym(n+1)$ in its action on the $1$-spaces spanned by the standard basis elements.

For each integer $k$ with $1\le k\le n+1$, the exterior power $\bigwedge^k V$ is clearly a module for $G$
 under the action
  $g(v_1\wedge \cdots\wedge v_k)=gv_1\wedge \cdots\wedge g v_k$.
For any non-empty subset $\eJ\sbe \eI$, let $a_J=\wedge_{j\in J} a_j$, where the $a_j$ appear with increasing subscripts.
Recall that  $\cA^k=\{a_J\mid J\sbe \{1,2,\ldots,n+1\}, |J|=k\}$
 is a basis for $\bigwedge^k V$.
One verifies that $v^+=a_{\{1,2,\ldots,k\}}$ is a vector of weight $\lambda_k$ that is stabilized by the subgroup $U^+$ of unipotent upper triangular matrices (see e.g.~\cite[Ch. 13]{Car2005}).
Thus $v^+$ is a vector of highest weight $\lambda_k$.
It is easy to see that $\lambda_k$ is minuscule: the collection of $1$-spaces spanned by elements from $\cA^k$ forms a single orbit under $W=\Sym(n+1)$. Thus, $\bigwedge^k V$ is the 
 irreducible Weyl module $V(\lambda_k)_\FF=\bV=L(\lambda_k)_\FF$.

\mn
The standard parabolic subgroup $P_{\eI-\{k\}}$ is precisely the stabilizer of the $k$-space $p=\langle a_1,\ldots,a_k\rangle_V$, which is a $k$-object 
of $\Delta$. The $k$-shadow space $\Gamma$ of $\De$ is the geometry whose points are the $k$-objects and where
 each line is the collection of points incident to some
 $\{k-1,k+1\}$-flag.
 
In accordance with Theorem~\ref{thm:Weyl embedding}, the standard embedding of $\Gamma$ into $\bigwedge^k V$ is given by
 $$\begin{array}{rl}
\cP & \to \PP(\bigwedge^k V)\\
p & \mapsto \wedge^k p\\
 \end{array}$$
where $\wedge^k p=\langle p_1\wedge \cdots \wedge p_k\rangle$ for some
 basis $p_1,\ldots,p_k$ of $p$. It is often called the {\em Grassmann embedding}.
This is well-defined since if $g\in\SL(V)$ induces a change of basis for $p$
 we have $g(p_1\wedge \cdots\wedge p_k)=d p_1\wedge \cdots\wedge p_k$,
  where $d$ is the determinant of the restriction of $g$ to $p$.

\mn
Next, we identify the form $\beta$ on the module $\bV$.
Since $\lambda_k$ is minuscule, Section~\ref{section:minuscule weights} tells us that $\cA^k$ forms an orthonormal basis for $\bigwedge^kV$ with respect to $\beta$.
Thus, if $\beta_1$ is the form on $V(\lambda_1)_\FF$, then 
 $\beta=\beta_1^\wedge$, as described in Lemma~\ref{lem:form on exterior power}.

Let us also identify $\tau$. The anti-involution $\tau$ of $G$ as described in Section~\ref{section:minimal polarized embedding} satisfies
 $x_\alpha(t)^\tau=x_{-\alpha}(t)$ for any $t\in \FF$ and $\alpha\in \Phi$.
In the present $A_n$ case, the root system is  
 $\Phi=\{\alpha_{i,j}\mid i,j\in \eI, i\ne j\}$ where 
  $\alpha_{i,j}=-\alpha_{j,i}$ with respect to the fundamental system 
   $\Pi=\{\alpha_{i,i+1}\mid i=1,\ldots,n\}$. 
With respect to the BN-pair chosen above we have
 $x_{\alpha_{i,j}}(t)=I_{n+1}+t E_{i,j}$, where $E_{i,j}$ is the elementary matrix whose entries $e_{k,l}$ satisfy 
  $e_{k,l}=\delta_{i k}\delta_{j l}$.
Thus, $\tau$ is simply the transposition map. 

\section{Polar Grassmannians}\label{section:polar grassmannians}
In this subsection $\Gamma$ is a polar $k$-Grassmannian of a building $\De$ of type $\sM_n$ over $\FF$, where $\sM_{n,k}(\FF)$ is as listed in Table~\ref{table:polar grassmannians}.
The building $\De$ is constructed from a non-degenerate reflexive sesquilinear or quadratic form $\zeta$ of Witt index $n$ on a vector space $V$ of dimension $m$ over the field $\FF$. 
The type of $\zeta$ is given in the table and $m$ is the subscript of the group, which is the full linear isometry group of $\zeta$.
In case $\zeta$ is $\sigma$-hermitian, we restrict to the case where $\sigma\in\Aut(\FF)$ has order $2$, $\FF$ is a quadratic extension over the fixed field $\FF^\sigma=\{x\in\FF\mid x^\sigma=x\}$, and the norm $N_\sigma\colon\FF\to\FF^\sigma$ is surjective.
\begin{table}[h]
$$
\begin{array}{clccc}
\sM_{n,k}(\FF)  & \mbox{ $\zeta$ }   & \mbox{group} & n & k \\
\hline
B_{n,k}(\FF) & \mbox{ parabolic orthogonal }               & \Spin_{2n+1}(\FF) &\ge 2& 1\le k <n\\
C_{n,k}(\FF) &\mbox{ symplectic }  & \Sp_{2n}(\FF) & \ge 2 & 1\le k\le n \\
D_{n,k}(\FF)  &\mbox{ hyperbolic orthogonal }  & \Om_{2n}^+(\FF) & \ge 3 & 1\le k\le n-2 \\
\wA_{2n,k}(\FF)  &\mbox{ $\sigma$-hermitian }&  U_{2n+1}(\FF) & \ge 2 & 1\le k < n \\
\wA_{2n-1,k}(\FF) &\mbox{ $\sigma$-hermitian } & U_{2n}(\FF)& \ge 2 & 1\le k < n\\
\wD_{n+1,k}(\FF)  & \mbox{ elliptic orthogonal }& \SO_{2n+2}^-(\FF)& \ge 2 & 1\le k< n \\
\hline
\end{array}$$
\caption{Polar Grassmannians}\label{table:polar grassmannians}
\end{table}

\mn
We first present a way to see that the Grassmann embeddings for these polar Grassmannians are polarized. Then we shall analyze $\beta$ and $\tau$ for the untwisted cases ($B_n$, $C_n$, and $D_n$).

\mn 
The points and lines of $\Gamma$ are also points and lines of the projective $k$-Grassmannian 
 $\bar{\Gamma}$ of type $A_{m-1,k}(\FF)$ associated to $V$ (See Section~\ref{section:projective grassmannians}).
The Grassmann embedding $\bar{e_\gr}$ of $\bar{\Gamma}$ restricts to a full projective embedding
  $e_\gr$ of $\Gamma$ into some subspace $V_\gr$ of the exterior power $\bigwedge^k V$.
This is called the {\em Grassmann embedding} of $\Gamma$.

\bpr\label{prop:grassmann embeddings are polarized}
Let $\Gamma$ be a polar $k$-Grassmannian as in Table~\ref{table:polar grassmannians}. Then the Grassmann embedding of $\Gamma$ is polarized.
\epr
\pf
Let $H$ be a singular hyperplane of $\Gamma$.
Since $\opp_\eI$ is the identity on $\eI$, we have $\Gamma=\Gamma^*$
 and so $H=H(p^*)$ consists of all points $q$ of $\Gamma$ not opposite to some point $p^*$, which also belongs to $\Gamma$. 
Viewing points of $\Gamma$ as $k$-spaces in $V$, this means that 
$H$ consists of all points $q$ of $\Gamma$ such that $q\cap (p^*)^\perp\ne 0$. Here $\perp$ denotes orthogonality with respect to $\zeta$.

Keeping in mind that $\zeta$ is non-degenerate we find that
 $\bar{p}^*=(p^*)^\perp$ is an $(m-k)$-space of $V$, that is a dual point in $(\bar{\Gamma})^*$.
Let $\bar{H}=H(\bar{p}^*)$ be the singular hyperplane of $\bar{\Gamma}$
 defined by $\bar{p}^*$.
Then $\bar{H}$ consists of {\em all} $k$-spaces $\bar{q}$ of $V$ with $\bar{q}\cap \bar{p}^*\ne 0$.
Thus, $\langle H\rangle_{e_\gr}\le \langle \bar{H}\rangle_{\bar{e_\gr}}$.
Since the Grassmann embedding $\bar{e_\gr}$ of $\bar{\Gamma}$ is polarized, 
 $\langle \bar{H}\rangle_{\bar{e_\gr}}$ is a hyperplane of $\bigwedge^k V$ that induces $\bar{H}$.
Since $\zeta$ is non-degenerate, there is a point $p$ in $\Gamma$
 opposite to $p^*$, that is, not contained in $(p^*)^\perp=\bar{p}^*$. 
Thus we find that the codomain $V_\gr$ of $\Gamma$ under $e_\gr$ is not entirely contained in the hyperplane $\langle \bar{H}\rangle_{\bar{e_\gr}}$.
Hence $\langle H\rangle_{e_\gr}$ is contained in the hyperplane $V_\gr\cap \langle \bar{H}\rangle_{\bar{e_\gr}}$ of $V_\gr$.
By Proposition~\ref{close=maxhyp} and Lemma~\ref{lem:induced hyperplanes}, $H$, is induced by $V_\gr$.
That is, $e_\gr$ is polarized.
\qed

\mn
Recall from Lemma~\ref{lem:form on exterior power} that  $\zeta$ induces a form $\zeta^\wedge$ on $\bigwedge^k V$ and hence on $V_\gr$ as follows. Namely, for $u_1,\ldots,u_k,v_1,\ldots,v_k\in V$, let 
 $$\zeta^{\wedge }(u_1\wedge \cdots\wedge u_k,v_1\wedge\cdots\wedge v_k)=\det(\zeta(u_i,v_j)).$$
Let $\perp^\wedge$ denote the orthogonality relation on $\bigwedge^k V$ with respect to $\zeta^\wedge$.
For any subspace $U\le \bigwedge^k V$, let 
 $\Rad(U,\zeta^\wedge)=U^{\perp^\wedge}\cap U$.
By Lemma~\ref{lem:G preserves zeta wedge}, since $G$ preserves $\zeta$, it also preserves $\zeta^\wedge$ and hence $\Rad(V_\gr,\zeta^\wedge)$.

\ble\label{lem:polar radical is zeta wedge radical}
Let $\Gamma$ be a polar $k$-Grassmannian as in Table~\ref{table:polar grassmannians}.
\begin{itemize}
\AI{a}
For any (dual) point $p^*\in \Gamma^*=\Gamma$ we have
 $\langle H(p^*)\rangle_{e_\gr}=e_\gr(p^*)^{\perp^\wedge}$;
\AI{b} as a consequence $\cR_{e_\gr}=\Rad(V_\gr,\zeta^\wedge)$.
\end{itemize}
\ele
\pf
(a) First note that, by Lemma~\ref{lem:form on exterior power}, $\zeta^\wedge$ is non-degenerate sesquilinear on $\bigwedge^k V$, so that $e_\gr(p^*)^{\perp^\wedge}$ is a proper hyperplane of $\bigwedge^k V$. 
On the other hand, by Proposition~\ref{prop:grassmann embeddings are polarized}, we know that $\langle H(p^*)\rangle_{e_\gr}$ is a proper hyperplane of $V_\gr$.
Thus, it suffices to prove that $\langle H(p^*)\rangle_{e_\gr}\sbe e_\gr(p^*)^{\perp^\wedge}$. 

Let $u=u_1\wedge\cdots\wedge u_k$ represent $p^*$
 and let $v=v_1\wedge\cdots\wedge v_k$ represent some point $q$ of $\Gamma$. 
 Now $\zeta^\wedge(u,v)= 0$ if and only if 
  the columns of the matrix $(\zeta(u_i,v_j))$ are linearly dependent, which happens if and only if $q\cap (p^*)^\perp\ne 0$ in $V$ and the latter is equivalent to saying that $q\in H(p^*)$.
In particular, $\langle H(p^*)\rangle_{e_\gr}\sbe e_\gr(p^*)^{\perp^\wedge}$. 
\qed

\mn
For the untwisted cases, we choose $\zeta$ and the basis $\cA$ for $V$ in the following way.
$$\begin{array}{cccc}
\sM_n\colon & B_n  & C_n & D_n\\[1ex]
\zeta \colon& 
\left(\begin{array}{@{}ccc@{}}
 2 & 0 & 0 \\
 0 & O_n & I_n\\
 0 & I_n & O_n\\
\end{array}\right)
&
\left(\begin{array}{@{}cc@{}}
 O_n & I_n\\
 -I_n & O_n\\
\end{array}\right)
&
\left(\begin{array}{@{}cc@{}}
 O_n & I_n\\
 I_n & O_n\\
\end{array}\right)
\\
&&&\\
\cA\colon & 
\{a_0,\ldots,a_{2n}\} & 
\{a_1,\ldots,a_{2n}\} &
\{a_1,\ldots,a_{2n}\}
\\
\end{array}
$$

\ble\label{lem:grassmann=weyl}
Let $\Gamma$ be a polar $k$-Grassmannian as in Table~\ref{table:polar grassmannians}, where $\sM_n$ is untwisted, that is, it is one of $B_n$, $C_n$ or $D_n$.
Then, the Grassmann embedding is the Weyl embedding.
\ele
\pf
The Grassmann embedding is the restriction of the Grassmann embedding for the projective $k$-Grassmannian into $\bigwedge^k V$. The codomain $V_\gr$ is by definition the subspace of $\bigwedge^k V$ spanned by the images of the points of $\Gamma$. Since we have $G\le \SL(V)$, the space $V_\gr$ is naturally a $G$-module. Transitivity of $G$ on the point set of $\Gamma$ shows that $V_\gr$ is the $G$-submodule of $\bigwedge^k V$ generated by $e_\gr(p)$ for any given point $p$ of $\Gamma$.

More precisely, under the Grassmann embedding the point $p=\langle a_1,\ldots,a_k\rangle$ of $\Gamma$ is sent to the $1$-space of $\bigwedge^k V$ spanned by $a_{1,\ldots,k}$.
Since $G\le \SL(V)$ in all cases, we see that the Grassmann  embedding is given by 
 $$\begin{array}{rl}
 e_\gr\colon \cP & \to \PP(\bigwedge^k V)\\
 g p &\mapsto \langle g a_{1,\ldots,k}\rangle, \mbox{ for any }g\in G. \\
 \end{array}$$
We now consider the Weyl embedding. 
First let $\FF=\CC$.
From~\cite{Car2005} we see that $v^+=a_{1,\ldots,k}\in \bigwedge^k V$ is a vector of highest weight $\lambda_k$. Now $V(\lambda)_\CC=V(\lambda)_\CC^0$ is the $G_\CC$-module generated by $v^+$.

Passing to an arbitrary field $\FF$, we see that $V(\lambda)_\FF$ is also a $G_\FF$-submodule of $\bigwedge^k V$ (where $V$ is now an $\FF$-vector space) that contains the $1$-dimensional highest weight space $V(\lambda)_{\FF,\lambda}=\FF v^+$.
The Weyl module $\bV$ is by definition the $G_\FF$-submodule of $V(\lambda)_\FF$ generated by $v^+$.
\qed

\bre
Note that we now have two proofs of the fact that the Grassmann embedding of a polar $k$-Grassmannian associated to a polar space of type $B_n$, $C_n$, or $D_n$ is polarized.
Namely, Proposition~\ref{prop:grassmann embeddings are polarized} gives a direct geometric proof, whereas another proof  comes from combining Lemma~\ref{lem:grassmann=weyl} and Proposition~\ref{prop:Weyl embedding is polarized}.
\ere

\mn
We now describe $\tau$ in terms of the description, given in~\cite{Car2005}, of the Lie algebra $\mfg_\CC$ and $G$ in its action on the natural module.
From the description of the root spaces in the Lie algebras of type $C_n$, and $D_n$ in loc.~cit., we see that $\tau$ is given by the transpose map on $\mfg_\CC$, and hence also on $G$. In the $B_n$ case $\tau$ is given by 
 $g\mapsto h^{-1} g^t h$, where $t$ denotes transpose and
  $h$ is the diagonal matrix $\diag\{2,1,\cdots,1\}$. 
This formula, which is initially computed for $\Char(\FF)=0$,  remains valid when $\Char(\FF)\ne 2$.

\mn
\bre
The forms $\zeta^\wedge$ and $\beta$ are not equal, but it follows from Lemmas~\ref{lem:polar radical is zeta wedge radical}~and~\ref{lem:grassmann=weyl} that for each $p\in \Gamma$
 there is a $p^*\in \Gamma^*=\Gamma$ such that
  $e_\gr(p^*)^{\perp^\wedge}=e_\gr(p)^{\perp^\beta}$, where
   $\perp^\beta$ denotes orthogonality with respect to the contravariant form $\beta$. In particular this means that 
    $\Rad(V_\gr,\beta)=\Rad(V_\gr,\zeta^\wedge)$.
\ere

\section{Projective Flag-Grassmannians}\label{section:projective flag grassmannians}
We continue the setup from Section~\ref{section:projective grassmannians}
 except that now $\Gamma$ is a $\eK$-shadow space of $\De$, for some 
  arbitrary non-empty subset $\eK\sbe \eI=\{1,2,\ldots,n\}$.
By Theorem~\ref{thm:Weyl embedding}, $\Gamma$ embeds into the Weyl module $\bV=V(\lambda_\eK)_\FF^0$ for $\mfg_\FF=\mfsl_{n+1}(\FF)$, which is the submodule of $V(\lambda_\eK)_\FF$  generated by the highest weight vector $v^+$. In this section, we construct the Weyl module from the natural $\mfsl_{n+1}(\FF)$-module $V$, that is, of highest weight $\lambda_1$. We give two constructions and study the form $\beta$.

\mn
Write $\lambda=\lambda_\eK$.
Our construction of $V(\lambda)_\FF$ follows~\cite{Fu1997}.
Starting from a $\DD$-module $V$, where $\DD$ is some integral domain, we shall give two descriptions of a space $V^\lambda_\DD$,  with the property that $V^\lambda_\FF=V(\lambda)_\FF$, when $\FF$ is a field (see Theorem~\ref{thm:Schur = Weyl}). To unburden our notation we shall drop the $\DD$ unless strictly needed.
The condition that $\DD$ be an integral domain is not necessary in all that follows (see loc. cit.), but it is all we need.

These constructions are valid for any weight $\lambda=\sum_{k\in \eK}l_k\lambda_k$, where $\lambda_k$ is the $k$th fundamental dominant weight. To this end we identify $\lambda$ with the Young diagram that has  $l_k$ {\em columns} of length $k$.  That is, the partition corresponding to the (rows of the) transpose of the diagram $\lambda$ is $\mu=(k^{l_k})_{k\in \eK}$, where $k$ runs through $\eK$ in decreasing order; we sometimes write $\mu=(\mu_1,\ldots,\mu_l)$, where $l=\sum_{k\in \eK}l_k$. In the constructions below we shall in fact assume $\eK\sbe [n+1]=\{1,2,\ldots,n+1\}$ so as to also include those modules $V^\lambda$ involving determinantal representations.

\myparagraph{A universal description of $V^\lambda$}
Assume that $\lambda$ has $d$ boxes.
Let $V^{\times \lambda}$ be a cartesian product of $d$ copies of $V$ indexed by 
the boxes of $\lambda$. 
Consider maps $\phi\colon V^{\times \lambda}\to W$, where $W$ is some $\DD$-module, with the following properties:
\begin{itemize}
\item[(i)] $\phi$ is $\DD$-linear in each argument;
\item[(ii)] $\phi$ is alternating in each column of $\lambda$;
\item[(iii)] for any $x\in V^{\times\lambda}$, $\phi(x)=\sum \phi(y)$, where
 $y$ runs through all vectors obtained from $x$ by an exchange between two given columns of $\lambda$ with a given set of boxes in the right chosen column.
\end{itemize} 
If $x$ and $y$ are vectors indexed by the boxes of the same Young diagram $\lambda$, then an {\em exchange} between $x$ and $y$ is determined by choosing, in $\lambda$, two distinct columns $i$ and $j$ along with a set of $s\ge 1$ boxes in both columns; $y$ is now obtained from $x$ by interchanging the coordinates of $x$ corresponding to the $s$ boxes in columns $i$ and $j$ of $\lambda$, while preserving their relative order within those columns.
For an illustration, see Figure~\ref{fig:exchange}.
For example, taking $\lambda=\lambda_3+\lambda_2$ and using rules (i), (ii), and (iii) with $i=1$, $j=2$, and $s=2$, for 
 $v_1,\ldots,v_5\in V$ and $\alpha\in\DD$, we require that
 $\phi(\alpha v_1,v_2,v_3,v_4,v_5)=
\alpha \phi(v_1,v_4,v_5,v_2,v_3)-
\alpha \phi(v_2,v_4,v_5,v_1,v_3)+
\alpha \phi(v_3,v_4,v_5,v_1,v_2)$.

\begin{figure}
$$\begin{array}{rr}
\begin{picture}(140,80)
\put(0,40){\makebox(10,0){
$x=(v_1,v_2,v_3,v_4,v_5)\leftrightarrow$}}
\multiput(70,30)(0,20){3}{\line(1,0){40}}
\multiput(70,10)(0,20){1}{\line(1,0){20}}
\multiput(70,10)(20,0){2}{\line(0,1){60}}
\multiput(110,30)(20,0){1}{\line(0,1){40}}
\put(80,20){\makebox(0,0){$v_1$}}
\put(80,40){\makebox(0,0){$v_2$}}
\put(80,60){\makebox(0,0){$v_3$}}
\put(100,40){\makebox(0,0){$v_4$}}
\put(100,60){\makebox(0,0){$v_5$}}
\end{picture}
&
\begin{picture}(140,80)
\multiput(10,30)(0,20){3}{\line(1,0){40}}
\multiput(10,10)(0,20){1}{\line(1,0){20}}
\multiput(10,10)(20,0){2}{\line(0,1){60}}
\multiput(50,30)(20,0){1}{\line(0,1){40}}
\put(20,20){\makebox(0,0){$v_4$}}
\put(20,40){\makebox(0,0){$v_2$}}
\put(20,60){\makebox(0,0){$v_5$}}
\put(40,40){\makebox(0,0){$v_1$}}
\put(40,60){\makebox(0,0){$v_3$}}
\put(100,40){\makebox(30,0){$\leftrightarrow (v_4,v_2,v_5,v_1,v_3)=y$}}
\end{picture}\\
\end{array}$$
\caption{An exchange between two sets of two boxes in columns $1$ and $2$}\label{fig:exchange}
\end{figure}
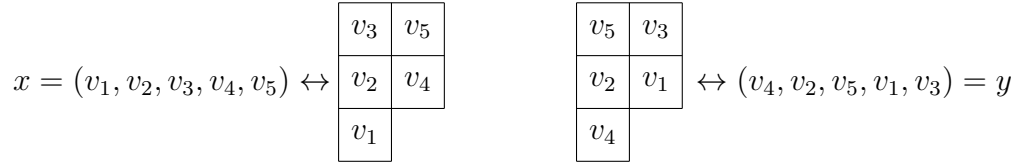We shall denote the universal target of such maps $\phi$ by $V^\lambda$. That is, there is a map $i\colon V^{\times \lambda}\into V^\lambda$ satisfying (i)-(iii) such that, given any map $\phi\colon V^{\times\lambda}\to W$ satisfying (i)-(iii), there is a $\DD$-linear map $\hat{\phi}\colon V^\lambda\to W$ with $\phi=\hat{\phi}\after i$.
This is called the {\em Schur module of shape $\lambda$}.
In the next paragraph we shall show the existence of $V^\lambda$. Moreover, in Theorem~\ref{thm:Schur = Weyl} we shall see that, for a field $\FF$, we have $V(\lambda)_\FF=V^\lambda_\FF$ if $\lambda$ has at most $n+1$ rows. In our situation this is satisfied always since $\eK\sbe [n+1]$.

\mn
\myparagraph{A concrete construction of $V^\lambda$}
All tensor and exterior products taken over boxes in $\lambda$ will be taken in order of the column word associated to $\lambda$; this word is obtained by concatenating the columns from left to right and 
 by ordering the boxes in each column from bottom to top\footnote{In~\cite{Fu1997} the columns of $\lambda$ are numbered top to bottom.}.

We shall refer to this construction as the {\em quotient construction of $V^\lambda$} since it realizes $V^\lambda$ as a quotient of $V^{\otimes\lambda}=V^{\otimes d}$ as follows.
To enforce rules (i) and (ii), let $V^{\wedge\lambda}=\otimes_{k\in \eK} (\bigwedge^k V)^{\otimes l_k}$ and consider the canonical map $\pi^\times_\wedge\colon V^{\times\lambda}\to V^{\wedge\lambda}$ (which is the composition of a $\DD$-$d$-linear map $\pi^{\times}_\otimes\colon V^{\times\lambda}\to V^{\otimes\lambda}$
 and a $\DD$-linear map $\pi^\otimes_\wedge\colon V^{\otimes\lambda}\to V^{\wedge\lambda}$). To enforce rule (iii), let $Q^\lambda(V)$ be the subspace of $V^{\wedge\lambda}$ generated by all vectors $\pi^\times_\wedge(x)-\sum\pi^\times_\wedge(y)$, where
  $y$ runs through all vectors obtained from $x$ by an exchange between two given columns of $\lambda$ with a given set of boxes in the right chosen column. Then $Q^\lambda(V)$ is the kernel of the  surjective $\DD$-linear map 
 $$\pi^\wedge_\lambda\colon V^{\wedge\lambda}\to  V^\lambda.$$
We shall denote by $\pi^\otimes_\lambda\colon V^{\otimes\lambda}\to V^\lambda$ the canonical map $\pi^\wedge_\lambda\after\pi^\otimes_\wedge$. 

\myparagraph{Properties of $V^\lambda$}
A filling of $\lambda$ from $[m]=\{1,2,\ldots,m\}$ 
 is a function $T\colon\lambda\to [m]$.
Given any  ordered set of vectors $\cB=\{b_1,\ldots,b_{n+1}\}$, and a filling $T$ of $\lambda$ from $[n+1]$,
 we get an element $b^\times_T$ of $V^{\times \lambda}$ by setting, for each box $d$ of $\lambda$, its entry indexed by box $d$ equal to $b_{T(d)}$.
Let 
$b_T^\otimes=\pi^\times_\otimes (b_T^\times)\in V^{\otimes\lambda}$, $b_T^\wedge=\pi^\otimes_\wedge(b_T^\otimes)\in V^{\wedge\lambda}$, and $b_T=\pi^\wedge_\lambda(b_T^\wedge)\in V^{\lambda}$.

\mn
We shall employ Lemma~\ref{lem:schur is Lie module} in the case where $V$ is the natural module for $\mfsl_{n+1}(\DD)$ and $\DD$ is $\ZZ$ or a field $\FF$.

\goodsofar

\ble\label{lem:schur is Lie module}
Suppose $V$ is a module for some Lie algebra $\mfg_\DD$ over $\DD$.
Then, the $\DD$-modules $V^{\otimes\lambda}$, $V^{\wedge\lambda}$, and $V^\lambda$ are $\mfg_\DD$-modules as well.
\ele
\pf
By equations~(\ref{eqn:Lie action on tensors}) and~(\ref{eqn:Lie action on exterior powers}) of Section~\ref{section:grassmannians} it follows that $V^{\otimes\lambda}$ and $V^{\wedge\lambda}$  are $\mfg_\DD$-modules. 
Now we show that $Q^\lambda(V)$ is a $\mfg_\DD$-submodule of $V^{\wedge \lambda}$.

Let $g\in \mfg_\DD$ and fix a subset $\cA=\{a_1,\ldots,a_{n+1}, a_{n+2}=ga_1,\ldots,a_{2n+2}=ga_{n+1}\}$ of $V$.
Given a filling $T\colon \lambda\to[n+1]$, and a box $d\in \lambda$, let $T^d$ be obtained from $T$ by setting, for all boxes $e$ of $\lambda$: 
$$T^d(e)=\left\{\begin{array}{ll} 
T(e) & \mbox{ if }e\ne d;\\
T(d)+n+1 &\mbox{ if }e=d.\\
\end{array}\right.$$
Thus, $a^\wedge_{T^d}$ is obtained from $a^\wedge_T$ by replacing the contents of box $d$ by its $g$-image, so  we have  $g a^\wedge_T=\sum_{d\in \lambda} 
a^\wedge_{T^d}$.
Now let $\xi$ be an exchange of $\lambda$, viewed as a permutation of the boxes in $\lambda$ and let $T'$ be the filling obtained by composing $\xi$ and $T$, that is $T'=T\after\xi\colon \lambda\to [n+1]$.
Given a filling $F\colon \lambda\to [n+1]$ we let $\xi\cdot a_F=a_{F\after \xi}$.
Then,
$$\begin{array}{lll}
 \xi g\cdot a^\wedge_T&=\xi \sum_{d\in\lambda}a^\wedge_{T^d} & =\sum_{d\in\lambda} a^\wedge_{T^d\after \xi}\\
g \xi\cdot a^\wedge_T&=g a^\wedge_{T\after \xi}&=\sum_{d'\in \lambda} a^\wedge_{(T\after \xi)^{d'}}\\
 \end{array}$$
It is straightforward to verify that $T^d\after \xi=(T\after \xi)^{\xi^{-1}(d)}$. Since $\xi$ is a permutation of $\lambda$ we can replace the last sum by
 $\sum_{\xi^{-1}(d)}a^\wedge_{(T\after \xi)^{\xi^{-1}(d)}}$ and conclude that
 $\xi g \cdot a^\wedge_T= g\xi \cdot a^\wedge_T$.
Thus, $g (a_T-\sum_\xi \xi a_T)=g a_T-\sum_\xi g\xi a_T=g a_T-\sum_\xi \xi g a_T$, where the sum is taken over all exchanges between two given columns of $\lambda$ with a given set of boxes in the right chosen column. 
Hence, $Q^\lambda(V)$ is a $\mfg_\DD$-module.
It follows that $V^\lambda=V^{\wedge\lambda}/Q^\lambda(V)$ is a $\mfg_\DD$-module as well.
\qed

\mn
From now on we shall assume that $V$ is free over $\DD$.

\mn
We shall now discuss bases. Given a basis $\cB$ for $V$, it is clear that $V^{\times\lambda}$ and $V^{\otimes\lambda}$ are free over $\DD$. In fact, 
$\cB^{\otimes\lambda}=\{b^\otimes_T\mid T\}$, where $T\colon \lambda\to [n+1]$ ranges over all possible fillings, is a $\DD$-basis for $V^{\otimes\lambda}$. Moreover, it is clear that $V^{\wedge\lambda}$ is free over $\DD$ with basis 
 $\cB^{\wedge\lambda}=\{b^\wedge_{T}\mid T\}$, where $T\colon\lambda\to [n+1]$ ranges over all fillings that are increasing down each column.
Next, we identify certain bases of $V^\lambda$.
A {\em Young tableau of shape $\lambda$} is a filling $T$ of $\lambda$ that is weakly increasing along each row of $\lambda$ and strictly increasing down each column of $\lambda$.
Let $\cB^\lambda$ be the collection of all $b_T$, where
 $T$ ranges over all Young tableaux of shape $\lambda$ with entries from $[n+1]$.
\ble\label{lem:V-lambda basis} {\rm~\cite[\S 8.1, Theorem 1]{Fu1997}}
Suppose $V$ is free over $\DD$ with basis $\cB=\{b_1,\ldots,b_{n+1}\}$, then $V^\lambda$ is free over $\DD$ with basis $\cB^\lambda$.
\ele

\ble\label{lem:V-lambda weights}
Let $V$ be a module for a Lie algebra $\mfg_\DD$ over $\DD$.
Suppose that $\cB=\{b_1,\ldots,b_{n+1}\}$ is a set of vectors such that $b_i$ has weight $\theta_i$ for each $i\in [n+1]$.
Then, for any filling $T\colon \lambda\to [n+1]$, the vector $b_T\in V^\lambda$ has weight $\sum_{i=1}^{m} t_i\theta_i$, where $t_i$ is the number of occurrences of $b_i$ in $b_T$
 and $m$ is the number of boxes in $\lambda$.
\ele
\pf
This is a straightforward calculation.
\qed

\mn
The constructions of $V^\lambda$ as well as Lemmas~\ref{lem:schur is Lie module},~\ref{lem:V-lambda basis}, and~\ref{lem:V-lambda weights}, show that if $V$ is a module for a (simple) Lie algebra $\mfg_\CC$ over $\CC$ and $A$ is an admissible lattice, then we can construct $A^{\times \lambda}$, $A^{\otimes\lambda}$,  $A^{\wedge\lambda}$, $Q^\lambda(A)$ and $A^\lambda$ replacing $V$ by $A$, taking $\DD=\ZZ$ and viewing $A$ as a module for the Lie algebra $\mfg_\DD=\mfU_\ZZ$ or $\mfg_\DD=\mfg(\lambda)_\ZZ$.
We summarize this and a little more in the following result.
 
\ble\label{lem:A lambda is admissible}
Let $V$ be a module for a simple Lie algebra $\mfg_\CC$ over $\CC$ and suppose $V_\FF=\FF\otimes_\ZZ A$ for some admissible lattice $A$.
\begin{itemize}
\AI{a} 
The $\ZZ$-modules $A^{\times \lambda}$, $A^{\otimes\lambda}$, $A^{\wedge\lambda}$, and $A^\lambda$ are 
 admissible lattices in $V^{\times \lambda}$, $V^{\otimes\lambda}$, $V^{\wedge\lambda}$, and $V^\lambda$ respectively.

\AI{b} For any field $\FF$ we have 
$(\FF\otimes_\ZZ A)^\lambda=\FF\otimes_\ZZ A^\lambda$.
In particular, if we construct $V(\lambda)_\FF=\FF\otimes_\ZZ A^\lambda$, then $V(\lambda)_\FF\cong V_\FF^\lambda$.
\AI{c} If $\cA=\{a_1,\ldots,a_{n+1}\}$ is a $\ZZ$-basis of weight vectors for $A$, then $\cA^\lambda$ is a $\ZZ$-basis of weight vectors 
for $A^\lambda$, where the weights are as described in Lemma~\ref{lem:V-lambda weights}. 
\end{itemize}
\ele
\pf
The admissible lattice $A$ is a free $\ZZ$-module as well as a 
 module for the universal enveloping algebra $\mfU_\ZZ$.
Thus, repeating the construction of $V^\lambda$ above replacing $V$ by $A$ and taking $\DD=\ZZ$, shows that 
$A^{\times\lambda}$, $A^{\otimes\lambda}$, $A^{\wedge\lambda}$ and $A^\lambda$ can be constructed.
Moreover, taking $\mfg_\DD=\mfU_\ZZ$, it follows from Lemma~\ref{lem:schur is Lie module} that these are $\mfU_\ZZ$-modules.

Let $\cA$ be the basis in (c). 
As we saw in the discussion preceding Lemma~\ref{lem:V-lambda basis}, $\cA$ canonically gives rise to bases 
$\cA^{\otimes\lambda}$ and $\cA^{\wedge\lambda}$ for 
$A^{\otimes \lambda}$, and 
  $A^{\wedge\lambda}$ respectively. 
By Lemma~\ref{lem:V-lambda basis} since $\cA$ is a basis for $V$,  $\cA^\lambda$ is a basis for $V^\lambda$ so in particular, $\cA^\lambda$ is independent over $\ZZ$ and of appropriate cardinality. The fact that $\cA^\lambda$ spans $A^\lambda$ over $\ZZ$ follows from the fact that
 $\cA^{\otimes\lambda}$ spans $A^{\otimes\lambda}$ and $A^\lambda$ is a quotient of $A^{\otimes\lambda}$.
Thus (a) follows. 
Part (b) is a special case of the remark preceding Lemma 1 in~\cite[\S 8.1]{Fu1997}.
Part (c) follows from Lemma~\ref{lem:V-lambda weights}. 
\qed

\mn
\bre\label{rem:v+as young tableau}
If in Lemma~\ref{lem:A lambda is admissible} we have $\mfg_\CC=\mfsl_{n+1}(\CC)$, then the natural module $V$ has a basis $\cA=\{a_1,\ldots,a_{n+1}\}$ whose weights satisfy
 $\theta_1\succ\theta_2\succ\cdots\succ \theta_{n+1}$ in the natural ordering on weights. 
 Namely, $\theta_i-\theta_{i+1}=\alpha_i$, the $i$-th fundamental root, for $i=1,\ldots,n$.
It then follows that the highest weight vector of $V^\lambda$ is $a_T$, where $T$ is the tableau whose $i$-th row is filled with $i$'s only.
Since $\theta_1+\cdots+\theta_k=\lambda_k$ for any $1\le k\le n$, $a_T$ has weight $\lambda$.
\ere

\bth\label{thm:Schur = Weyl}{\rm (Theorem 8.2 of \cite{Fu1997})}
If $\lambda$ has at most $n+1$ rows and $\FF=\CC$, then $V^\lambda$ is an irreducible representation of highest weight $\lambda$ for $\GL_{n+1}(\FF)$. These are all irreducible polynomial representations of $\GL_{n+1}(\FF)$.
\eth 

\bco
For any field $\FF$, the $\mfsl_{n+1}(\FF)$-modules
 $V(\lambda)_\FF$ and $V_\FF^\lambda$ are isomorphic.
\eco
\pf
Since $\lambda$ was constructed from the subset $\eK\sbe [n+1]$, all $V^\lambda_\CC$
 constructed in this section are irreducible $\GL_{n+1}(\CC)$-modules.
Since $\GL_{n+1}(\CC)$ is a central extension of $\SL_{n+1}(\CC)$ it follows that
 $V^\lambda_\CC$ is also irreducible as an $\SL_{n+1}(\CC)$-module.
Hence $V^\lambda_\CC$ is naturally an irreducible $\mfsl_{n+1}(\CC)$-module.
 
The well-known classification of finite dimensional irreducible modules for simple complex Lie algebras in particular ensures that, for each $\lambda$ (which in our case is integral and dominant), there is a unique irreducible $\mfsl_{n+1}(\CC)$-module of highest weight $\lambda$. Since $V_\CC^\lambda$ is an irreducible $\mfsl_{n+1}(\CC)$-module of highest weight $\lambda$, $L(\lambda)_\CC=V(\lambda)_\CC=V^\lambda_\CC$.

By Lemma~\ref{lem:A lambda is admissible} if $A$ is an admissible lattice in $V$, then $A^\lambda$ is an admissible lattice in $V^\lambda$.
The result follows, since by part (b) of that lemma, 
 $V(\lambda)_\FF=\FF\otimes_\ZZ A^\lambda=(\FF \otimes_\ZZ A)^\lambda=V_\FF^\lambda$, for any field $\FF$.
\qed

\mn
Note that by Lemma~\ref{lem:V-lambda basis}~(and~\ref{lem:A lambda is admissible}) the dimension of $V^\lambda$ is independent of the field $\FF$. This is a special case of Proposition~\ref{prop:U_ZZ-modules} parts (b) and (c).

As already noted in Remark~\ref{rem:not cyclic}, $V^\lambda$ is in general not cyclic.
\myparagraph{The Weyl embedding}
Now let $\lambda=\lambda_K$  be as in Subsection~\ref{subsec:Weyl module}, and $\FF$ some field.
Consider the following point of $\Gamma$:
 $p=(A_k)_{k\in \eK}$, where $A_k=\langle a_1,\ldots,a_k\rangle_V$.
Then by Theorem~\ref{thm:Weyl embedding} the Weyl embedding satisfies
$$\begin{array}{rl}
e_W\colon \cP & \to \PP((V^\lambda)^0)\\
 p & \mapsto a_T
\end{array}$$
where $T$ is the Young tableau whose $i$-th row is filled with $i$'s only.
That is $a_T=\pi^\wedge_\lambda(\bigotimes_{k\in \eK}\wedge_{i=1}^k a_i)$.
For any subset $\eK\sbe \eI$, let $e_\eK$ be the Weyl embedding of the 
 $\eK$-shadow space of $\Delta$.
Then we have
\beq\label{eqn:e_K factors}
e_\eK(p)=\pi^\wedge_\lambda(\bigotimes_{k\in \eK}e_{k}(A_k)).
\eeq
Now note that $\pi^\wedge_\lambda$ is a $G$-module homomorphism and that $G$ is transitive on the points of $\Gamma$.
Therefore we have equality~(\ref{eqn:e_K factors}) for any point $p=(A_k)_{k\in \eK}$ of $\Gamma$.

\goodsofar

\myparagraph{The form $\beta$ and a second construction of $V^\lambda$}
We first show that the contravariant form $\beta$ on $V^{\wedge\lambda}$ and $V^\lambda$ can be obtained from the contravariant form on $V$.
Recall that $\beta$ is constructed initially on the module $V(\lambda)_\CC$, then restricted to an admissible lattice, and then tensored with $\FF$ to get its equivalent on $V(\lambda)_\FF$.
Thus in order to describe $\beta$, we can and shall work over $\CC$. Our aim is to find a result similar to part (a) of Lemmas~\ref{lem:form on tensor product}~and~\ref{lem:form on exterior power}. To this end it is convenient to use the Schur functor (see e.g.~\cite[\S8.3]{Fu1997}~and~\cite{FuHa1991}).

\mn
We first use the Schur functor (and its dual) to describe 
 $V^\lambda$ and $(V^*)^\lambda$.
We shall consider $(V^*)^\lambda$ to be constructed from $V^*$ in the same way as $V^\lambda$ is constructed from $V$ and denote the maps corresponding to $\pi^\times_\lambda$, $\pi^\wedge_\lambda$, $\pi^\otimes_\lambda$ by putting a bar over them, so that we get $\bar{\pi}^\times_\lambda$, $\bar{\pi}^\wedge_\lambda$, and $\bar{\pi}^\otimes_\lambda$.
Fix $\lambda$ and let it have $d$ boxes.
Let $\bC=\CC[S_{d}]$, where $S_d$ is the symmetric group on 
 the set $[d]$.
In the sequel, for any right $\bC$-module $A$ and left $\bC$-module $B$, the tensor product
 $A\otimes_\bC B$ is the quotient of the usual tensor product $A\otimes_\CC B$ by the subspace generated by all expresssions
  $$(a\cdot \sigma)\otimes b-a\otimes (\sigma\cdot b)\mbox{ where } a\in A, \sigma\in S_d, b\in B.$$
We view  $V^{\otimes {d}}$ (resp.\ $(V^*)^{\otimes {d}}$) as a left (resp.~right) $\CC$-vector space  and a right (resp.~left) $\bC$-module, where $S_{d}$ naturally (resp.\ reverse) permutes the components of $V^{\otimes d}$ (resp.\ $(V^*)^{\otimes d}$). 

Then, we have natural isomorphisms 
$V^{\otimes {d}}\cong V^{\otimes {d}}\otimes_\bC \bC$
and 
$(V^*)^{\otimes {d}}\cong  \bC\otimes_\bC (V^*)^{\otimes {d}}$.
Note that the natural pairing $p\colon V\times V^*\to \CC$ given by 
 $\langle v,f\rangle = f(v)$ gives rise to the pairing
  $p^\otimes\colon V^{\otimes {d}}\times (V^*)^{\otimes {d}}\to \CC$
   by setting $\langle v_1\otimes \cdots\otimes v_{d},f_1\otimes \cdots \otimes f_{d}\rangle=\Pi_{i=1}^{d} f_i(v_i)$, thus identifying
  $(V^*)^{\otimes {d}}$ and $(V^{\otimes {d}})^*$.
Note here that we use the conventional notation $f(v)$ even though $V^*$ is a right vector space.
It follows from the above definitions that we now have
 \beq\label{eqn:dual C action}
 \langle v\otimes c\sigma,1\otimes f\rangle=\langle v\otimes 1,c\sigma\otimes f\rangle
 \eeq
for any $f\in (V^*)^{\otimes {d}}$,  
 $v\in V^{\otimes {d}}$, $c\in \CC$, and $\sigma\in S_{d}$.
Here the pairing $\langle\cdot,\cdot\rangle$ is induced by $p^\otimes$ via the natural isomorphisms $V^{\otimes {d}}\cong V^{\otimes {d}}\otimes_\bC \bC$
and 
$(V^*)^{\otimes {d}}\cong  \bC\otimes_\bC (V^*)^{\otimes {d}}$.

\mn
A {\em numbering} $U$ of $\lambda$ is a filling from $[{d}]$ without repeated entries.  
Let $U$ be the numbering of $\lambda$ that agrees with the natural ordering of boxes of $\lambda$ taken in the quotient construction of $V^\lambda$, that is, so that the column word of $U$ is $1,2,\ldots,{d}$. Let $R(U)$ and $C(U)$ be the {\em row group} and {\em column group} of $U$; that is, $R(U)$ (resp. $C(U)$) is the subgroup of $S_d$ that simultaneously preserves the subsets of numbers in $U$ associated with the rows (resp. columns) of $\lambda$.
Let 
$$\begin{array}{llll}
\rho_U=\sum_{\rho\in R(U)}\rho, &
\gamma_U=\sum_{\gamma\in C(U)}\sign(\gamma)\gamma, &
\sigma_U=\gamma_U\rho_U, 
&\bar{\sigma}_U=\rho_U\gamma_U.\\
\end{array}$$

The Specht module $S^\lambda$ with diagram $\lambda$ can be identified with
 the right $\bC$-module ${\sigma}_U \bC$ as well as the left $\bC$-module $\bC{\sigma}_U$
  (See e.g.~\cite[Ch. 4]{FuHa1991}).
Using the Schur functor, i.e.\ tensoring with the Specht module, we have isomorphisms
$$\begin{array}{@{}lrl}
&V^\lambda&\cong V^{\otimes {d}}\otimes_\bC \bC{\sigma}_U\\
&\pi^\otimes_\lambda(v_1\otimes\cdots\otimes v_{d})&\mapsto v_1\otimes\cdots\otimes v_{d}\otimes_\bC\sigma_U\\
\mbox{ and }&&\\
&(V^*)^\lambda&\cong {\sigma}_U\bC\otimes_\bC (V^*)^{\otimes {d}}\\
&\bar{\pi}^\otimes_\lambda(f_1\otimes\cdots\otimes f_{d})&\mapsto \frac{1}{k_\lambda}{\sigma}_U\otimes_\bC f_1\otimes\cdots\otimes f_{d},\\
\end{array}$$
where $k_\lambda$ is the integer such that $\sigma_U^2=k_\lambda\sigma_U$.
The quotient construction of $V^\lambda$ corresponds exactly to applying the Schur functor as above since by choice of $U$, we have
$$\begin{array}{lll}
V^{\otimes \lambda}=V^{\otimes d}\otimes_\bC \bC
&\stackrel{\pi^\otimes_\wedge}{\to} V^{\wedge \lambda}=V^{\otimes d}\otimes_\bC \bC\gamma_U &\stackrel{\pi^\wedge_\lambda}\to V^\lambda=V^{\otimes d}\otimes_\bC \bC\sigma_U\\
v_1\otimes\cdots\otimes v_d\otimes_\bC 1 & \mapsto 
v_1\otimes\cdots\otimes v_d\otimes_\bC \gamma_U & \mapsto
v_1\otimes\cdots\otimes v_d\otimes_\bC \gamma_U\rho_U\\
\end{array}$$
and $\pi^\otimes_\lambda=\pi^\wedge_\lambda\after\pi^\otimes_\wedge$ (cf.~\cite[\S 8.3]{Fu1997}).
In the construction of $(V^*)^\lambda$ one could use  $\bar{\sigma}_U$ instead of $\frac{1}{k_\lambda}\sigma_U$. Our choice is more natural, as we'll see in the proof of the following general lemma.

\ble\label{lem:V lambda*=V* lambda}
We have an isomorphism $(V^\lambda)^*\cong (V^*)^\lambda$ induced  by the pairing: $p^\lambda\colon V^\lambda\times (V^*)^\lambda\to \CC$ given by 
 $$\langle \pi^\otimes_\lambda(v_1\otimes\cdots\otimes v_{d}), 
 \bar{\pi}^\otimes_\lambda(f_1\otimes\cdots\otimes f_{d})\rangle=\sum_{\rho\in R(U)}\sum_{\gamma\in C(U)}\sign(\gamma)\Pi_{i=1}^{d} f_{i}(v_{(\gamma\rho)(i)}),$$ where
 $U$ is as above.
\ele

\pf
Consider the isomorphisms
$$\begin{array}{rl}
V^\lambda&\cong V^{\otimes {d}}\otimes_\bC \bC{\sigma}_U\\
(V^*)^\lambda&\cong {\sigma}_U\bC\otimes_\bC (V^*)^{\otimes {d}}\\
\end{array}$$
It can be shown that  $\sigma_U^2=k_\lambda\sigma_U$, for some non-zero $k_\lambda\in \ZZ$, so that, for $v\in V^{\otimes d}$ and $f\in (V^*)^{\otimes d}$, we have 
\beq\label{eqn:k_lambda}
\langle v\otimes_\bC\sigma_U,1\otimes_\bC f\rangle
=
\langle v\otimes_\bC \frac{1}{k_\lambda}{\sigma}_U^2,1\otimes_\bC f\rangle
\eeq 
It follows from~(\ref{eqn:dual C action}) that the pairing $p^\otimes$ 
 has the property that, 
\beq\label{eqn:sigma on pairing}
\langle v\otimes_\bC \frac{1}{k_\lambda}{\sigma}_U^2,1\otimes_\bC f\rangle
=
\langle v\otimes_\bC\sigma_U, \frac{1}{k_\lambda}{\sigma}_U\otimes_\bC f\rangle.
\eeq
We have maps:
$$
\begin{array}{lll}
V^{\otimes d} \otimes_\bC \bC\sigma_U
& \stackrel{i}{\into}  
V^{\otimes d} \otimes_\bC \bC
 & \stackrel{\pi}{\onto} 
V^{\otimes d} \otimes_\bC \bC\sigma_U
\end{array}
$$
where $i$ is the identity and $\pi$ is right multiplication by 
 $\sigma_U$.
The surjective map $\pi$ induces an injection
 $\pi^*\colon (V^{\otimes d}\otimes_\bC \bC\sigma_U)^*\into (V^{\otimes d}\otimes_\bC \bC )^*$.
Since $\sigma_U^2=k_\lambda\sigma_U$, we have
 $\pi\after i=k_\lambda\id$ so that for every $g\in  (V^{\otimes d}\otimes_\bC \bC\sigma_U)^*$ we have $g=\frac{1}{k_\lambda} i^*\after \pi^*(g)$.
Thus $i^*\colon   
(V^{\otimes d} \otimes_\bC \bC)^*\to (V^{\otimes d} \otimes_\bC \bC\sigma_U)^*$, given by $f\mapsto f\after i$, induces the inverse to $\frac{1}{k_\lambda}\pi^*$ on the image of $\pi^*$.
To see what this means, we identify $V^{\otimes d}\otimes_\bC \bC\sigma_U$ with its image under $i$.
Then $g=\frac{1}{k_\lambda} i^*\after \pi^*(g)$ says that  $g\in (V^{\otimes d}\otimes_\bC \bC\sigma_U)^*$
is simply the restriction of some element of $(V^{\otimes d}\otimes_\bC \bC)^*$
 (namely $\frac{1}{k_\lambda}\pi^*(g)$).
Now use the pairing $p^\otimes$ to view 
$\bC\otimes_\bC(V^*)^{\otimes d}\cong (V^{\otimes d}\otimes_\bC \bC)^*$. Namely define a map $b^\otimes$ by setting $c\otimes_\bC f\mapsto \langle \cdot,c\otimes_\bC f\rangle$, for all $c\in \bC$ and $f\in (V^*)^{\otimes d}$.
Then, from Equations~(\ref{eqn:k_lambda})~and~(\ref{eqn:sigma on pairing}) we can see that if this element $g$ is represented by 
 some $f\in \bC\otimes_\bC(V^*)^{\otimes d}$, then 
it is also represented by $\frac{1}{k_\lambda}{\sigma}_U\otimes f\in 
\sigma_U\bC\otimes_\bC (V^*)^{\otimes d}$.
Hence, 
 $i^*\after b^\otimes$ is an isomorphism between 
 ${\sigma}_U\bC\otimes_\bC (V^*)^{\otimes d}$ and 
 $(V^{\otimes d}\otimes_\bC \bC{\sigma}_U)^*$.
The definition of $\pi^\otimes_\lambda$, $\bar{\pi}^\otimes_\lambda$ and $b^\otimes$ imply that, for $v\in V^{\otimes d}$ and $f\in (V^*)^{\otimes d}$ we have
$\langle \pi^\otimes_\lambda(v),\bar{\pi}^\otimes_\lambda(f)\rangle
 = \langle v\otimes\sigma_U,1\otimes f \rangle$.
The conclusion of the lemma now follows from the meaning of $v_1\otimes\cdots\otimes v_d\otimes \sigma_U$ and the standard pairing $p^\otimes$.
\qed

\goodsofar

\mn
If $V$ is a module for a simple Lie algebra $\mfg_\CC$ over $\CC$,  and $\beta$ is a $\tau$-contravariant form on $V$, we shall denote the forms on $V^{\otimes\lambda}$ and $V^{\wedge\lambda}$ defined using Lemmas~\ref{lem:form on tensor product}~and~\ref{lem:form on exterior power} by $\beta^{\otimes\lambda}$ and $\beta^{\wedge\lambda}$.
It follows that $\beta^{\wedge\lambda}$ is $\tau$-contravariant.
We now extend this to a form on $V^\lambda$.

\ble\label{lem:form on V^lambda}
Let $V$ be a finite dimensional $\FF$-vector space, let 
 $\sigma$ be an automorphism of $\FF$ of order at most $2$ and let 
  $\zeta$ be a $\sigma$-sesquilinear form on $V$. Let $k\in \NN_{\ge 1}$.
 Moreover, let $\bC=\FF[S_{d}]$ and let $U$ and $\sigma_U$ be as above.
Then
\begin{itemize}
\AI{a} there is a unique $\sigma$-sesquilinear form $\zeta^{\lambda}$
 on $V^{\otimes {d}}\otimes_\bC \bC\sigma_U$ given by \\
  $$\begin{array}{l}
  \zeta^{\lambda }(v_1\otimes \cdots\otimes v_{d}\otimes_\bC\sigma_U, u_1\otimes\cdots\otimes u_{d}\otimes_\bC \sigma_U)\\
  =\sum_{\gamma\in C(U)}\sum_{\rho\in R(U)}\sign(\gamma)\sum_{\tilde{\gamma}\in C(U)}\sign(\tilde{\gamma})\Pi_{i=1}^{d} \zeta(v_{(\gamma\rho \tilde{\gamma})(i)},u_i);\\
  \end{array}$$
\AI{b}
if $\zeta$ is symmetric bilinear, so is $\zeta^{\lambda}$;
\AI{c}
if $V$ is a finite dimensional module for a simple Lie algebra $\mfg_\FF$ and $\zeta$ is $\tau$-contravariant, then so is $\zeta^{\lambda}$.
\end{itemize}
\ele
\pf
(a)
The form $\zeta$ can be obtained by composing
 the standard pairing $V\times V^*\to \FF$ given by $\langle v,f\rangle= f(v)$ with a
  $\sigma$-semilinear homomorphism $\phi\colon V\to V^*$ so that
  $\zeta(v,u)=\langle v,\phi(u)\rangle$.
Now the following composition that we shall  call $\phi^\lambda$ is again a $\sigma$-semilinear homomorphism; note that if $\phi$ is an isomorphism, then so is $\phi^\lambda$.
$$\begin{array}{rlcl}
 V^{\otimes {d}}\otimes_\bC \bC\sigma_U
 &\to &  V^{\otimes {d}}\otimes_\bC \bC\bar{\sigma}_U
 &\to {\sigma}_U\bC\otimes_\bC (V^*)^{\otimes {d}}
 \\

 u_1\otimes\cdots \otimes u_d\otimes_\bC \sigma_U
 &\to & u_1\otimes\cdots \otimes u_d\otimes_\bC\frac{1}{k_\lambda}\gamma_U\bar{\sigma}_U 
 &\to \frac{1}{k_\lambda}\sigma_U\gamma_U\otimes_\bC \phi(u_1)\otimes\cdots\otimes_\bC\phi(u_d)
 \\
\end{array}$$ 
Namely, the first map is the isomorphism given by right multiplication by $\frac{1}{k_\lambda}\gamma_U$, and whose inverse is given by multiplication on the right by $\rho_U$. 
The second map combines $\phi^{\otimes d}$ with the isomorphism $\sum_{x\in S_d}a_x x\mapsto \sum_{x\in S_d}a_x x^{-1}$ between the left and right regular representation of $\bC$. Note that $C(U)$ and $R(U)$ are closed under taking inverses so that 
$\gamma_U\rho_U\gamma_U$ is invariant under taking inverses.
Note here that 
 $$\phi^\lambda(u_1\otimes\cdots\otimes u_d\otimes_\bC\sigma_U)=\sum_{\tilde{\gamma}\in C(U)}\sign(\tilde{\gamma})\bar{\pi}^\lambda(\phi(u_{\tilde{\gamma}^{-1}(1)})\otimes\cdots\otimes\phi(u_{\tilde{\gamma}^{-1}(d)})).$$
We then set $\zeta^\lambda=p^\lambda(\cdot,\phi^\lambda(\cdot))$.
It is immediate from this construction that $\zeta^\lambda$ is  $\sigma$-sesquilinear. 
In order to obtain the formula we use that $\Pi_{i=1}^d \zeta(v_{(\gamma\rho)(i)},\phi(u_{\tilde{\gamma}^{-1}(i)}))=\Pi_{j=1}^d\zeta(v_{(\gamma\rho\tilde{\gamma})(j)},\phi(u_j))$. 

(b) Note that if $\zeta$ is symmetric, then $\zeta(v_{(\gamma\rho\tilde{\gamma})(i)},u_i)=\zeta(u_{(\tilde{\gamma}^{-1}\rho^{-1}\gamma^{-1})(j)},v_j)$ for $j=\gamma\rho\tilde{\gamma}(i)$.
Since the sums are taken over all $\gamma,\tilde{\gamma}\in C(U)$
 and $\rho\in R(U)$ and these are closed under taking inverses, we find that $\zeta^\lambda$ is symmetric as well.
(c) 
This follows from part (e) of Lemma~\ref{lem:form on tensor product}.
\qed

\ble\label{lem:zeta lambda at v+}
Suppose $\cA=\{a_1,\ldots,a_n\}$ is an orthonormal basis for $V$ with respect to $\zeta$ and 
 $T$ is the Young tableau of shape $\lambda$ whose $i$-th row is filled with $i$'s only.
Then, $\zeta^\lambda(a_T,a_T)=k_\lambda$. Moreover, 
 $\zeta^\lambda(a_T,a_{T'})=0$ for any Young tableau $T'\ne T$ of shape $\lambda$.
\ele
\pf
In view of Lemma~\ref{lem:form on V^lambda} we have
 $$
 \begin{array}{rl}
 \zeta^\lambda(a_T,a_T)
 &=
 \sum_{\gamma\in C(U)}\sum_{\rho\in R(U)}\sum_{\tilde{\gamma}\in C(U)}\sign(\gamma)\sign(\tilde{\gamma})\Pi_{i=1}^d\zeta(a_{(\gamma\rho\tilde{\gamma})(i)},a_i)\\
&= \sum_{\gamma\in C(U)}\sum_{\rho\in R(U)}\sum_{\tilde{\gamma}\in C(U)}\sign(\gamma)\sign(\tilde{\gamma}),\\
 \end{array}$$
where the sum is taken over those $\gamma,\rho,\tilde{\gamma}$ such that $\gamma\rho\tilde{\gamma}\tilde{\rho}=1$ for some 
 $\tilde{\rho}\in R(U)$. This is because $\zeta(a_{(\gamma\rho\tilde{\gamma})(i)},a_i)=1$
   if and only if $a_{(\gamma\rho\tilde{\gamma})(i)}$ and $a_i$ are in the same row and $0$ otherwise.
Now consider the equation $\gamma_U\rho_U\gamma_U\rho_U=k_\lambda\gamma_U\rho_U$. The coefficient of $1\in \bC$ on the right hand side is exactly the sum above and it clearly equals $k_\lambda$. 

Now if $T'\ne T$ is a Young tableau with the same Young diagram as $\lambda$ then for any permutation $\sigma\in S_d$, the fillings $\sigma T$ and $T'$ differ in some box. Since $\cA$ is orthonormal, this means that 
 $\Pi_{i=1}^d \zeta(a_{\sigma(i)},a_i)=0$. It follows that 
 $\zeta^\lambda(a_T,a_{T'})=0$.
\qed

\bco\label{cor:contravariant form on V^lambda}
Let $\beta$ be the symmetric bilinear $\tau$-contravariant form on $V$ as defined in Subsection~\ref{subsec:contravariant form}. Then $\frac{1}{k_\lambda}\beta^\lambda$ is the symmetric bilinear $\tau$-contravariant form on $V^\lambda$.
\eco 
\pf
This follows from Remark~\ref{rem:v+as young tableau}, Lemmas~\ref{lem:form on V^lambda}~and~\ref{lem:zeta lambda at v+}.
\qed

\subsection{An illustration}
We illustrate what happens in this section with a well-known example~\cite{ThaVan1999a}. Let $\De$ be the building of type $A_2$ and let $\Gamma$ be the $\{1,2\}$-shadow space of $\De$ over a field $\FF$.

We have $\mfg_\CC=\mfsl_3(\CC)$, which acts on its natural module $V$ by matrix-vector multiplication from the left, where
 vectors are coordinate vectors with respect to a basis $\cA=\{a_1,a_2,a_3\}$.
Writing $X_{i,j}$ for the elementary matrix whose non-zero entry is a  $1$ in the $(i,j)$ position, and setting $H_{i,j}=X_{i,i}-X_{j,j}$
 and $Y_{i,j}=X_{j,i}$ for $i<j$, we have a Chevalley basis 
 $$\sfC=\{X_{i,j}\mid i<j\}\cup \{H_{1,2},H_{2,3}\}\cup \{Y_{i,j}\mid i<j\}.$$
The multiplication is given by $[X_{i,j},X_{k,l}]=\delta_{j,k}X_{i,l}-\delta_{l,i}X_{k,j}$ for all $i,j,k,l$.
Now, in terms of Young diagrams, we have $\lambda=(2,1)$ and $V^\lambda_\CC$ is the adjoint representation, i.e.\ it is $\mfsl_3(\CC)$ itself under the action 
 $\ad{x}(y)=xy-yx$. The highest weight vector is $v^+=X_{1,3}$ and the minimal admissible lattice $A_{\min}$ is the $\ZZ$-span of the Chevalley basis $\sfC$.
 Writing $a_{i,j,k}=a_i\otimes a_j\otimes a_k\otimes\sigma_U$, and comparing the $\mfsl_3(\CC)$ action on itself with that on $V^\lambda_\CC$, we find that the Chevalley basis elements are identified with elements of  $V^{\otimes 3}\otimes_\bC\bC\sigma_U$ as follows:
$$\begin{array}{rlrlrl}
X_{1,3}&=a_{2,1,1} &  &  &  Y_{1,3}&=a_{3,2,3} \\
X_{1,2}&=-a_{3,1,1} &  H_{1,2}&=a_{2,1,3}-2a_{3,1,2} &  Y_{1,2}&=a_{3,2,2} \\
X_{2,3}&=a_{2,1,2} &  H_{2,3}&=a_{2,1,3}+a_{3,1,2} &  Y_{2,3}&=-a_{3,1,3} \\
\end{array}$$
We also have $H_{1,3}=2a_{2,1,3}-a_{3,1,2}$. 
If $\cA$ is orthonormal with respect to $\beta$, then by direct calculation from Lemma~\ref{lem:form on V^lambda} one verifies that $\beta^\lambda$ is given with respect to the Chevalley basis $\sfC$ by the following matrix
$$B=\left(\begin{array}{@{}cccrrccc@{}}
1&&&&&&&\\
&1&&&&&&\\
&&1&&&&&\\
&&&2&-1&&&\\
&&&-1&2&&&\\
&&&&&1&&\\
&&&&&&1&\\
&&&&&&&1\\
\end{array}\right).$$
Note that $B$ has full rank in all characteristics except $3$, where it has rank $7$. In that case $H_{1,2}-H_{2,3}$ spans the radical of $\beta^\lambda$. Note also that this form can be given as $\beta^\lambda(x,y)=\trace(xy^\tau)$, for all $x,y\in \mfsl_3(\CC)$ where $y^\tau$ denotes the transpose of $y$.

Now consider the geometric picture. Let $\FF$ be an arbitrary field. The shadow space $\Gamma$ is the geometry of point-line flags $(p,l)$ of $\PG(V_\FF)$. Identifying $V^\lambda_\FF$ with the adjoint module $\mfsl_3(\FF)$, the Weyl embedding is given by 
$$\begin{array}{rl}
\cP &\into \PG(V^\lambda_\FF)\\
(p,l)&\mapsto \langle w_{p,l}=v_p n_l^\tau\rangle,\\
\end{array}$$
where $v_p$ is some (column) vector spanning $p$ and $n_l$ is some (column) vector that is orthogonal to $l$ with respect to the inner product $\beta$ and $\tau$ denotes transposition.
Note that $\trace(w_{p,l})=0$ since $p$ lies on $l$. Considering the apartment $\Sigma$ of $\De$ given by the basis $\cA$, we find that 
$(\langle a_i\rangle,\langle a_i,a_k\rangle)\mapsto \langle a_k\otimes a_i\otimes a_i\otimes\sigma_U\rangle=\langle a_{k,i,i}\rangle=\langle X_{i,j}\rangle$, where $\{i,j,k\}=\{1,2,3\}$. In particular $(\langle a_1\rangle,\langle a_1,a_2\rangle)$ is sent to the space  $\langle a_{2,1,1}\rangle=\langle v^+\rangle$ of highest weight $\lambda=\lambda_1+\lambda_2$.

We now look at singular hyperplanes.
Two point-line flags $(p_1,l_1)$ and $(p_2,l_2)$ are opposite if and only if $p_1\not\in l_2$ and $p_2\not\in l_1$.
Considering the matrices 
$w_{p_1,l_1}=v_{p_1}n_{l_1}^\tau$ and 
$w_{p_2,l_2}=v_{p_2}n_{l_2}^\tau$ we find that $(p_1,l_1)$ and $(p_2,l_2)$ are opposite if and only if 
 $\trace(w_{p_1,l_1} w_{p_2,l_2})\ne 0$.
The singular hyperplane $H(p,l)$ of $\Gamma$ of points not opposite to $(p,l)\in \Gamma^*=\Gamma$ is therefore induced by the hyperplane
 $\ker(\trace(w_{p,l}*\cdot))$ of $\PG(V^\lambda_\FF)$.
It follows therefore that the polar radical of $\PG(V^\lambda_\FF)$ is the radical of the symmetric bilinear form
 $\zeta(x,y)=\trace(xy)$. 
Since the transpose map $\tau$ simply sends the image of $\Gamma^*=\Gamma$ to that of $\Gamma$, we see once again that $\Rad(\zeta)=\Rad(\beta^\lambda)$. 

As for minimal polarized embeddings for $\Gamma$, we note that in characteristic $3$, the identity matrix $I_3=H_{1,2}-H_{2,3}$ satisfies $\trace(I_3w_{q,m})=0$ for all $(q,m)\in \Gamma$. It follows that in this case the minimal polarized embedding is $\PG(6,\FF)$. In all other characteristics there is no radical and the minimal polarized embedding is $\PG(7,\FF)$.

\myparagraph{Acknowledgements}
I would like to thank my friends and colleagues Antonio Pasini and Ilaria Cardinali for their comments on earlier versions of this paper as well as for many discussions on this topic during my visits to the University of Siena.
I would also like to thank two anonymous referees for their invaluable comments.

\bibliographystyle{alpha}

\end{document}